\documentclass{amsart}
\usepackage{amsmath} 
\usepackage{amssymb}
\usepackage{mathrsfs}

\newtheorem{theorem}{Theorem}[section] 
\newtheorem{claim}[theorem]{Claim}

\newtheorem{trt}[theorem]{The Representation Theorem}

\newtheorem{conclusion}[theorem]{Conclusion}
\newtheorem{observation}[theorem]{Observation}

\theoremstyle{definition}
\newtheorem{definition}[theorem]{Definition}

\newtheorem{fact}[theorem]{Fact}

\newtheorem{convention}[theorem]{Convention}
\newtheorem{discussion}[theorem]{Discussion}

\newtheorem{hypothesis}[theorem]{Hypothesis}

\theoremstyle{remark}
\newtheorem{remark}[theorem]{Remark}
\newtheorem{question}[theorem]{Question}
\newtheorem{notation}[theorem]{Notation}

\newcommand{\lf}{{\rm lf}}
\newcommand{\dbl}{{\rm dbl}}
\newcommand{\CH}{{\rm CH}}

\newcommand{\col}{{\rm col}}
\newcommand{\Ded}{{\rm Ded}}

\newcommand{\TV}{{\rm TV}}
\newcommand{\Lim}{{\rm Lim}}

\newcommand{\fnq}{{\rm fnq}}

\newcommand{\pp}{{\rm pp}}

\newcommand{\acc}{{\rm acc}}

\newcommand{\suc}{{\rm suc}}
\newcommand{\odd}{{\rm odd}}
\newcommand{\CF}{{\rm CF}}
\newcommand{\AAM}{{\rm AAM}}
\newcommand{\JEP}{{\rm JEP}}

\newcommand{\AMM}{{\rm AMM}}
\newcommand{\IND}{{\rm IND}}
\newcommand{\ZFC}{{\rm ZFC}}
\newcommand{\Sym}{{\rm Sym}}
\newcommand{\Symc}{{\rm Symc}}
\newcommand{\uf}{{\rm uf}}

\newcommand{\ortp}{{\rm ortp}}

\newcommand{\even}{{\rm even}}

\newcommand{\otp}{{\rm otp}}

\newcommand{\NF}{{\rm NF}}
\newcommand{\AM}{{\rm AM}}
\newcommand{\tp}{{\rm tp}}

\newcommand{\seb}{{\rm sb}}
\newcommand{\Seq}{{\rm Seq}}

\newcommand{\lex}{{\rm lex}}

\newcommand{\at}{{\rm at}}
\newcommand{\Ord}{{\rm Ord}}

\newcommand{\aut}{{\rm aut}}

\newcommand{\uniq}{{\rm uniq}}
\newcommand{\bs}{{\rm bs}}
\newcommand{\cm}{{\rm cm}}
\newcommand{\tr}{{\rm tr}}
\newcommand{\pr}{{\rm pr}}
\newcommand{\bd}{{\rm bd}}

\newcommand{\qf}{{\rm qf}}

\newcommand{\cg}{{\rm cg}}

\newcommand{\id}{{\rm id}}

\newcommand{\LST}{{\rm LST}}

\newcommand{\dom}{{\rm dom}}

\newcommand{\exlf}{{\rm exlf}}
\newcommand{\slf}{{\rm slf}}

\newcommand{\Dom}{{\rm Dom}}

\newcommand{\Rang}{{\rm Rang}}

\newcommand{\rest}{{\restriction}}

\newcommand{\wilog}{{\rm without loss of generality}}
\newcommand{\Wilog}{{\rm Without loss of generality}}

\newcommand{\then}{{\underline{then}}}
\newcommand{\when}{{\underline{when}}}

\newcommand{\oor}{{\underline{or}}}
\newcommand{\Then}{{\underline{Then}}}

\newcommand{\If}{{\underline{if}}}
\newcommand{\Iff}{{\underline{iff}}}
\newcommand{\mn}{{\medskip\noindent}}
\newcommand{\sn}{{\smallskip\noindent}}

\newcommand{\gB}{{\mathfrak B}}
\newcommand{\ga}{{\mathfrak a}}

\newcommand{\cE}{{\mathscr E}}
\newcommand{\varp}{{\varepsilon}}

\newcommand{\cH}{{\mathscr H}}

\newcommand{\cF}{{\mathscr F}}
\newcommand{\cG}{{\mathscr G}}

\newcommand{\bbL}{{\mathbb L}}

\newcommand{\cM}{{\mathscr M}}

\newcommand{\gs}{{\mathfrak s}}

\newcommand{\cP}{{\mathscr P}}

\newcommand{\gS}{{\mathfrak S}}
\newcommand{\gk}{{\mathfrak k}}

\newcommand{\bbZ}{{\mathbb Z}}

\newcommand{\cS}{{\mathscr S}}
\newcommand{\cT}{{\mathscr T}}

\newcommand{\gt}{{\mathfrak t}} 
\newcommand{\cU}{{\mathscr U}}
\newcommand{\cX}{{\mathscr X}}

\newcommand{\cf}{{\rm cf}}

\newcount\skewfactor
\def\mathunderaccent#1#2 {\let\theaccent#1\skewfactor#2
\mathpalette\putaccentunder}
\def\putaccentunder#1#2{\oalign{$#1#2$\crcr\hidewidth
\vbox to.2ex{\hbox{$#1\skew\skewfactor\theaccent{}$}\vss}\hidewidth}}

\newbox\noforkbox \newdimen\forklinewidth
\setbox0\hbox{$\textstyle\bigcup$}
\setbox1\hbox to \wd0{\hfil\vrule width \forklinewidth depth \dp0
                        height \ht0 \hfil}
\setbox\noforkbox\hbox{\box1\box0\relax}
\def\unionstick{\mathop{\copy\noforkbox}\limits}
\def\nonfork#1#2_#3{#1\unionstick_{\textstyle #3}#2}
\def\nonforkin#1#2_#3^#4{#1\unionstick_{\textstyle #3}^{\textstyle
    #4}#2}
\setbox0\hbox{$\textstyle\bigcup$}
\setbox1\hbox to \wd0{\hfil{\sl /\/}\hfil}
\setbox2\hbox to \wd0{\hfil\vrule height \ht0 depth \dp0 width
                                \forklinewidth\hfil}
\newbox\doesforkbox
\setbox\doesforkbox\hbox{\box1\box0\relax}
\def\nunionstick{\mathop{\copy\doesforkbox}\limits}

\def\fork#1#2_#3{#1\nunionstick_{\textstyle #3}#2}
\def\forkin#1#2_#3^#4{#1\nunionstick_{\textstyle #3}^{\textstyle
    #4}#2}

\newcommand{\stickT}{%
\setbox255=\hbox{\raise1ex\hbox{$\hspace{0.2pt}\,\bullet\,$}}
\mathord{\rlap{\hbox to\wd255{\hss\hbox{$|$}\hss}}
\box255}
}
\newcommand{\stickS}{%
\setbox255=\hbox{\raise0.6ex\hbox{$\scriptstyle\bullet$}}
\mathord{\rlap{\hbox to\wd255{\hss\hbox{$\scriptstyle|$}\hss}}
\box255}
}

\newenvironment{PROOF}[2][\proofname.]
   {\begin{proof}[#1]\renewcommand{\qedsymbol}{\textsquare$_{\rm #2}$}}
   {\end{proof}}

\begin{document}

\title {lf groups, aec amalgamation, few automorphisms}
\author {Saharon Shelah}
\address{Einstein Institute of Mathematics\\
Edmond J. Safra Campus, Givat Ram\\
The Hebrew University of Jerusalem\\
Jerusalem, 9190401, Israel\\
 and \\
 Department of Mathematics\\
 Hill Center - Busch Campus \\ 
 Rutgers, The State University of New Jersey \\
 110 Frelinghuysen Road \\
 Piscataway, NJ 08854-8019 USA}
\email{shelah@math.huji.ac.il}
\urladdr{http://shelah.logic.at}
\thanks{The author thanks Alice Leonhardt for the beautiful typing.
  First typed February 18, 2016.  In References
  \cite[0.22=Lz19]{Sh:312} means \cite[0.22]{Sh:312} has label z19
  there, L stands for label; so will help if \cite{Sh:312} will
  change. Paper 1098}







\date {October 31, 2017}

\begin{abstract}
In \S1 we deal with amalgamation bases, e.g. we define when an a.e.c. $\gk$ has
$(\lambda,\kappa)$-amalgamation which means ``many" $M \in K^{\gk}_\lambda$ are
amalgamation bases.  We then consider what happens for the class of lf groups. 
In \S2 we deal with weak
definability of $a \in N \backslash M$ over $M$, for $\bold
K_{\exlf}$.  In \S3 we deal with indecomposable members of
$\bold K_{\exlf}$ and with the existence of universal members of
$K^{\gk}_\mu$, for $\mu$ strong limit of cofinality $\aleph_0$.  Most
note worthy: if $\bold K_{\lf}$ has a universal model in $\mu$ then it
has a canonical one similar to the special models, (the parallel to
saturated ones in this cardinality).  
In \S4 we prove ``every $G \in \bold K^{\lf}_{\le \lambda}$ 
can be extended to a complete $(\lambda,\theta)$-full $G$" for
many cardinals.  In a continuation we 
may consider ``all the cardinals" or at least ``almost all the
cardinals"; also, we may consider a priori fixing the outer automorphism group.
\end{abstract}

\maketitle
\numberwithin{equation}{section}
\setcounter{section}{-1}
\newpage

\centerline {Anotated Content}
\bigskip

\noindent
\S0 \quad Introduction, (label w), pg.\pageref{0}
\bigskip

\noindent
\S1 \quad Amalgamation Basis, (label a), pg.\pageref{1}
\mn
\begin{enumerate}
\item[${{}}$]  [Consider an a.e.c. $\gk$, e.g. $\bold K_{\lf}$. We
  define $\AM_{\gk} = \{(\lambda,\kappa):\lambda \ge \kappa =
  \cf(\kappa),\lambda \ge \LST_{\gk}$ and the $\kappa$-majority
of $M \in K^{\gk}_\lambda$ are amalgamation bases$\}$, on
  ``$\kappa$-majority" see below.  What pairs have to be there?  I.e.
  for all a.e.c. $\gk$ with $\LST_{\gk} \le \theta$.  One
  case is when $M \in K^{\gk}_\lambda$ is $(< \kappa)$-existentially
  closed, $\kappa \in (\LST_{\gk},\lambda]$ is a compact cardinal or
  just satisfies what is needed for $M$.  This implies
  $(\lambda,\kappa) \in \AM_{\gk}$.  A similar argument gives ``$\kappa$ weakly
  compact $> \LST_{\gk} \Rightarrow (\kappa,\kappa) \in \AM_{\gk}$".
  Those results are naturally expected but surprisingly there are 
considerably more cases: if $\lambda$ is
  strong limit singular of cofinality $\kappa$ and $\kappa$ is a measurable
  cardinal $> \LST_{\gk}$ then $(\lambda,\kappa) \in \AM_{\gk}$.
  Moreover if also $\theta \in (\LST_{\gk},\lambda]$ is a measurable cardinal
    then $(\lambda,\theta) \in \AM_{\gk}$.]
\end{enumerate}
\bigskip

\noindent
\S2 \quad Definability, (label n), pg.\pageref{2}
\mn
\begin{enumerate}
\item[${{}}$]  [For an a.e.c. $\gk$, we may say $b_1$ is
  $\gk$-definable in $N$ over $M$ when $M \le_{\gk} N,b_1 \in N
  \backslash M$ and for no $N_*,b_1,b_2$ do we have $M \le_{\gk} N_*,b_1
  \ne b_2 \in N_*$ and $\ortp(b_\ell,N,N_*) = \ortp(b,M,N)$; there are
  other variants.  We clarify the situation for $\bold K_{\lf}$.]
\end{enumerate}
\bigskip

\noindent
\S3 \quad  $\bold K_{\lf}$, (label b), pg.\pageref{3}

\S(3A) \quad Indecomposability, pg.\pageref{3a}
\mn
\begin{enumerate}
\item[${{}}$]  [We say $M \in \bold K_{\gk}$ is
  $\theta$-indecomposable (or $\theta \in \CF(M)$) 
\when \, there is no strictly
  $<_{\gk}$-increasing sequence $\langle M_i:i < \theta \rangle$ with
  union $M$.  We clarify the situation for $\bold K_{\lf}$.]
\end{enumerate}

\S(3B) \quad Universality, pg.\pageref{3b}
\mn
\begin{enumerate}
\item[${{}}$]  [Let $\mu$ be strong limit of cofinality $\aleph_0$.
  We characterize when there is a universal member of $\bold K^{\lf}_\mu$
and assuming this, prove the existence of a substitute for ``special
model in $\mu$".  This works for any suitable a.e.c.]
\end{enumerate}
\bigskip

\noindent
\S4 \quad Complete $H$ are dense in $\bold K^{\exlf}_\lambda$ for some
$\lambda'\, s$, (label c), pg.\pageref{4}
\mn
\begin{enumerate}
\item[${{}}$]  [Our aim is to find out when for $\mu = \lambda$ (or
  just $\mu \le \lambda$) every
  $G \in \bold K^{\lf}_\mu$ can be extended to a complete $H \in
\bold K^{\exlf}_\lambda$, moreover $(\lambda,\sigma)$-full one.  We prove
  this for many $\lambda$'s.]
\end{enumerate}
\bigskip

\noindent
\S5 \quad More Uncountable Cardinals, (label e), pg.\pageref{5}
\mn
\begin{enumerate}
\item[${{}}$]  [In \S(5A) we outline a new proof for all $\lambda >
  \beth_\omega$ or just $\lambda = \lambda^{\langle
    \theta;\aleph_0\rangle},\theta = \cf(\theta) \in
  [\aleph_1,\lambda]$.  In it we consider an increasing continuous
  sequence $\langle G_{\bold m,i}:i \le \theta\rangle$ of groups; we
  try to get: if $\pi$ is an automorphism of $G_{\bold m,\theta}$ then
  for some club $E$ of $\theta$ we have $i \in E \Rightarrow 
\pi(G_{\bold m,i}) = G_{\bold m,i}$.  For this we use more group theory.  Still
 to deal with $\pi \rest G_{\bold m,i}$ for 
$i \in S_{\gs} \cap E$ we should use
  $\Pr_*(\lambda,\lambda,\aleph_0,\aleph_0)$, colouring principle.  In \S(5B)
we try to prove the theorem for every $\lambda > \mu$
  (see \ref{b16}) compared to \S(4A):
\newline
\underline{in \ref{c20}}:  weak $\theta = \lambda^+$(?)
\newline
\underline{in \ref{c23}}:  use $(S_1,S_3)$ only(?)
\newline
In the theorem we concentrate on the $\lambda = \mu^+$ case; works for
$\lambda > \mu$.]
\end{enumerate}

\S(5A) \quad Regular $\lambda$, pg.\pageref{5A}
\mn
\begin{enumerate}
\item[${{}}$]  [We prove that for ``most" regular $\lambda$, the
  complete groups are dense in $\bold K^{\exlf}_\lambda$.  We also
  give a more detailed proof of \ref{c29}.]
\end{enumerate}

\S(5B) \quad Singular $\lambda$, pg.\pageref{5B}
\mn
\begin{enumerate}
\item[${{}}$]  [We deal with non-strong limit ones.]
\end{enumerate}

\S(5C) \quad Strong Limit Singular and a New Property, pg.\pageref{5C}
\mn
\begin{enumerate}
\item[${{}}$]  [We do this by continuing \cite{Sh:331}.]
\end{enumerate}

\S(5D) \quad Continuation, pg.\pageref{5D}
\mn
\begin{enumerate}
\item[${{}}$]  [In \ref{n14} we try to weaken the construction for
  $\theta$-indecomposability using two rows of $k$.  Now \ref{n23},
  \ref{n72} are relatives of \ref{i8}.  Next Definition \ref{e3},
  \ref{e31} are relatives of $\Pr_*$ and Claim \ref{e7}, \ref{e16},
  \ref{e22}, \ref{e40} are relatives of \ref{i8} inside \ref{e13},
  successor of singulars with \ref{e22} on singular not strong limit;
  but unlike $\Pr_*$ the $u_\varp$ are of cardinality $< \partial$
  instead of finite and Definition \ref{e10}, define $\bold M_{2.1}$ 
and $\Pr_{2.7}$ relative of $\Pr_*$.  Now Claim \ref{e34} a try to draw
complete 5.  In the proof of \ref{e13} use ``$|u_\varp| < \partial$"
appear.  Question: in what point needed?]
\end{enumerate}
\newpage

\section {Introduction} \label{0}
\bigskip

\subsection {Review}\
\bigskip

We deal mainly with the class of locally finite groups.  We continue
\cite{Sh:312}, see history there.

\noindent
We wonder:
\begin{question}
\label{w6}
1) Is there a universal $G \in \bold K^{\lf}_\lambda$, e.g. for $\lambda =
\aleph_1 < 2^{\aleph_0}$, i.e. consistently?

\noindent
2) Is there a universal $G \in \bold K^{\lf}_\lambda$, e.g. for $\lambda =
\beth_\omega$?  Or just $\lambda$ strong limit of cofinality $\aleph_0$
(which is not above a compact cardinal)?
\end{question}
\bigskip

See the end of \S(3B).  This leads to questions on the
existence of amalgamation bases.  We give general claims on existence of 
amalgamation bases in \S1.

That is, we ask:
\begin{question}
\label{w8}
For an a.e.c. $\gk$ or just a universal class (justified by \S(0C)) we
ask:

\noindent
1) For $\lambda \ge \LST_{\gk}$, are the amalgamation bases (in
$K^{\gk}_\lambda$) dense in $K^{\gk}_\lambda$?  (Amalgamation basis in
$\gk_\lambda$, of course, see \ref{w35}).

\noindent
2) For $\lambda \ge \LST_{\gk}$ and $\kappa = \cf(\kappa)$ is the
$\kappa$-majority of $M \in K^{\gk}_\lambda$ amalgamation bases? (On
$\kappa$-majority, see \ref{w8}).  The set of such $(\lambda,\kappa)$ is
called $\AM_{\gk}$. For weakly compact $\lambda > \LST_{\gk}$.
\end{question}
\bigskip

Using versions of existentially closed models in $K^{\gk}_\lambda$,
for $\lambda$ weakly compact we get
$(\lambda,\lambda) \in \AM_{\gk}$.  But surprisingly there are other
cases: $(\lambda,\kappa)$ when $\lambda$ is strong limit singular,
with $\cf(\lambda) > \LST_{\gk}$ measurable and $\kappa =
\cf(\lambda)$ or just $\lambda > \kappa > \LST_{\gk},\kappa$ is
measurable.

This is the the content of \S1.

In \S(3B) we return to the universality problem for 
$\mu = \beth_\omega$ or just strong limit of
cofinality $\aleph_0$.  We prove for $\bold K_{\lf}$ and similar
classes that if there is a universal model of cardinality $\mu$, \then
\, there is something like a special model of cardinality $\mu$, in
particular, universal and unique up to isomorphism.  This relies on
\S(3A), which proves the existence of so-called
$\theta$-indecomposable (i.e. $\theta$ is not a possible cofinality)
 models in $\bold K_{\lf}$ of various cardinalities.

Returning to Question \ref{w6}(2), a possible avenue is to try to prove 
the existence of universal members in $\mu$ when $\mu =
\Sigma \mu_n$ each $\mu_n$ measurable $< \mu$, i.e. maybe for some
reasonable classes this holds.

In \S2 we deal with the number of $a \in G_2$ definable over $G_1
\subseteq G_2$ in the orbital sense and find a ZFC bound.

\noindent
We consider in \S4:
\begin{question}
\label{w2}
For each pair $\lambda \ge \mu + \aleph_1$ or even 
$\lambda = \mu > \aleph_1$, does every 
$G \in \bold K^{\lf}_{\le \mu}$ have a complete extension
in $\bold K^{\exlf}_\lambda$.

The earlier results assume more than $\lambda > \mu$,
e.g. $\lambda = \mu^+ \wedge \mu^{\aleph_0} = \mu$ and $(\lambda,\mu)
= (\aleph_1,\aleph_0)$.  Note that for $\bold K_{\lf}$, the statement is
stronger when fixing $\lambda$ we 
increase $\mu$ (because every $G_1 \in \bold K^{\lf}_\mu$
has an extension in $\bold K^{\lf}_\lambda$ when $\lambda \ge \mu$).
We intend to deal in \S4,\S5 with proving it for most pairs
$\lambda \ge \mu + \aleph_1$, even when $\lambda = \mu$.  
Note that if $\lambda = \mu^+$ and we construct a sequence
$\langle G_i:i < \lambda\rangle$ of members 
from $\bold K_\mu$ increasing continuous,
$G_0 = G$ with union of cardinality $\lambda$ 
then any automorphism $\pi$ of $H = \bigcup\limits_{i <
  \lambda} G_i$ satisfies $\{\delta < \lambda:\pi$ maps $G_\delta$
onto $G_\delta\}$ is a club, this helps.  But as we like to have
$\lambda = \mu$ we can use only $\langle G_i:i < \theta\rangle$, with
$\theta = \cf(\theta) \le \lambda$, to be chosen appropriately.
However, we have a substitute: if for unboundedly many 
$i < \theta,\theta \notin \CF(G_i)$, i.e. $G_i$
is not $\theta$-decomposable and $\theta = \cf(\theta) > \aleph_0$,
then for any automorphism $\pi$ of $G_\theta = \bigcup\limits_{i} G_i$
the set $E = \{\delta < \theta:\pi(G_\delta) = G_\delta\}$ is a club
of $\theta$.  See, e.g. Shelah-Thomas \cite[\S(3A)]{ShTh:524} on $\CF(G)$, the
cofinality spectrum of $G$.

An additional point is that we like our $H$ to be more existentially
closed, this is interpreted as being $(\lambda,\theta)$-full.  We also need to
``have a list of $\lambda$ countable subsets which is dense
enough", for this we use $\lambda =
\lambda^{\aleph_0}$ or $\lambda = \lambda^{\langle
  \theta;\aleph_0\rangle}$.

See a related work with Paolini which deals with universal lf
classes.  Where some of the cases here fit, we deal there with other cases.
\end{question}
\bigskip

\subsection {Amalgamation Spectrum} \label{0B}
\bigskip

On a.e.c. see \cite{Sh:88r}, \cite{Sh:E53}.
We note below that the versions of the amalgamation spectrum are the
same (fixing $\lambda \ge \kappa$) for:
\mn
\begin{enumerate}
\item[$(*)$]
\begin{enumerate}
\item[(a)]  a.e.c. $\gk$ with $\kappa = \LST_{\gk},\lambda = \kappa +
  (\tau_{\gk})$;
\sn
\item[(b)]  universal $\bold K$ with $\kappa = \sup\{\|N\|:N \in \bold
  K$ is f.g.$\},\lambda = \kappa + |\tau_{\gk}|$; 
\end{enumerate}
\end{enumerate}
\bigskip

\noindent
Recall
\begin{trt}
\label{w32}
Let $\lambda \ge \kappa \ge \aleph_0$.

\noindent
1) For every a.e.c. $\gk$ with $|\tau_{\gk}| \le \lambda$ and
$\LST_{\gk} \le \kappa$ there is $\bold K$ such that:
\mn
\begin{enumerate}
\item[(a)]  
\begin{enumerate}
\item[$(\alpha)$]  $\bold K$ is a universal class;
\sn
\item[$(\beta)$]  $|\tau_{\bold K}| \le \lambda,\tau_{\bold K}
  \supseteq \tau_{\gk},|\tau_{\bold K} \backslash \tau_{\gk}| \le \kappa$;
\sn
\item[$(\gamma)$]  any f.g. member of $\bold K$ has cardinality $\le
\kappa$.
\end{enumerate}
\item[(b)]  $K_{\gk} = \{N \rest \tau_{\gk}:N \in \bold K\}$, moreover:
\sn
\item[(b)$^+$]   if $(\alpha)$ and $(\beta)$, then $(\gamma)$, where:
\sn
\begin{enumerate}
\item[$(\alpha)$]  $I$ is a well founded partial order such that
  $s_1,s_2 \in I$ has a mlb (= maximal lower bound) called $s_1 \cap s_2$;
\sn
\item[$(\beta)$]  $\bar M = \langle M_s:s \in I\rangle$ satifies $s
  \le_I t \Rightarrow M_s \le_{\gk} M_t$ and $M_{s_1} \cap M_{s_2} =
  M_{s_1 \cap s_2}$; 
\sn
\item[$(\gamma)$]  there is $\bar N$ such that:
\begin{itemize}
\item  $\bar N = \langle N_s:s \in I\rangle$;
\sn
\item  $N_s \in \bold K$ expand $M_s$;
\sn
\item  $s \le_I t \Rightarrow N_s \subseteq N_t$.
\end{itemize}
\end{enumerate}
\sn
\item[(b)$^{++}$]  Moreover, in clause (b)$^+$, if $I_0 \subseteq I$
  is downward closed and $\bar N^0 = \langle N^0_s:s \in I_0\rangle$
  is as required in (b)$^+$ on $\bar N \rest I_0$, then we can demand
  there that $\bar N \rest I_0 = \bar N^0$.
\end{enumerate}
\end{trt}

\begin{PROOF}{\ref{w32}}
As in \cite{Sh:88r}.
\end{PROOF}

\begin{definition}
\label{a32f}
1) We say $\bold K$ is a universal class \when \,:
\mn
\begin{enumerate}
\item[(a)]  for some vocabulary $\tau,\bold K$ is a class of
  $\tau$-models;
\sn
\item[(b)]  $\bold K$ is closed under isomorphisms;
\sn
\item[(c)]  for a $\tau$-model $M,M \in \bold K$ iff every finitely
  generated submodel of $M$ belongs to $\bold K$.
\end{enumerate}
\end{definition}

\begin{claim}
\label{w35}
For $\gk,\bold K$ as in \ref{w32} and see Definition \ref{a31}.

\noindent
1) If $N \in \bold K_{\lambda_0},M = N \rest \tau_{\gk}$, \then \,:
$N$ is a $(\lambda_1,\lambda_2)$-amalgamation base in $\bold K$ \Iff
\, $M$ is a $(\lambda_1,\lambda_2)$-amalgamation base in $\gk$.

\noindent
2) $\bold K$ has $(\lambda_0,\lambda_1,\lambda_2)$-amalgamation \Iff
\, $\gk$ has $(\lambda_0,\lambda_1,\lambda_2)$-amalgamation.

\noindent
3) $\AM_{\bold K} = \AM_{\gk}$ see Definition \ref{a31}(5).
\end{claim}

\begin{observation}
\label{w38}
If $\bold K$ is a universal class, $\kappa \ge \sup\{\|N\|:N \in \bold
K$ is finitely generated$\},\lambda \ge \kappa + |\tau_{\bold K}|$,
\then \, $\gk = (\bold K,\subseteq)$ and $\bold K$ are as in the conclusion
of \ref{w32}.
\end{observation}
\bigskip

\subsection {Preliminaries} \label{0C}\
\bigskip

\noindent
Recall (this is used only in \S5 and see more there)
\begin{definition}
\label{w11}
Assume $\lambda \ge \mu \ge \sigma + \theta_0 + \theta_1,\bar\theta
   = (\theta_0,\theta_1)$; if $\theta_0 = \theta_1$ we may write
   $\theta_0$ instead of $\bar\theta$.

\noindent
1) Let $\Pr_0(\lambda,\mu,\sigma,\bar\theta)$ mean that there is $\bold
c:[\lambda]^2 \rightarrow \sigma$ witnessing it which means:
\mn
\begin{enumerate}
\item[$(*)_{\bold c}$]  if (a) then (b) where:
\sn
\begin{enumerate}
\item[$(a)$]
\begin{enumerate}
\item[$(\alpha)$]   for $\iota = 0,1$ and $\alpha <
  \lambda$ we have $\bar\zeta^\iota = \langle
  \zeta^\iota_{\alpha,i}:\alpha < \mu,i < \bold i_\iota\rangle$,  a
  sequence without repetitions of ordinals $< \lambda$  
\sn
\item[$(\beta)$]   $\bold i_0 < \theta_0,\bold i_1 < \theta_1$;
\sn
\item[$(\gamma)$]  $h:\bold i_0 \times \bold i_1 \rightarrow \sigma$
\end{enumerate}
\sn
\item[$(b)$]   for some $\alpha_0 < \alpha_1 < \mu$ we have:
\sn
\begin{itemize}
\item  if $i_0 < \bold i_0$ and $i_1 < \bold i_1$ then
$\bold c \{\zeta^0_{\alpha_0,i_0},\zeta^1_{\alpha_1,i_1}\} =
  h(i_0,i_1)$.
\end{itemize}
\end{enumerate}
\end{enumerate}
\mn
2) For $\iota \in \{0,1\}$ let
$\Pr_{0,\iota}(\lambda,\mu,\sigma,\bar\theta)$ be defined similarly
but we replace $(a)(\beta)$ and $(b)$ by:
\mn
\begin{enumerate}
\item[$(a)$] 
\begin{enumerate}
\item[$(\beta)'$]  $h:\bold i_\iota \rightarrow \sigma$; 
\end{enumerate}
\sn
\item[$(b)'$]   for some $\alpha_0 < \alpha_1 < \mu$ we have
\begin{itemize}
\item  if $i_0 < \bold i_0$ and $i_1 < \bold i_1$ then
$\bold c \{\zeta^0_{\alpha_0,i_0},\zeta^1_{\alpha_1,i_1}\} = h(i_\iota)$.
\end{itemize}
\end{enumerate}
\mn
3) Let $\Pr^{\uf}_{0,\iota}(\lambda,\mu,\sigma,\bar\theta)$ mean that some
$\bold c:[\lambda]^2 \rightarrow \sigma$ witness it which means:
\mn
\begin{enumerate}
\item[$(*)'_{\bold c}$]  if (a) then (b) where:
\mn
\begin{enumerate}
\item[$(a)$] 
\begin{enumerate}
\item[$(\alpha)$]   as above
\sn
\item[$(\beta)$]  $h:\bold i_\iota \rightarrow \sigma$ and
  $D$ is an ultrafilter on $\bold i_{1-\iota}$
\end{enumerate}
\sn
\item[$(b)$]   for some $\alpha_0 < \alpha_1 < \mu$ we have
\sn
\begin{itemize}
\item   if $i < \bold i_\iota$ then $\{j < \bold i_{1-\iota}:\bold
  c\{\zeta^\iota_{\alpha_\iota,i},\zeta^{1-\iota}_{\alpha_{1-\iota,j}}\}
  = h(i)\}$ belongs to $D$.
\end{itemize}
\end{enumerate}
\end{enumerate}
\mn
4) Let $\Pr_1(\lambda,\mu,\sigma,\bar\theta)$ be defined similarly to
part (1) but inside $(*)_{\bold c}$ for some 
$\gamma < \sigma$ the function $h$ is constantly $\gamma$.

\noindent
4A) Let $\Pr^*_0(\lambda,\sigma,\theta)$ mean that there is
$\bold c:[\lambda]^2 \rightarrow \sigma$ witnessing it, which means:
\mn
\begin{enumerate}
\item[$(*)^1_{\bold c}$]   like $(*)_{\bold c}$ in part (1) but
\sn
\begin{itemize}
\item  \underline{in $(a)(\gamma)$}: for
  some $\gamma < \sigma,h$ is constantly zero except that $h(0,0)=\gamma$
\end{itemize}
\end{enumerate}
\mn
5) $\lambda \rightarrow [\mu]^2_\theta$ means some $\bold
c:[\lambda]^2 \rightarrow \theta$ witness it which maens: if $A \in
[\lambda]$ and $i < \theta$ then for some $\alpha_1 < \alpha_2$ from
$\bar A$ we have $\bold c\{\alpha_1,\alpha_2\} = i$; follows from
$\Pr_1(\lambda,\lambda,\sigma,\theta),\sigma \ge 1$.

\noindent
6) $\lambda \rightarrow [\lambda;\lambda]^2_\theta$ means some $\bold
c:[\lambda]^2 \rightarrow \theta$ witness it meaning: if $A_1 \in
[\lambda]^\lambda,A_2 \in [\lambda]^\lambda$ and $i < \theta$ then for
some $\alpha_1 \in A,\alpha_2 \in A_2$ we have $\bold
c\{\alpha_1,\alpha_2\} = 1$ (note this follows from
$\Pr_1(\lambda,\lambda,\sigma,\theta),\sigma \ge 1$.
\end{definition}

\noindent
After Todorcevic \cite{To}:
\begin{claim}
\label{w14}
1) If $\lambda > \aleph_0$ is regular, \then \,
$\Pr_1(\lambda^+,\lambda^+,\kappa,\kappa)$.

\noindent
2) $\aleph_1 \nrightarrow (\aleph_1;\aleph_1)^2_{\aleph_1}$.
\end{claim}

\begin{PROOF}{\ref{w14}}
1) By \cite[Ch.IV]{Sh:e}, see history there.

\noindent
2) By Moore \cite{Mo06}.
\end{PROOF}

\noindent
Relatives we shall use are \ref{c28}, \ref{e3}.

\begin{notation}
\label{w17}
1) For a group $G$ let $\seb_G(A) = \seb(A,G)$ be the subgroup of $G$
generated by $A$.

\noindent
2) Let $\bold C_G(A) := \{g \in G:ag=ga$ for every $a \in G\}$.
\end{notation}

\noindent
The following will be used in \S(3A),\S4.
\begin{definition}
\label{w22}
Let $\lambda \ge \theta \ge \sigma$.

\noindent
1) Let $\lambda^{[\theta;\sigma]} = \min\{|\cP|:\cP
\subseteq [\lambda]^\sigma$ and for every $u \in [\lambda]^\theta$ we can
find $\bar u = \langle u_i:i < i_*\rangle$ such that $i_* <
\theta,\cup\{u_i:i < i_*\} = u$ and $[u_i]^\sigma \subseteq \cP\}$;
if $\lambda = \lambda^\sigma$ then $\cP = [\lambda]^\sigma$ witness
$\lambda = \lambda^{[\theta;0]}$ trivially.

\noindent
2) Let $\lambda^{\langle \theta;\sigma\rangle} = \min\{|\cP|:\cP
\subseteq [\lambda]^\sigma$ and for every $u \in [\lambda]^\theta$
there is $v \in [u]^\sigma$ which belongs to $\cP\}$.

\noindent
3) Let $\lambda^{(\theta;\sigma)} = \min\{|\cP|:\cP
\subseteq [\lambda]^\sigma$ and for every $u \in [\lambda]^\theta$
there is $v \in \cP$ such that $|v \cap u| = \sigma\}$.

\noindent
4) If $\mu^+ < \lambda$ and no cardinal in the interval 
$(\mu^+,\lambda)$ is a fix
point then for some regular $\sigma \le \theta \in (\mu,\lambda)$ we
have $\lambda^{\langle \theta;\sigma\rangle} = \lambda$.
\end{definition}
\newpage

\section {Amalgamation Bases}\label{1}
\bigskip

We try to see if there are amalgamation bases
$(K^{\gk}_\lambda,\le_{\gk})$ and if they are dense in a strong sense:
determine for which regular $\kappa$, the $\kappa$-majority of $M \in
K^k_\lambda$ are amalgamation bases.

Another problem is $\Lim_{\gk} = \{(\lambda,\kappa)$: there is a medium
limit model in $K^{\gk}_\lambda\}$, see \cite{Sh:88r}.  This seems
close to the existence of $(\lambda,\kappa)$-limit models, see
\cite{Sh:900}, \cite{Sh:906} and \cite{Sh:877}.  In particular, 
can we get the following:
\bigskip

\begin{question}
\label{a0}
If the set of $M \in \bold K_\lambda$, which
are an amalgamation base, is dense in $(\bold K_\lambda,\subseteq)$, \then
\, in $(\bold K_\lambda,\subseteq)$ there is a
$(\lambda,\aleph_0)$-limit model.

We shall return to this in \S(3C). 
\end{question}

\begin{convention}
\label{a5}
1) $\gk = (K_{\gk},\le_{\gk})$ is an a.e.c. but for simplicity we
allow an empty model, $\le_{\gk}$ than anybody else.

\noindent
2) $\bold K = K_{\gk}$, but we may write $\bold K$ instead of $\gk$
when not said otherwise.
\end{convention}

\begin{definition}
\label{a11}
1) For $M \in K_{\gk}$ and $\mu \ge \LST_{\gk}$ and ordinal $\varp$ we
define an equivalence relation $E_{M,\mu,\varp} = E^M_{\mu,\varp} =
E^M_\varp = E_\varp$ by induction on $\varp$.
\medskip

\noindent
\underline{Case 1}:  $\varp = 0$.

$E^M_\varp$ is the set of pairs $(\bar a_1,\bar a_2)$ such that:
$\bar a_1,\bar a_2 \in {}^{\mu >}M$ have the same length and realize the
same quantifier free type, moreover, for $u \subseteq \ell g(\bar
a_1)$ we have $M \rest (\bar a_1 \rest u) \le_{\gk} M 
\Leftrightarrow M \rest (\bar a_2
\rest u) \le_{\gk} M$.
\medskip

\noindent
\underline{Case 2}:  $\varp$ is a limit ordinal.

$E_\varp = \cap\{E_\zeta:\zeta < \varp\}$.
\medskip

\noindent
\underline{Case 3}:  $\varp = \zeta +1$.

$\bar a_1 E^M_\varp \bar a_2$ \Iff \, for every $\ell \in
\{1,2\},\alpha < \mu$ and $\bar b_\ell \in {}^\alpha M$ there is $\bar
b_{3-\ell} \in {}^\alpha M$ such that $(\bar a_1 \char 94 \bar b_1)
E_\zeta (\bar a_2 \char 94 \bar b_2)$.
\end{definition}

\begin{definition}
\label{a14}
For $\mu > \LST_{\gk}$ and ordinal $\varp$ we define $K_{\gk,\varp} =
\bold K_\varp,K_{\gk,\mu,\varp} = \bold K_{\mu,\varp}$ 
by induction on $\varp$ by (well the notation $\bold K_\varp$ from
here and $\bold K_\lambda = \{M \in \bold K:\|M\| = \lambda\}$ are in
conflict, but usually clear from the context):
\mn
\begin{enumerate}
\item[$(a)$]  $\bold K_\varp = \bold K_{\gk}$ for $\varp=0$;
\sn
\item[$(b)$]   for $\varp$ a limit ordinal $\bold K_\varp
  = \cap\{\bold K_\zeta:\zeta < \varp\}$;
\sn
\item[$(c)$]  for $\varp = \zeta +1$, let $\bold K_\varp$
  be the class of $M_1 \in \bold K_\zeta$ such that: if 
$M_1 \subseteq M_2 \in \bold K_\zeta,\bar a_1 \in {}^{\mu>} M_1,
\bar b_2 \in {}^{\mu >}(M_2)$ \then \, for some 
$b_1 \in {}^{\mu >}M_1$ we have $\bar a
  \char 94 \bar b_1 E^{M_1}_\zeta \bar a \char 94 \bar b_2$.
\end{enumerate}
\end{definition}

\begin{claim}
\label{a17}
For every $\varp$:
\mn
\begin{enumerate}
\item[$(a)$]  for every $M_1 \in \bold K_{\gk}$ there
  is $M_2 \in \bold K_\varp$ extending $H$;
\sn
\item[$(b)$]  $E^M_\varp$ has $\le \beth_{\varp +1}(\mu)$
equivalence classes, hence in clause (a) we can\footnote{We can improve 
the bound a little, e.g.if $\mu = \chi^+$ then 
$\beth_{\varp +1}(\chi)$ suffices.} add $\|M_2\| \le \|M_1\| +
\beth_{\varp +1}(\mu)$;
\sn
\item[$(c)$]  $M_1 \in \bold K_{\mu,\varp}$ \when \, $\bold
  K_\varp$ has amalgamation and $M_1 \subseteq
  M_2,M_2 \in \bold K_\varp$ implies:
\sn
\begin{itemize}
\item    if $\zeta < \varp,\bar a \in {}^{\mu >}(M_1),
\bar b_2 \in {}^{\mu >}(M_2)$ then there is $\bar b_1 \in
{}^{\ell g(\bar b)}(M_1)$ 
such that $\bar a \char 94 \bar b_1 E^{M_2}_{\mu,\zeta} \bar a 
\char 94 \bar b_2$;
\end{itemize}
\sn
\item[$(d)$]   if $I$ is a $(< \mu)$-directed partial order
  and $M_s \in \bold K_\varp$ is $\subseteq$-increasing 
with $s \in I$, \then \, $M = \bigcup\limits_{s} M_s \in \bold K_\varp$;
\sn
\item[$(e)$]   if $H_1 \subseteq H_2$ are from $\bold K_\varp$ then $H_1
  \prec_{\bbL_{\infty,\mu,\varp}(\gk)} H_2$;
\sn
\item[$(f)$]   if $\varp = \mu,\mu = \cf(\mu)$ or $\varp = \mu^+$, 
and $H_1 \subseteq H_2$ are from $\bold K_{\mu,\varp}$, \then \,
$H_1 \prec_{\bbL_{\mu,\mu}} H_2$.
\end{enumerate}
\end{claim}

\begin{PROOF}{\ref{a17}}
We can prove this by induction on $\varp$.  The details should be clear.
\end{PROOF}

\begin{definition}  
\label{a31}
1) We say $M_0 \in \bold K_\lambda$ is a $\bar\chi$-amalgamation base
\when \,: $\bar\chi = (\chi_1,\chi_2)$ and $\chi_\ell \ge \|M\|$ and
if $M_0 \le_{\gk} M_\ell \in \bold K_{\chi_\ell}$ for $\ell=1,2$,
\then \, for some $M_3 \in \bold K_{\gk}$ 
which $\le_{\gk}$-extend $M$, both $M_1$ and $M_2$
can be $\le_{\gk}$-embedded into $M_3$ over $M_0$.

\noindent
2) We may replace ``$\chi_\ell$" by ``$< \chi_\ell$" with obvious
meaning (so $\chi_\ell > \|M_0\|$).  If $\chi_1 = \chi_2$ we may write
$\chi_1$ instead of $(\chi_1,\chi_2)$.  If $\chi_1 = \chi_2 = \lambda$ we
may write ``amalgamation base".

\noindent
3) We say $\bold K_{\gk}$ has $(\bar\chi,\lambda,\kappa)$-amalgamation bases
\when \, the $\kappa$-majority of $M \in \bold K_\lambda$ is a
$\bar\chi$-amalgamation base where:

\noindent
3A) We say that the $\kappa$-majority of $M \in \bold K_\lambda$ satisfies
$\psi$ \when \,  some $F$ witnesses it, which means:
\mn
\begin{enumerate}
\item[$(*)$]
\begin{enumerate}
\item[(a)]  $F$ is a function with\footnote{We may use $F$
    with domain $\{\bar M:M = \langle M_i:i < j\rangle$ is increasing,
    each $M_i \in \bold K$ has universe an ordinal $\alpha \in
    [\lambda,\lambda^+)\}$; see \cite{Sh:88a}.}
 domain $\{M \in \bold K_{\gk}:M$ has universe an 
ordinal $\in [\lambda,\lambda^+)\}$;
\sn
\item[(b)]   if $M \in \Dom(F)$ then $M \le_{\gk} F(M) \in \Dom(F)$;
\sn
\item[(c)]  if $\langle M_\alpha:\alpha \le
  \kappa\rangle$ is increasing continuous, $M_\alpha \in \Dom(F)$ and
$M_{2 \alpha+2} = F(M_{2\alpha +1})$ for every $\alpha < \kappa$,
\then \, $M_\kappa$ is a $\bar\chi$-amalgamation base.
\end{enumerate}
\end{enumerate}
\mn
4) We say the pair $(M,M_0)$ is an $(\chi,\mu,\kappa)$-amalgamation
base (or amalgamation pair) \when \,: $M \le_{\gk} M_0 \in \bold K_{\gk},
\|M\| = \kappa,\|M_0\| = \mu$ and if $M_0 \le_{\gk} M_\ell \in
\bold K_{\le \chi}$ for $\ell=1,2$, \then \, for some
$M_3,f_1,f_2$ we have $M_0 \le_{\gk} M_3 \in \bold K_{\gk}$ and
$f_\ell \le_{\gk}$-embeds $M_\ell$ into $M_3$ over $M_0$.

\noindent
5) Let $\AM_{\bold K} = \AM_{\gk}$ be the class of pairs
$(\lambda,\kappa)$ such that $\bold K$ has
$((\lambda,\lambda),\lambda,\kappa)$-amalgamation bases.
\end{definition}

\begin{definition}
\label{a32}
1) For $\gk,\bar\chi,\lambda,\kappa$ as above and 
$S \subseteq \lambda^+$ (or $S
\subseteq \Ord$ but we use $S \cap \lambda^+$) we say $\gk$ has
$(\bar\chi,\lambda,\kappa,S)$-amalgamation bases \when \, there is a
function $F$ such that:
\mn
\begin{enumerate}
\item[$(*)_F$] 
\begin{enumerate}
\item[(a)]  $F$ is a function with domain $\{\bar M:\bar
  M$ is a $\le_{\gk}$-increasing continuous
sequence of members of $\bold K_{\gk}$ each with universe 
an ordinal $\in [\lambda,\lambda^+)$ 
and length $i+1$ for some $i \in S\}$;
\sn
\item[(b)]   if $\bar M = \langle M_i:i \le j\rangle \in \Dom(F)$ then:
\begin{enumerate}
\item[$(\alpha)$]  $F(\bar M) \in \bold K_{\gk}$;
\sn
\item[$(\beta)$]  $M_j \le_{\gk} F(\bar M)$;
\sn
\item[$(\gamma)$]  $F(\bar M)$ has universe an ordinal
  $\in[\lambda,\lambda^+]$;
\end{enumerate}
\sn
\item[(c)]  if $\delta = \sup(S \cap \delta) < \lambda^+$ has
  cofinality $\kappa$ and 
$\bar M = \langle M_i:i \le \delta\rangle$ is 
$\le_{\gk}$-increasing continuous and for every $j < \kappa$ 
we have $j \in S \Rightarrow \bar M_{j+1} = F(\bar M \rest (j+1))$
  hence $\bar M \rest (j+1) \in \Dom(F)$ 
\then \, $M_\delta$ is a $\bar\chi$-amalgamation base.
\end{enumerate}
\end{enumerate}
\mn
2) We say $\gk$ has the weak
$(\bar\chi,\lambda,\kappa,S)$-amalgamation bases 
\when \, above we replace clause (c) by:
\mn
\begin{enumerate}
\item[(c)$'$]  if $\langle M_i:i < \lambda^+\rangle$ is
  $\le_{\gk}$-increasing and $j \in S \cap \lambda^+ \Rightarrow
  M_{j+1} = F(\bar M \rest (j+1))$ \then \, for some club $E$ of
  $\lambda^+$ we have $\delta \in E$ and $\cf(\delta) = \kappa
  \Rightarrow M_\delta$ is a $\bar\chi$-amalgamation base.
\end{enumerate}
\mn
3) We say $\gk$ has $(\bar\chi,\lambda,W,S)$-amalgamation bases when
$W \subseteq \lambda^+$ is stationary and in part (2) we replace (in
the end of $(c)'$, ``$\delta \in E$ and $\cf(\delta) = \kappa$" by
``$\delta \in E \cap W$".
\end{definition}

\begin{PROOF}{\ref{a34}}
Easy.
\end{PROOF}

\begin{claim}
\label{a37}
1) If $\lambda = \kappa > \LST_{\gk}$ is a weakly compact cardinal and $M \in
\bold K_{\kappa,1}$, see Definition \ref{a14} 
\then \, $M$ is a $\kappa$-amalgamation base.

\noindent
2) If $\kappa$ is compact cardinal and $\lambda = \lambda^{< \kappa}$
and $M \in \bold K_{\kappa,1}$ has cardinality $\lambda$, \then \, $M$
is a $(< \infty)$-amalgamation base; so $\gk$ has $(<
\infty,\lambda,\ge \kappa)$-amalgamation bases.

\noindent
3) In part (2), $\kappa$ has to satisfy only: if $\Gamma$ is a set $\le
\lambda$ of sentences from $\bbL_{\LST(\gk)^+,\aleph_0}$ and every $\Gamma'
\in [\Gamma]^{< \kappa}$ has a model, then $\Gamma$ has a model.
\end{claim}

\begin{PROOF}{\ref{a37}}
Use the representation theorem for a.e.c. from \cite[\S1]{Sh:88r}
which is quoted in \ref{w32} here and the definitions.
\end{PROOF}

\begin{conclusion}
\label{a41}
If the pair $(\lambda,\kappa)$ is as in \ref{a37}, \then \, $\gk$ has
$(\lambda,\kappa)$-amalgamation bases; see \ref{a31}(3).
\end{conclusion}

\begin{claim}
\label{a41}
If $\gk,\bold K$ are as in \ref{w32} and the universal class $\bold
K$, i.e. $(\bold K,\subseteq)$ have $(\bar\chi,\lambda,\kappa)$
amalgamation and $\lambda \ge \LST(\gk)$, \then \, so does $\gk$.
\end{claim}

\begin{PROOF}{\ref{a41}}
Easy.
\end{PROOF}

\noindent
A surprising result says that in some singular 
cardinals we have ``many" amalgamation bases.
\begin{claim}  
\label{a44}
If $\mu$ is a strong limit cardinal and $\cf(\mu) > \LST_{\gk}$ is 
a measurable cardinal (so $\mu$ is measurable or $\mu$ is singular but
the former case is covered by \ref{a37}(1)) \then\, $\bold K$ 
has $(\mu,\cf(\mu))$-amagamation bases.
\end{claim}

\begin{PROOF}{\ref{a44}}  
By \ref{a41} \wilog \, $\gk$ is a universal class $\bold K$.
\Wilog \, $\mu$ is a singular cardinal (otherwise the result
follows by Claim \ref{a37}).
Let $\kappa = \cf(\mu),D$ a normal ultrafilter on
$\kappa$ and let $\langle \mu_i:i < \kappa\rangle$ be an increasing
sequence of cardinals with limit $\mu$ such that $\mu_0
\ge \LST_{\gk} + \kappa$.

We choose $\bold u$ such that:
\mn
\begin{enumerate}
\item[$(*)_1$]
\begin{enumerate}
\item[(a)]   $\bold u = \langle \bar u_\alpha:\alpha < \mu^+\rangle$;
\sn
\item[(b)]  $\bar u_\alpha = \langle u_{\alpha,i}:i < \kappa\rangle$;
\sn
\item[(c)]  $u_{\alpha,i} \in [\alpha]^{\mu_i}$ is
  $\subseteq$-increasing with $i$;
\sn
\item[(d)]  $\alpha = \bigcup\limits_{i < \kappa} u_{i,\kappa}$;
\sn
\item[(e)]  if $\alpha < \beta < \mu^+$, \then \, $u_{\alpha,i}
  \subseteq \alpha_{\beta,i}$ for every $i < \kappa$ large enough
\end{enumerate}
\end{enumerate}
\mn
For transparency we allow $=^M$ to be non-standard, i.e. just a
congruence relation on $M$.

We now choose functions $F,G$ by:
\mn
\begin{enumerate}
\item[$(*)_2$] 
\begin{enumerate}
\item[(a)]  $\dom(F) = \{M \in \bold K_{\gk}:M$ has universe some
  $\alpha \in [\mu,\mu^+)\}$;
\sn
\item[(b)]   for $\alpha \in [\mu,\mu^+)$
let $\cM_\alpha = \{M \in \bold K_{\gk}:M$ has universe $\alpha$
\sn
\item[(c)]  for $M \in \cM_\alpha,u \subseteq \alpha$ let $M[u] = M
  \rest c \ell(u,M)$ and let $M^{[i]} = M[u_{\alpha,i}]$, hence $u
  \subseteq \alpha \Rightarrow M[u] \le_{\gk} M$;
\sn
\item[(d)]   if $M \in \dom(F)$ has universe $\alpha$ \then \, $M^+ =
  F(M)$ satisfies:
\sn
\begin{enumerate}
\item[$(\alpha)$]  $M \subseteq M^+ \in \bold M_{\alpha + \lambda}$
  (equivalently $M \le_{\gk} M^+ \in \cM_{\alpha + \lambda}$)
\sn
\item[$(\beta)$]  if $i < \kappa$ and $M^{[i]} \le N \in \bold
  K_{\mu_i}$, then exactly one of the following occurs:
\begin{itemize}
\item  there is an embedding of $N$ into $M^+$ over $M^{[i]}$
\sn
\item  there is no $M' \in \bold K$ extending $M^+$ and an embedding
  of $N$ into $M'$ over $M^{[i]}$
\end{itemize}
\end{enumerate}
\end{enumerate}
\end{enumerate}
\mn
This is straightforward.  
It is enough to prove that $F$ witnesses that $\gk$ has
$(\mu,\kappa)$-amalgamation bases, i.e. using $F(\langle M_i:i \le
j\rangle) = F(M_j)$.

For this it suffices:
\mn
\begin{enumerate}
\item[$(*)_3$]  $M^1,M^2$ can be amalgamated over $M_\kappa$ \when \,:
\sn
\begin{enumerate}
\item[$(a)$]  $\langle M_i:i \le \kappa\rangle$ is $\le_{\gk}$-increasing
  continuous;
\sn
\item[$(b)$]  $M_i \in \bold K^{\gk}_\mu$ has universe $\alpha_i,M_{2i +2} \in
  \cM_{\alpha_{2i+2}}$;
\sn
\item[$(c)$]  $F(M_{2i+1}) = M_{2i+2}$;
\sn
\item[$(d)$]  $M_\kappa \subseteq M^1 \in \bold K^{\gk}_\mu$ and $M_\kappa
  \subseteq M^2 \in \bold K^{\gk}_\mu$.
\end{enumerate}
\end{enumerate}
\mn
We can find an increasing (not necessarily continuous) sequence
$\langle \varp(i):i < \kappa\rangle$ of ordinals $< \kappa$ such that
$i < j < \kappa \Rightarrow u_{\alpha_{\varp(i)},j} \subseteq
u_{\alpha_{\varp(j)},j}$ and $u_i = u_{\alpha_{\varp(i)},i},u^*_i \in
[\alpha_\kappa]^{\le\mu_i}$ is $\subseteq$-increasing.

\Wilog \, $M',M''$ has universe $\beta = \alpha_\kappa + \mu$.

Now,
\mn
\begin{enumerate}
\item[$(*)$]  let $\langle u^*_i:i < \kappa\rangle$ be
  $\subseteq$-increasing with union $\beta$ such that:

$i < \kappa \Rightarrow u_i \subseteq u^*_i$.
\end{enumerate}
\mn
Notice that:
\mn
\begin{enumerate}
\item[$\boxplus$]  it suffices to prove that: for every $i <
  \kappa,M^1[u^*_i],M^2[u^*_i]$ can be 
$\le_{\gk}$-embedded into $M_\kappa$ over $M_\kappa[u_i]$ 
(you can use its closure); say $h^\iota_i$ is a
 $\le_{\gk}$-embedding of $M^\iota[u^*_i]$ into $M_\kappa$ over
$M_\kappa[u_i]$.
\end{enumerate}
\mn
It suffices to prove $\boxplus$ by taking ultra-products, i.e. let
$N_i$ be $(\mu^+,M_\kappa,M^\iota,M^\iota[u_\iota]u_i,
h^\iota_i)_{\iota=1,2}$ and let $D$ be a normal ultrafilter on
$\kappa$ and ``chase arrows" in $\prod\limits_{i < \kappa} N_i/D$.  It 
is possible to prove $\boxplus$ by the choice of $F$.
\end{PROOF}

\begin{claim}
\label{a47}
1) Assume $\kappa > \theta > \LST_{\gk},\theta$ is a measurable and
$\kappa$ is weakly compact.  \then \, $\gk$ has
$(\kappa,\theta)$-amalgamation bases.

\noindent
2) Assume $\kappa,\theta$ are measurable cardinals $> \LST_{\gk}$ and
$\mu > \kappa + \theta$ is strong limit singular of cofinality
$\kappa$.  \Then \, $\gk$ has $(\mu,\theta)$-amalgamation bases.

\noindent
3) If $\kappa > \theta > \LST_{\gk},\theta$ is a measurable cardinal
and $\{M \in \bold K^{\gk}_\kappa:M$ is a $(\chi_1,\chi_2)$-amagamation
base$\}$ is $\le_{\gk}$-dense in $\bold K^{\gk}_\kappa$, \then \, $\gk$ has
$(\chi_1,\chi_2,\kappa,\theta)$-amalgamation bases. 
\end{claim}

\begin{PROOF}{\ref{a47}}
1) As $\gk$ has $(\kappa,\kappa)$-amalgamation bases by \ref{a37}(1)
we can apply part (3) of \ref{a47} with
$(\kappa,\kappa,\kappa,\theta)$ here standing for
$(\chi_1,\chi_2,\kappa,\theta)$ there.

\noindent
2) Similarly to part (1) using \ref{a44} instead of \ref{a37}(1).

\noindent
3) Similar to the proof of \ref{a44}, that is, we replace $\boxplus$
by Claim \ref{a48} and $(*)_2$ by:
\mn
\begin{enumerate}
\item[$(*)^1_2$]  if $M \in \bold K_\alpha$, then $F(M)$ is a member of
  $K_{\gk}$ which is a $\bar\chi$-amaglamation base and $M \le_{\gk}
  F(M)$.
\end{enumerate}
\end{PROOF}

\begin{claim}
\label{a48}
Assume $\bar M = \langle M_i:i \le \kappa\rangle$ is
$\le_{\gk}$-increasing (not necessarily continuous) and $M_\kappa :=
\bigcup\limits_{i < \kappa} M_i$ is of cardinality $\le
\min\{\chi_1,\chi_2\}$ and each $M_i$ is a $\bar\chi$-amalgamation
base.  \Then \, $M_\kappa$ is a $\bar\chi$-amalgamation base.
\end{claim}

\begin{claim}
\label{a49}
1) In \ref{a44}, we can replace ``$(\mu,\cf(\mu))-amalgamation base" by
``(\mu,\cf(\mu),S)$-amalgamation base" for any unbounded subset $S$ of
$S$.

\noindent
2) Similarly in \ref{a47}.
\end{claim}

\begin{question}
\label{a50}
1) What can $\AM_{\gk} = \{(\lambda,\kappa):\gk$ has
$(\lambda,\kappa)$-amalgamation, $\lambda > \LST_{\gk}\}$ be?

\noindent
2) What is $\AM_{\gk}$ for $\gk = \bold K_{\exlf}$?

\noindent
3) Suppose we replace $\kappa$ by stationary $W \subseteq \{\delta <
\lambda^+:\cf(\delta) = \kappa\}$.  How much does this matter?
\end{question}

\begin{discussion}
\label{a53}
1) May be helpful for analyzing $\AM_{\bold K_{\lf}}$ but also of self
interest is analyzing $\gS_{k,n}[\bold K]$ with $k,n$ possibly
infinite, see \cite[\S4]{Sh:312}.

\noindent
2) In fact for \ref{a50}(3) we may consider Definition \ref{a56}.
\end{discussion}

\begin{definition}
\label{a56}
For a regular $\theta$ and $\mu \ge \alpha$ let:
\mn
\begin{enumerate}
\item[$(A)$]  $\Seq^0_{\mu,\alpha}$ is in the class of $\bar N$ such
  that:
\sn
\begin{enumerate}
\item[$(a)$]  $\bar N = \langle N_i:i \le \alpha\rangle$ is
  $\le_{\gk}$-increasing continuous except for $i=0$;
\sn
\item[$(b)$]  $i \ne 0 \Rightarrow \|N_i\| = \mu$;
\sn
\item[$(c)$]  $|N_0| = 0$;
\end{enumerate}
\sn
\item[$(B)$]  $\Seq^1_{\mu,\alpha} = \{\bold n = (\bar N^1,\bar
  N^2):\bar N^\iota \in \Seq^0_{\mu,\alpha +1}$ and $\beta \le \alpha
  \Rightarrow N^1_\beta = N^2_\beta$ so let $N_\beta = N_{\bold
    n,\beta} = N^1_\beta$;
\sn
\item[$(C)$]  we define the game $\Game_{\bar N,\bold n}$ for $\bold n \in
\Seq^1_{\mu,\alpha}$;
\sn
\begin{enumerate}
\item[$(a)$]  a play last $\alpha +1$ moves and is between $\AAM$ and $\AM$;
\sn
\item[$(b)$]  during a play a sequence $\langle (M_i,M'_i,f_i):i \le
  \alpha\rangle$ is chosen such that:
\begin{enumerate}
\item[$(\alpha)$]  $M_i \in \bold K_\lambda$ is
  $\le_{\gk}$-increasing continuous;
\sn
\item[$(\beta)$]  $f_i$ is a $\le_{\gk}$-embedding of
  $N_{\bold n,i}$ into $M_i$ and even $M'_i$;
\sn
\item[$(\gamma)$]  $f_i$ is increasing continuous for limit
  $i,f_\delta = \bigcup\limits_{i < \delta} f_i$ and $f_0$ is empty;
\sn
\item[$(\delta)$]  $M_i \le_k M'_{i+1} \le M_{i+1}$ and for $i$
  limit or zero $M'_i = M_i$;
\end{enumerate}
\sn
\item[(c)]
\begin{enumerate}
\item[$(\alpha)$]   if $i=0$ in the $i$-th move
  \underline{first} $\AM$ chooses $M_0$ and
\newline
 \underline{second} $\AAM$ chooses $f_0 = \emptyset,M'_0 = M_0$;
\sn
\item[$(\beta)$]  if $i=j+1$, in the $i$-th move
\underline{first} $\AM$ chooses $f_i,M'_i$ and 
\underline{second} $\AA$  chooses $M_i$;
\sn
\item[$(\gamma)$]   if $i$ is a limit ordinal:
  $M_i,f_i,M'_i$ are determined;
\sn
\item[$(\delta)$]   if $i=\alpha+1$,
\underline{first} $\AAM$ chooses $N_i \in \{N^1_\alpha,N^2_\alpha\}$
and then this continues as above;
\end{enumerate}
\sn
\item[(d)]  the player $\AMM$ wins when $\AM$ has no legal move;
\end{enumerate}
\sn
\item[(D)]  let $\Seq_{\gk}$ be the set of $\lambda,\mu,\theta$ 
such that there is $\bold n$ satisfying:
\sn
\begin{enumerate}
\item[(a)]  $\bold n \in \Seq^1_{\mu,\theta}$;
\sn
\item[(b)]  $N^1_{\bold n,\theta +1},N^2_{\bold n,\theta +2}$ cannot
  be amalgamated over $N_{\bold n,\theta} (= N^\iota_{\bold
    n,\theta},\iota=1)$;
\sn
\item[(c)]  in the game $\Game_{\bold n}$, the player $\AM$ has a
  winning strategy.
\end{enumerate}
\end{enumerate}
\end{definition} 

\begin{question}
\label{a59}
1) What can be $\Seq_{\gk}$ for $\gk$ an a.e.c. with $\LST_{\gk} =
\chi$?

\noindent
2) What is $\Seq_{\bold K_{\lf}}$?
\end{question}

\begin{claim}
\label{a34}
Let $S$ be the class of odd ordinals.

\noindent
1) If $\gk$ has $(\bar\chi,\lambda,\kappa,S)$-amalgamation \then \, $\gk$ has
$(\bar\chi,\lambda,\kappa)$-amalgamation.

\noindent
2) If $\lambda = \lambda^{< \kappa}$ \then \, also the inverse holds.
\end{claim}
\newpage

\section {Definability} \label{2}

The notion of ``$a \in M_2 \backslash M_1$ is definable over $M_1$" is
clear for first order logic, $M_1 \prec M_2$.  But in a class like
$\bold K_{\lf}$ we may wonder.
\begin{claim}
\label{n60}
Below (i.e. in \ref{n61} - \ref{n70}) we can replace $\bold K_{\lf}$ by:
\mn
\begin{enumerate}
\item[$(*)$]  $\gk$ is a a.e.c. and one of the following holds:
\sn
\begin{enumerate}
\item[(a)]  $\gk$ is a universal, so $\bold k_1 = \gk \rest \{M \in
  K_{\gk}:M$ is finitely good$\}$ determine $\gk$;
\sn
\item[(b)]  like (a) but $\gk_1$ is closed under products;
\sn
\item[(c)]  like (a), but in addition:
\begin{enumerate}
\item[$(\alpha)$]  $0_{\gk} = 0_{\gk_1}$ is an individual constant;
\sn
\item[$(\beta)$]   if $M_1,M_2 \in K_{\gk_1}$ then $N = M_1 \times M_2
  \in K_{\gk_1}$; moreover $f_\ell:M_\ell \rightarrow N$ is a
  $\le_{\gk_1}$-embedding for $\ell=1,2$ where:
\begin{itemize}
\item  $f_1(a_1) = (a_1,0_{M_2})$;
\sn
\item  $f_2(a_2) = (0_{M_1},a_2)$.
\end{itemize}
\end{enumerate}
\end{enumerate}
\end{enumerate}
\end{claim}

\begin{discussion}
\label{n60d}
Can we in (c) define types as in \ref{n61} such that they behave
suitably (i.e. such that \ref{n67}, \ref{n70} below works?)  We need $c
\ell(A,M)$ to be well defined.
\end{discussion}

\begin{definition}  
\label{n61}
1) For $G \subseteq H \in \bold K_{\lf}$ we let $\uniq(G,H) = \{x \in
H$: if $H \subseteq H^+ \in \bold K_{\lf},y \in H^+$ and
$\tp_{\bs}(y,G,H^+) = \tp_{\bs}(x,G,H)$ then $y=x\}$.

\noindent
1A) Above we let $\uniq_\alpha(G,H) = \uniq^1_\alpha(G,H) =
\{\bar x \in {}^\alpha H$: if $H
\subseteq H^+ \in \bold K_{\lf}$, then no $\bar y \in {}^\alpha(H^+)$
realizes $\tp_{\bs}(\bar x,G,H)$ in $H^+$ and satisfies $\Rang(\bar y) \cap
\Rang(\bar x) \subseteq G\}$.

\noindent
1B) Let $\uniq^2_\alpha(G,H)$ be defined as in (1A) but in the end
``$\Rang(\bar x) = \Rang(\bar  y)$".

\noindent
1C) Let $\uniq^3_\alpha(G,H)$ be defined as in (1A) but in the end
``$\bar x = \bar y$".

\noindent
2) For $G_1 \subseteq G_2 \subseteq G_3 \in \bold K_{\lf}$ let
$\uniq(G_1,G_2,G_3) = \{x \in G_2$: if $G_3 \subseteq G \in \bold
K_{\lf}$ then for no $y \in G \backslash G_2$ do we have
$\tp_{\bs}(y,G_1,G) = \tp_{\bs}(x,G_1,G_2)$.
\end{definition}

\begin{question}  
\label{n64}
1) Given $\lambda$, can we bound $\{|\uniq(G,H)|:G \subseteq H \in
\bold K_{\lf}$ and $|G| \le \lambda\}$.

\noindent
2) Can we use the definition to prove ``no $G \in \bold
K^{\lf}_{\beth_\omega}$ is universal"?
\end{question}

\noindent
To answer \ref{n64}(1) we prove $2^\lambda$ is a bound and more; toward this:
\begin{claim}  
\label{n67}
If (A) then (B), where:
\mn
\begin{enumerate}
\item[(A)]
\begin{enumerate}
\item[(a)]  $G_n \in \bold K_{\lf}$ for $n < n_*;n_*$ may
  be any ordinal but the set $\{G_n:n < n_*\}$ is finite;
\sn
\item[(b)]  $h_{\alpha,n}:I \rightarrow G_n$ for $\alpha
  < \gamma_*,n < n_*$;
\sn
\item[(c)]  if $s \in I$,  \then \, the set 
$\{(G_n,h_{\alpha,n}(s)):\alpha < \gamma_*$ and $n<n_*\}$ 
is finite;
\end{enumerate}
\sn
\item[(B)]  there is $(H,\bar a)$ such that:
\sn
\begin{enumerate}
\item[(a)]  $H \in \bold K_{\lf}$;
\sn
\item[(b)]  $\bar a = \langle a_s:s \in I \rangle$ generates $H$;
\sn
\item[(c)]  if $s_0,\dotsc,s_{k-1} \in I$ \then
\newline
$\tp_{\at}(\langle a_{s_\ell}:\ell < k\rangle,\emptyset,H) = 
\bigcap\limits_{n,\alpha} \tp_{\at}(\langle
  h_{\alpha,n}(s_0),\dotsc,h_{\alpha,n}(s_{k-1}\rangle,\emptyset,G_n)$;
\sn
\item[(d)]  the mapping $b_s \rightarrow a_s$ for $s \in I_*$ 
embeds $H_*$ into $H$ \when\, :
\begin{enumerate}
\item[$(*)$]  $H_* \subseteq G_n$ for $n<n_*,I_* \subseteq I,\langle
b_s:s \in I_*\rangle$ list the elements of 
$H_*$ (or just generates it) and $\alpha < \gamma_* \wedge s \in I_* \wedge 
n < n_* \Rightarrow h_{\alpha,n}(s) = b_s$.
\end{enumerate}
\end{enumerate}
\end{enumerate}
\end{claim}

\begin{PROOF}{\ref{n67}}  
Note that:
\mn
\begin{enumerate}
\item[$(*)_1$]   there are $H$ and $\bar a$ such that:
\sn
\begin{enumerate}
\item[$(a)$]  $H$ is a group;
\sn
\item[$(b)$]  $\bar a = \langle a_s:s \in I\rangle$;
\sn
\item[$(c)$]  $a_s \in H$;
\sn
\item[$(d)$]  $\bar a$ generates $H$;
\sn
\item[$(e)$]  for any finite $u \subseteq I$ and atomic
  formula $\varphi(\bar x_{[u]})$ we have $H \models \varphi(\bar a_{[u]})$ \Iff
  \, for every $n < n_*$ we have $G_n \models
  \varphi[\ldots,h_{\alpha,n}(s),\ldots]_{s \in u}$.
\end{enumerate}
\end{enumerate}
\mn
[Why?  Let $G_{\alpha,n} = G_n$ for $\alpha < \gamma_*,n < n_*$ and let $H =
\Pi\{G_{\alpha,n}:n < n_*,\alpha < \gamma_*\}$ and let $a_s = \langle
h_{\alpha,n}(s):(\alpha,n) \in (\gamma_*,n_*)\rangle$ for $s \in I$
and, of course, $\bar a = \langle a_s:s \in I\rangle$.]
\mn
\begin{enumerate}
\item[$(*)_2$]  \Wilog \, $\bar a$ generates $H$, i.e. (B)(d) holds.
\end{enumerate}
\mn
[Why?  Just read $(*)_1$.]
\mn
\begin{enumerate}
\item[$(*)_3$]  If $u \subseteq I$ is finite, then $c \ell(\bar
  a_{[u]},H)$ is finite (and for \ref{n60} it belongs to $K_{\gk}$)
\end{enumerate}
\mn
[Why?  By Clause (A)(c) of Claim \ref{n67}; and for \ref{n60} recalling
\ref{n60}(d).]
\mn
\begin{enumerate}
\item[$(*)_4$]  $H \in \bold K_{\lf}$ (i.e. (B)(a) holds).
\end{enumerate}
\mn
[Why?  By $(*)_2 + (*)_3$; for \ref{n60} use also \ref{n60}(d).]
\mn
\begin{enumerate}
\item[$(*)_5$]  Clause (B)(c) holds.
\end{enumerate}
\mn
[Why?  By $(*)_1(e)$.]
\mn
\begin{enumerate}
\item[$(*)_6$]  Clause (B)(d) holds.
\end{enumerate}
\mn
[Why?  Follows from our choices.]
\end{PROOF}

\begin{claim}  
\label{n70}
If $G_1 \in \bold K^{\lf}_{\le \lambda}$ and $G_1 \subseteq G_2 \in \bold
K_{\lf}$ has cardinality $\le \mu = \mu^\lambda$ (e.g. $G_1,G_2 \in
\bold K^{\lf}_\lambda,\mu = 2^\lambda$), \then\, for some
pair $(G_3,X)$ we have:
\mn
\begin{enumerate}
\item[$\oplus$]
\begin{enumerate}
\item[(a)]  $G_2 \subseteq G_3 \in \bold K^{\lf}_\mu$
\sn
\item[(b)]  $X \subseteq G_3$ has cardinality $\le 2^\lambda$
\sn
\item[(c)]   if $c \in G_3$, then exactly one of the following occurs:
\sn
\begin{enumerate}
\item[$(\alpha)$]  $c \in X$ and $\{b \in G_3:\tp_{\at}(b,G_1,G_3) =
  \tp_{\at}(c,G_1,G_3)\}$ 
is a singleton (and this holds also in $G_4$ when $G_3
\subseteq G_4 \in \bold K_{\lf}$);
\sn
\item[$(\beta)$]   there are $\|G_3\|$ elements of $G$ realizing
  $\tp_{\bs}(a,G_1,G_3)$;
\end{enumerate}
\sn
\item[(d)]  if $\alpha < \lambda^+,\bar a \in {}^\alpha(G_3)$ and
  $p(\bar x_{[\alpha]}) = \tp_{\at}(\bar a,G,G_3)$,
$p'(\bar x_{[\alpha]}) = \tp_{\bs}(\bar a,G,G_3)$, 
\then \, for some non-empty $\cP \subseteq \cP(\alpha)$ 
closed under the intersection of 2 to which $\alpha$
belongs we have:
\begin{enumerate}
\item[$(\alpha)$]  if $\bar a',\bar a'' \in
  {}^\alpha(G_3)$ realizes $p(\bar x_{[\alpha]})$ \then \,
$u := \{\beta < \alpha:(a'_\beta = a''_\beta)\} \in \cP$;
\sn
\item[$(\beta)$]   if $u \in \cP$ \then \, we can find 
$\langle \bar a_\varp:\varp < \|G_3\|\rangle$ a $\Delta$-system
with heart $u$ (i.e. $\bar a_{\varp_1,\beta_1} = 
\bar a_{\varp_2,\beta_2} \Leftrightarrow 
((\varp_1,\beta_1) = (\varp_2,\beta_2)) \vee (\beta_1
= \beta_2 \in u))$, each $\bar a_\varp$ realizing $p(\bar x_{[\alpha]}$
and even $p'(\bar x_{[\alpha]}))$.
\end{enumerate}
\end{enumerate}
\end{enumerate}
\end{claim}

\begin{remark}
1) Can we generalize the (weak) elimination of quantifiers in modules?

\noindent
2) An alternative presentation is to try $G^I_D/\cE$, where:
\mn
\begin{itemize}
\item  $\cE \subseteq \{E:E$ is an equivalence relation on $I$ such
  that $I/E$ is finite$\}$ and $(\cE \ge)$ is directed;
\sn
\item  $G^I_D$ is $G^I \rest \{f:f +G$ and there is $E \in \cE$ such
  that $s E t \Rightarrow f(s) = f(t)\}$.
\end{itemize}
\mn
3) For suitable $(I,D,\cE)$ we have: if $p$ is a set of 
$\le \mu$ basic formulas
with parameters from $G_1 = G^I_D/\cE$ we have: $p$
 is realized in $G_1$ iff every $\varphi_1,\dotsc,\varphi_n,\neg
\varphi_i \in p,\varphi_\ell$ atomic is realized in $G_1$.
\end{remark}

\begin{PROOF}{\ref{n70}}  
We can easily find $G_3$ such that:
\mn
\begin{enumerate}
\item[$(*)_1$]  
\begin{enumerate}
\item[(a)]  $G_2 \subseteq G_3 \in \bold K^{\lf}_\mu$;
\sn
\item[(b)]   if $G_3 \subseteq H \in \bold K_{\lf},
\gamma < \lambda^+,\bar a \in {}^\gamma H$ and
$u = \{\alpha < \gamma:a_\alpha \in G_3\}$, 
\then \, there are $\bar a^\varp \in {}^\gamma(G_3)$
for $\varp < \mu$ such that:
\begin{enumerate}
\item[$(\alpha)$]  $\tp_{\bs}(\bar a^\varp,G_1,G_3) =
  \tp_{\bs}(\bar a,G_1,G_3)$;
\sn
\item[$(\beta)$]   if $\varp,\zeta < \mu$ and
  $\alpha,\beta < \gamma$ and $a^\varp_\alpha = a^\zeta_\beta$ then
$((\varp,\alpha) = (\zeta,\beta)) \vee (\alpha = \beta \wedge 
a^\varp_\alpha = a_\alpha = a^\zeta_\alpha)$.
\end{enumerate}
\end{enumerate}
\end{enumerate}
\mn
We shall prove that
\mn
\begin{enumerate}
\item[$(*)_2$]  $G_3$ is as required in $\oplus$.
\end{enumerate}
\mn
Obviously this clearly suffices.  Clearly clause $\oplus(a)$ holds and
 clauses $\oplus(b)+(c)$ follows from clause $\oplus(d)$.

[Why?  \Wilog \, $G_1 = \bold K^{\lf}_\lambda$, let $\langle
a_\beta:\beta < \lambda\rangle$ list the elements of $G_1$.  For $c
\in G_3$ let $\bar a_c = \langle a_\beta:\beta < \lambda\rangle \char
94 \langle c \rangle$ and applying clause (d) we get $\cP_c \subseteq
\cP(\lambda+1)$ as there.  We finish letting 
$X := \{c \in G_3:\lambda \notin \cP_c\}$.]

Now let us prove clause $\oplus(d)$, so let $\alpha < \lambda^+,\bar a \in
{}^\alpha(G_3)$ and $p(\bar x_{[\alpha]}) = \tp_{\at}(\bar
a,G_1,G_3)$ and $p'(\bar x_{[\alpha]}) = \tp_{\bs}(\bar a,G_1,G_3)$; 
\wilog \, $\bar a$ is without repetitions but this is not used.

Define:
\mn
\begin{enumerate}
\item[$(*)_3$]  $\cP = \{u \subseteq \alpha$: there are $\bar a',\bar a'' 
\in {}^\alpha(G_3)$ realizing $p(\bar x_{[\alpha]})$ such that $u
= (\forall \beta < \alpha)(\beta \in u \equiv a'_\beta = a''_\beta)\}$.
\end{enumerate}
\mn
Now
\mn
\begin{enumerate}
\item[$(*)_4$]  $\alpha \in \cP$.
\end{enumerate}
\mn
[Why?  Let $\bar a' = \bar a'' = \bar a$.]
\mn
\begin{enumerate}
\item[$(*)_5$]  if $u_1,u_2 \in \cP$, then $u_1 \cap u_2 \in \cP$.
\end{enumerate}
\mn
[Why?  Let $\bar a'_\ell,\bar a''_\ell$ witness that $u_\ell \in \cP$, i.e. 
both $\bar a'_\ell,\bar a''_\ell$ realize $p(\bar x_{[\alpha]})$ in $G_3$ 
and $u_\ell = \{\beta < \alpha:a'_{\ell,\beta} = a''_{\ell,\beta}\}$.

Let $I = I_* + \sum\limits_{\varp < \mu} I_\varp$ be 
linear orders (so pairwise disjoint), where we chose the linear orders
such that $I_\varp \cong \alpha$ for $\varp < \mu$ and let
$s_{\varp,\beta}$ be the $\beta$-th member of $I_\varp$
and let $\langle c_s:s \in I_*\rangle$ list $G_3$ such that $c_{s(*)}
= e_{G_3}$ and $s(*) \in I_*$.

We shall now apply \ref{n67}, so let
\mn
\begin{enumerate}
\item[$(a)$]  $\gamma_* = 1 + \alpha + \alpha$ and $n_*=1$
\sn
\item[$(b)$]  for $\varp < \mu,\gamma < \gamma_*$ let 
$\langle h_{\gamma,0}(s_{\varp,\beta}):\beta < \alpha\rangle$ be equal to:
\sn
\begin{enumerate}
\item[$\bullet$]  $\bar a$ \If \, $\gamma =0$;
\sn
\item[$\bullet$]  $\bar a'_1$ \If \, $\gamma \in \{1 + \zeta:\zeta <
  \varp\}$;
\sn
\item[$\bullet$]  $\bar a''_1$ \If \, $\gamma \in [1+\zeta:\zeta \in
  [\varp,\alpha)\}$;
\sn
\item[$\bullet$]  $\bar a'_2$ \If \, $\gamma \in \{1 + \alpha
  +\zeta:\zeta < \varp\}$;
\sn
\item[$\bullet$]  $\bar a''_2$ \If \, $\gamma \in \{1 + \alpha +
  \zeta:\zeta \in [\varp,\alpha)\}$;
\end{enumerate}
\sn
\item[$(c)$]  $h_{\gamma,0}(s) = c_s$ for $s \in I,\gamma < \gamma_*$;
\sn
\item[$(d)$]  $G_3,G_3,I,I_*$ here stand for $G_0,H_*,I,I_*$ there.
\end{enumerate}
\mn
We get $(H,\bar a^*)$ as there, so by (B)(d) there essentially $G_3
\subseteq H$ and by (B)(c) there the $\bar a^* \rest I_\varp$
 realizes $p(\bar x_{[\alpha]})$; moreover, realizes $p'(\bar
 x_{[\alpha]})$; also $\langle \bar a^* \rest I_\varp:\varp <
 \mu\rangle$ is a $\Delta$-system with heart $u$.

The rest should be clear; we do not need to extend $G_3$ by $(*)_1$.]
\end{PROOF}
\newpage

\section {$\bold K_{\lf}$} \label{3}
\bigskip

\subsection {Indecomposability} \label{3a}
\bigskip

In this section we deal with indecomposability, equivalently $\CF(M)$,
see e.g. \cite{ShTh:524}; which is meaningful and of
interest also for other classes.

\begin{definition}
\label{b8}
1) We say $\theta \in \CF(M)$ or $M$ is $\theta$-decomposable 
\when \,: $\theta$ is regular
and if $\langle M_i:i < \theta\rangle$ is $\subseteq$-increasing with
union $M$, \then \, $M = M_i$ for some $i$.

\noindent
2) We say $M$ is $\Theta$-indecomposable \when \, it is
$\theta$-indecomposable for every $\theta \in \Theta$.

\noindent
3) We say $M$ is ($\ne \theta)$-indecomposable \when \,: $\theta$ is
regular and if $\sigma = \cf(\sigma) \ne \theta$ then $M$ is
$\sigma$-indecomposable. 

\noindent
4) We say $\bold c:[\lambda]^2 \rightarrow \chi$ is
$\theta$-indecomposable \when \,: if $\langle u_i:i < \theta\rangle$ is
$\subseteq$-increasing with union $\lambda$ then $\chi = 
\{\bold c\{\alpha,\beta\}:\alpha \ne \beta \in u_i\}$ for some $i <
\theta$; similarly for the
other variants.

\noindent
5) If we replace $\subseteq$ by $\le_{\gk},\gk$ an a.e.c., \then \, we write
$\CF_{\gk}(M)$ or ``$\theta-\gk$-indecomposable".
\end{definition}

\begin{definition}
\label{b9}
We say $G$ is $\theta$-indecomposable inside $G^+$ \when \,:
\mn
\begin{enumerate}
\item[(a)]  $\theta = \cf(\theta)$;
\sn
\item[(b)]  $G \subseteq G^+$;
\sn
\item[(c)]  if $\langle G_i:i \le \theta\rangle$ is
  $\subseteq$-increasing continuous and $G \subseteq G_\theta \subseteq
  G^+$ then for some $i < \theta$ we have $G \subseteq G_i$.
\end{enumerate}
\end{definition}

\begin{discussion}
\label{b11}
1) The point of the definition of indecomposable is the following
observation, \ref{b14}.

\noindent
2) Using cases of indecomposability, see \ref{b16}, help here to prove
density of complete members of $\bold K^{\lf}_\lambda$
 and improve characterization of the existence of universal members in
 e.g. cardinality $\beth_\omega$.
\end{discussion}

\begin{observation}
\label{b14}
1) If $\langle M_i:i < \delta\rangle$ is $\le_{\gk}$-increasing with
union $M$, each $M_i$ is $\theta-\gk$-indecomposable and $\cf(\delta) \ne
\theta$, \then \, $M$ is $\theta-\gk$-indecomposable.

\noindent
2) If $\langle M^\ell_i:i <\theta\rangle$ is $\le_{\gk}$-increasing
and for $\ell=1,2, \, 
\bigcup\limits_{i} M^1_i = M = \bigcup\limits_{i} M^2_i$ and each
$M^1_i$ is $\theta-\gk$-indecomposable, \then \, $\bigwedge\limits_{i <
  \theta} \, \bigvee\limits_{j < \theta} M^1_i \le_{\gk} M^2_j$.

\noindent
3) If $\langle M^\ell_i:i \le \delta\rangle$ is $\le_{\gk}$-increasing
continuous and each $M^\ell_{i+1}$ is $\theta-\gk$-indecomposable for $i <
\delta,\ell=1,2$ and $M^1_\delta = M^2_\delta$ and $\theta =
\cf(\delta) >\aleph_0$, \then \, $\{i < \delta:M^1_i = M^2_i\}$ is a
club of $\delta$.

\noindent
4) If $M$ is a Jonsson algebra of cardinality $\lambda$, \then \, $M$
is $(\ne \cf(\lambda))$-indecomposable.

\noindent
5) Assume $J$ is a directed partial order, $\langle M_s:s \in J\rangle$
is $\subseteq$-increasing and $J_* := \{s \in J:M_s$ is
$\theta-\gt$-indecomposable$\}$ is cofinal in $J$.  \Then \,
$\bigcup\limits_{s \in J} M_s$ is $\theta-\gk$-indecomposable provided
that:
\mn
\begin{enumerate}
\item[$(*)$]  if $\bigcup\limits_{i <\theta} J_i \subseteq J$ is
  cofinal and $\langle J_i:i < \theta\rangle$ is
  $\subseteq$-increasing, \then \, for
some $i,J_i$ is cofinal in $J$ or at least $\bigcup\limits_{s \in J_i}
M_s = \bigcup\limits_{s \in J} M_s$.
\end{enumerate}
\end{observation}

\begin{PROOF}{\ref{b14}}
Should be clear but we elaborate, e.g.:

\noindent
5) Toward contradiction let $\langle N_i:i < \kappa\rangle$ be
$\subseteq$-increasing with union $\bigcup\limits_{s \in J} M_s$.  For
each $s \in J_*$ there is $i(s) < \kappa$ such that $N_{i(s)}
\supseteq M_s$.  Let $J_i = \{t(s):s \in J_*$ and $i(s) \le i\}$ for $i <
\theta$.  Clearly $\langle J_i:i < \theta\rangle$ is as required in
the assumption of $(*)$, hence for some $i < \theta$ we have
$\bigcup\limits_{s \in J} M_i = \bigcup\limits_{s \in J_i} M_s$, 
so necessarily $N_i \supseteq \cup\{M_s:s \in J\}$, and thus equality holds.
\end{PROOF}

\noindent
We turn to $\bold K_{\lf}$.
\begin{claim}
\label{b15}
1) Assume $I$ is a linear order and $\bold c:[I]^2 \rightarrow \theta$ is
$\theta$-indecomposable, $G_1 \in \bold K_{\lf}$ and $a_i \in
G_1(i < \theta)$ are\footnote{The demand ``the $a_i$'s commute in
  $G_1$" is used in the proof of $(*)_8$, and the demand ``$a_{\beta_i}$ has
order 2" is used in the proof of $(*)_7$.} 
pairwise commuting and each of order 2.

\Then \, there is $G_2$ such that:
\mn
\begin{enumerate}
\item[(a)]  $G_2 \in \bold K_{\lf}$ extends $G_1$;
\sn
\item[(b)]  $G_2$ is generated by $G_1 \cup \bar b$ where $\bar b =  
\langle b_s:s \in I\rangle$;
\sn
\item[(c)]  $b_s$ commutes with $G_1$ and has order 2 for $s \in I$
\sn
\item[(d)]  if $s_1 \ne s_2$ are from $I$ then 
$[b_{s_1},b_{s_2}] = a_{\bold c\{s_1,s_2\}}$;
\sn
\item[(e)]  $G_2$ is generated by $G_1 \cup \bar b$ freely except the
  equations implicit in (a)-(d) above;
\sn
\item[(f)]  $\seb(\{a_i:i < \theta\},G_1)$ is $\theta$-indecomposable
  inside $G_2$; see Definition \ref{b9}.
\end{enumerate}
\mn
2) Assume $I$ a linear order is the disjoint union of $\langle
I_\alpha:\alpha < \alpha_*\rangle,u_\alpha \subseteq \Ord$ has
cardinality $\theta_\alpha$ and $\bold c_\alpha:[I_\alpha]^2
\rightarrow u_\alpha$ is $\theta_\alpha$-indecomposable for $\alpha <
\alpha_*,\langle u_\alpha:\alpha < \alpha_*\rangle$ is a sequence of
pairwise disjoint sets with union $u$ and $0 \notin u$ and $a_\varp \in
G_1$ and for $a_\varp,a_\zeta$ commute for $\varp,\zeta \in
u_\alpha,\alpha < \alpha_*$ and each $a_\varp$ has order 2, and $a_0 =
e$.

Let $\bold c:[I]^2 \rightarrow u \cup \{0\}$ extends each $\bold
c_\alpha$ and is zero otherwise.

\Then \, there is $G_2$ such that:
\mn
\begin{enumerate}
\item[(a)-(e)]  as above
\sn
\item[(f)]  if $\alpha < \alpha_*$ then $\seb(\{a_{\alpha,\varp}:\varp
  < \theta_\alpha\},\theta_1)$ is $\theta_\alpha$-indecomposable
  inside $G_2$. 
\end{enumerate}
\end{claim}

\begin{PROOF}{\ref{b15}}
For notational simplicity for $f$ part (2) \wilog \,:
\mn
\begin{enumerate}
\item[$(*)_0$]  each $I_\alpha$ is a convex subset of $I$.
\end{enumerate}
\mn
Let $G_2$ be the group generated freely by $G_1 \cup \bar b$ modulo
the equations in $\Gamma$ where:
\mn
\begin{enumerate}
\item[$(*)_1$]  $\Gamma = \bigcup\limits_{\ell < 4} \Gamma_\ell$
  where
\sn
\begin{enumerate}
\item[(a)]  $\Gamma_0$ is the set of equations $G_1$ satisfies
\sn
\item[(b)]  $\Gamma_1 = \{b_s b_s = e:s \in I\}$
\sn
\item[(c)]  $\Gamma_2 = \{b^{-1}_s ab_s = a:a \in G_1$ and $s \in I\}$
\sn
\item[(d)]  $\Gamma_3 = \{[b_s,b_t] = a_{\bold c\{s,t\}}:s <_I t\}$.
\end{enumerate}
\end{enumerate}
\mn
Next,
\mn
\begin{enumerate}
\item[$(*)_2$]
\begin{enumerate}
\item[(a)]  we define the relation $<$ on $[I]^{< \aleph_0}$ by $u_1
  <_I u_2 \Leftrightarrow (\forall s_1 \in u_1)(\forall s_2 \in
  u_1)(s_1 <_I s_2)$; we may replace $u_\ell$ by $s_\ell$ if $u_\ell =
  \{s_\ell\}$; note that by our choice, $u \in [I]^{< \aleph_0}
  \Rightarrow \emptyset < u < \emptyset$
\sn
\item[(b)]   if $u = \{s_0 > \ldots > s_{n-1}\},s_{n-1} > s$ let
  $a_{s,u} = a_{\bold c\{s,s_0\}} a_{\bold c\{s,s_1\}} \ldots
a_{\bold c\{s,s_{n-1}\}}$ so of order 2 because $\bold c\{s,s_\ell\} <
\theta$ so $a_{\bold c\{s,s_\ell\}} \in A := \{a_i:i < \theta\}$.
\end{enumerate}
\end{enumerate}
\mn
Let
\mn
\begin{enumerate}
\item[$(*)_3$]  $\cX = \{(u,a):u \subseteq I$ is finite and 
$a \in G_1\}$.
\end{enumerate}
\mn
Now we define a function $h$ from the set $Y = G_1 \cup \bar b$ of
generators of $G_2$ into $\Sym(\cX)$ as follows, we may write $h_x$
instead of $h(x)$.

First,
\mn
\begin{enumerate}
\item[$(*)_4$]  for $c \in G_1$ we choose $h_c \in \Sym(\cX)$ as
  follows: for $a \in G_1$ let $h_c(u,a) = (u,c a)$.
\end{enumerate}
\mn
Now clearly,
\mn
\begin{enumerate}
\item[$(*)_5$]
\begin{enumerate}
\item[(a)]   $h_c \in \Sym(\cX)$ for $c \in G_1$;
\sn
\item[(b)]   $c \mapsto h_c$ is an embedding of $G_1$ into
  $\Sym(\cX)$.
\end{enumerate}
\end{enumerate}
\mn
Next
\mn
\begin{enumerate}
\item[$(*)_6$]  for $t \in I$ we define $h_t:\cX \rightarrow \cX$ by
  defining $h_t(u,a)$ by induction on $|u|$ for $(u,a) \in \cX$ as
  follows:
\sn
\begin{enumerate}
\item[(a)]   if $u = \emptyset$ then $h_t(u,a) = (\{t\},a)$
\sn
\item[(b)]  if $u = \{s\}$ then $h_t(u,a)$ is defined as follows:
\begin{enumerate}
\item[$(\alpha)$]   if $t <_I s$ then $h_t(u,a) = (\{t,s\},a)$
\sn
\item[$(\beta)$]   if $t=s$ then $h_t(u,a) = (\emptyset,a)$
\sn
\item[$(\gamma)$]   if $s <_I t$ then $h_t(u,a) = (\{s,t\},
a_{\bold c\{s,t\}} a)$
\end{enumerate}
\sn
\item[(c)]  if $s_1 < \ldots < s_n$ list $u \in [I]^n$
  and $k \in \{0,\ldots,n\}$ and $s \in (s_k,s_{k+1})_I$ where we
stipulate $s_0 = - \infty,s_{n+1} = + \infty$ \then \, 
\newline
$h_t(u,a) = (u \cup \{t\},a_{\bold c\{s_0,t)} a_{\bold c\{s_1,t\}} \ldots 
a_{\bold c\{s_k,t\}} a)$
\sn
\item[(d)]  \If \, $s_1 < \ldots < s_n$ list $u \in [I]^n$
 and $k \in \{0,\ldots,n-1\}$ and $t = s_{k+1}$ \then \,
$h_t(u,a) = (u \backslash \{t\},a_{\bold c\{s_1,t\}} 
a_{\bold c\{s_2,t\}} \ldots a_{\bold c\{s_k,t\}} a)$.
\end{enumerate}
\end{enumerate}
\mn
Note that $(*)_6(b)(\alpha)$ is the same as $(*)_6(c)$ for $n=1$ and
$(*)_6(b)(\gamma)$ is the same as $(*)_6(d)$ for $n=1$.

Now clearly
\mn
\begin{enumerate}
\item[$(*)_7$]  if $t \in I$ then $h_t \in \Symc(\cX)$ has order 2.
\end{enumerate}
\mn
[It is enough to prove $h_t(h_t(\eta,a)) = (\eta,a)$.  We divide to
cases according to ``by which clause of $(*)_6$ is $h_t(\eta,a)$ defined".
\medskip

\noindent
\underline{If the definition} is by $(*)_6(a)$ then $h_t(\emptyset,a)
= (\{t\},a)$ and by $(*)_6(b)(\beta)$

\[
h_t h_t(\emptyset,a) = h_t(\{t\},a) = (\emptyset,a).
\]

\mn
\underline{If the definition} is by $(*)_6(b)(\beta)$, the proof is
similar.
\medskip

\noindent
\underline{If the definition} is by $(*)_6(b)(\gamma)$ then recalling
$(*)_6(d)$

\[
h_t(h_t(u,a)) =  h_t(h_t(\{s\},a)) = h_t(\{s,t\},a_{\bold
  c\{s,t\}} a) = (\{s\},a_{\bold c\{s,t\}} a_{\bold c\{s,t\}} a).
\]

\mn
But in the claim we assume $a_{\bold c\{s,t\}}$
 and every $a_i$ have order 2.
\medskip

\noindent
\underline{If the definition} is by $(*)(b)(\alpha)$, the proof is similar.

\noindent
\underline{If the definition} is by $(*)_6(c)$, then recall
$(*)_6(d)$ and compute similarly to the two previous cases, recalling
$\langle a_{\bold c\{s,t\}}:s \in I\rangle$ are pairwise commuting of order 2.
\medskip

\noindent
\underline{If the definition} is by $(*)_6(d)$ - this is just like 
the last case.]
\mn
\begin{enumerate}
\item[$(*)_8$]  $[h_{s_1},h_{s_2}] = h_{a_{\bold c\{s_1,s_2\}}}$ in
  $G_2$ for $s_1,s_2 \in I$.
\end{enumerate}
\mn
[Why?  We have to check by cases; here we use ``the $a_\varp$'s are
pairwise commuting for $s \in u_\alpha$" and the choice of $\bold c'$.]
\mn
\begin{enumerate}
\item[$(*)_9$]  if $c \in G_1$ and $t \in I$ then $h_c,h_t$ commute in
  $\Sym(\cX)$.
\end{enumerate}
\mn
[Why?  Just read $(*)_4 + (*)_6$.]
\mn
\begin{enumerate}
\item[$(*)_{10}$]
\begin{enumerate}
\item[(a)]  the mapping $c \mapsto h_c$ from $G_1 \cup \{b_s:s \in
  I\}$ respects all the equations from $\Gamma$
\sn
\item[(b)]  there is a homomorphism $\bold h$ from $G_2$ into
  $\Sym(\cX)$ such that $c \in G_1 \cup \bar b \Rightarrow \bold h(c)
  = h_c$.
\end{enumerate}
\end{enumerate}
\mn
[Why?  For clause (a), recalling $(*)_1$, 
for $\Gamma_0$ by $(*)_5(b)$, for $\Gamma_1$ by $(*)_7$, 
for $\Gamma_2$ by $(*)_9$ and for $\Gamma_3$ by $(*)_8$.
Clause (b) follows by clause (a).]
\mn
\begin{enumerate}
\item[$(*)_{11}$]   $G_1 \subseteq G_2$.
\end{enumerate}
\mn
[Why?  By $(*)_{10}(b)$ and $(*)_5(b)$.]
\mn
\begin{enumerate}
\item[$(*)_{12}$]  $\seb(\{a_i:i < \theta\},G_1)$ is
  $\theta$-indecomposable inside $G_2$.
\end{enumerate}
\mn
[Why?  Because the function $\bold c$ is $\theta$-indecomposable by assumption.]
\mn
\begin{enumerate}
\item[$(*)_{13}$]  $H$ is locally finite.
\end{enumerate}
\mn
[Why?  Check or see the proof of \ref{b16}.]

Together we are done.
\end{PROOF}

\noindent
The following will not be used (except in the proof of \ref{b15}).
\begin{claim}  
\label{b16}
If (A) \then \, there is a pair $(H,\bar b)$ such that (B) holds, where:
\mn
\begin{enumerate}
\item[(A)]
\begin{enumerate}
\item[(a)]  $\bold c:[\lambda]^2 \rightarrow \lambda$
 is $\Theta$-indecomposable, (e.g. $\Theta = (\ne \cf(\lambda)))$
and satisfies $\bold c(\{\alpha,\beta\}) \le \alpha \cap \beta$;
\sn
\item[(b)] $G \in \bold K^{\lf}_{\le \lambda}$;
\sn
\item[(c)]  $\bar a = \langle a_\alpha:\alpha < \lambda \rangle$
  is a sequence of elements of $G$ possibly with repetitions, and
  (used only for $(B)(c),(d)(\alpha)$) it generates $G$;
\sn
\item[(d)]  the $a_\alpha$'s are pairwise commuting and of order 2
\end{enumerate}
\sn
\item[(B)]
\begin{enumerate}
\item[(a)]  $H \in \bold K^{\lf}_\lambda$;
\sn
\item[(b)]  $G \subseteq H$;
\sn
\item[(c)]  
\begin{enumerate}
\item[$(\alpha)$]  $\seb(\{a_\alpha:\alpha < \lambda\},G)$ is 
$\Theta$-indecomposable inside $H$;
\sn
\item[$(\beta)$]   if $\bold c$ witness $\lambda 
\nrightarrow [\lambda]^2_\lambda$ \then \, $H$ is a Jonsson group;
\end{enumerate}
\sn
\item[(d)]
\begin{enumerate}
\item[$(\alpha)$]  $\bar b \in {}^\lambda H$ generates  $H$ over $G$;
\sn
\item[$(\beta)$]  each $b_\alpha$ has order 2;
\sn
\item[$(\gamma)$]  each $b_\alpha$ commutes with $G$;
\end{enumerate}
\sn
\item[(e)]  if $\alpha < \beta < \lambda$ and $\bold
  c(\{\alpha,\beta\}) = 2 \gamma$ then $[b_\alpha,b_\beta] = a_\gamma$;
\sn
\item[(f)]  if $\alpha < \beta < \lambda$ and 
$\bold c\{\alpha,\beta\} = 2 \gamma +1$ then $[b_\alpha,b_\beta] = b_\gamma$.
\end{enumerate}
\end{enumerate}
\mn
2) Like part (1), making the changes as in \ref{b15}.
\end{claim}

\begin{PROOF}{\ref{b16}}
For notational simplicity we stipulate $\bold c(\{\alpha,\alpha\}) =
\alpha$ for $\alpha < \lambda$.
Let $H$ be the group generated by $H \cup \{b_\alpha:\alpha < \lambda\}$ freely
except the set $\Gamma$ of equations: those stated in
$(B)(b),(d)(\beta),(\gamma)(e),(f)$.

Clearly,
\mn
\begin{enumerate}
\item[$(*)_1$] 
\begin{enumerate}
\item[$(a)$]   every element of $H$ can be 
represented as a product of the 
form $b_{\alpha_0} \ldots b_{\alpha_{n-1}} a$ 
with $\alpha_0 > \ldots > \alpha_{n-1}$ all $< \lambda$ and $a \in G$;
\sn
\item[$(b)$]  there is a unique homomorphism $\bold h$
  from $H$ onto $H_\lambda := \bigoplus\limits_{\alpha < \lambda}
  (\bbZ/z \bbZ) x_\alpha$ mapping $b_\alpha$ to $x_\alpha$ and every
  $a \in G$ to the unit element
\sn
\item[$(c)$]   moreover, the representation in clause (a) is unique.
\end{enumerate}
\end{enumerate}
\mn
[Why?  There is such a representation by the equations from
$(B)(d)(\beta),(\gamma)$ so clause (a) holds.  Considering the choice of $H$,
obviously clause $(*)_1(b)$ holds.  Lastly, for clause $(*)_1(c)$ assume
$b_{\alpha_{1,0}} \ldots b_{\alpha_{1,n_1-1}} a_1 = c =
b_{\alpha_{2,0}} \ldots b_{2,n_2-1} a_2$, both as in (a); using $h$
necessarily $n_1 =n_2,\alpha_{1,\ell} = \alpha_{2,\ell}$ (for
$\ell=0,\dotsc,n_i-1)$ and we can finish easily as in the proof of \ref{b15}.]

Hence,
\mn
\begin{enumerate}
\item[$(*)_2$]  $G \subseteq H$;
\sn
\item[$(*)_3$]  $H$ has cardinality $\lambda$.
\end{enumerate}
\mn
So clause (B)(b) holds and clauses (B)(c)-(f) are easy to check.  But what
about clause (a), i.e. why is $H$ locally finite?

Clearly for any finite subset $A$ of $H$ we can find
a finite $u \subseteq \lambda$ and finitely generated subgroup $K$ of $G$ such
that:
\mn
\begin{enumerate}
\item[$(*)_4$]  $A$ is included in the subgroup $L$ of $H$
  generated by $K \cup \{b_\alpha:\alpha \in u\}$.
\end{enumerate}
\mn
By clause (A)(a) of the claim's assumption, 
in particular, $\bold c(\{\alpha,\beta\}) \le \alpha
\cap \beta$, there is a finite $u_1$ such that $u \subseteq u_1
\subseteq \lambda$ and $\alpha,\beta \in u_1 \Rightarrow \bold
c(\{\alpha,\beta\}) \in u_1$; why? e.g. prove by induction on $\gamma$
that: \Iff \, $u \subseteq \lambda$ is finite and $\alpha,\beta \in u
\backslash \gamma \Rightarrow \bold c\{\alpha,\beta\} \in u$, \then \,
there is a finite $u_1$ such that $u \subseteq u_1 \subseteq
\lambda,u_1 \backslash \gamma = u \backslash \gamma$ then
$\alpha,\beta \in u_1 \Rightarrow \bold c(\{\alpha,\beta\})$.

Let $K_1$ be the subgroup of $G$ generated by $K \cup\{a_\gamma$: for
some $\alpha \ne \beta \in u_1$ we have $2 \gamma = 
\bold c(\{\alpha,\beta\})\}$.

Now,
\mn
\begin{enumerate}
\item[$(*)_5$]
\begin{enumerate}
\item[(a)]  $K_1$ is f.g. $\subseteq G$ and hence it is finite;
\sn
\item[(b)]  $L_1= \{b_{\alpha_0} \ldots b_{\alpha_{n-1}} a:a \in K_1,n
  < \omega$ and $\alpha_0 > \ldots > \alpha_{n-1}$ belongs to 
$u_1\}$ is a subgroup of $G$ which has $\le |K_1| \cdot 2^{|u_1|} 
< \aleph_0$ members;
\sn
\item[(c)]   if $\alpha,\beta \in u_1$ and $[b_\alpha,b_\beta] =
  b_\gamma$ then $\gamma \in u_1$;
\sn
\item[(d)]  if $\alpha,\beta \in u_1$ and 
$[b_\alpha,b_\beta] = a_\gamma \in G$ then $a_\gamma \in K_1$;
\sn
\item[(e)]  $L_1$ is a subgroup of $H$;
\sn
\item[(f)]  $L_1$ is finite.
\end{enumerate}
\end{enumerate}
\mn
Hence,
\mn
\begin{enumerate}
\item[$(*)_6$]  $L_1$ is a finite subgroup of $H$ containing $A$.
\end{enumerate}
\mn
As $A$ was any finite subset of $H$, clearly $H$ is locally finite.
\end{PROOF}

\begin{claim}
\label{b18}
If (A) and (B) then (C) \when \,:
\mn
\begin{enumerate}
\item[(A)]  
\begin{enumerate}
\item[(a)]  $\lambda \le \mu$
\sn
\item[(b)]  $\theta = \cf(\theta) < \lambda$
\sn
\item[(c)]  $\Theta$ is a set of regular cardinals $\le \mu$
such that at least one of the following holds:
\begin{enumerate}
\item[$(\alpha)$]  $\Theta = \{\theta\}$ and $\lambda =
  \lambda^\theta$ and some $\bold c:[\lambda]^2 \rightarrow \theta$ is
  $\theta$-indcomposable;
\sn
\item[$(\beta)$]  for every $\theta \in \Theta$ there is $\partial \in
  (\theta,\lambda]$ such that $\lambda^{\langle
  \partial;\theta\rangle} = \lambda$ and\footnote{This formulation is
  to stress that it would be better if we can omit ``$\partial =
  \theta$" but this requires strengthening \ref{b15} which we do not
  know how to do for $\bold K_{\lf}$.}
$\partial = \theta$ and some $\bold c_\theta:[\lambda]^2 \rightarrow \theta$ is
  $\theta$-indecomposable, see Definition \ref{w22}(c), \ref{b8}(4)
  respectively;
\end{enumerate}
\sn
\item[(d)]  $\theta_* = \cf(\theta_*)$ is $\le \mu$ but $\notin \Theta$.
\end{enumerate}
\sn
\item[(B)]  $(G_1,G^+_1) \in \bold K^{\lf}_{\lambda,\mu}$ which means:
\sn
\begin{enumerate}
\item[(a)]  $G_1 \subseteq G^+_1 \in \bold K_{\lf}$;
\sn
\item[(b)]  $G_1,G^+_1$ is of cardinality $\lambda,\mu$ respectively
\end{enumerate}
\sn
\item[(C)]  there are $G_2,G^+_2$ such that:
\sn
\begin{enumerate}
\item[(a)]  $(G_2,G^+_2) \in \bold K^{\lf}_{\lambda,\mu}$ and even
  $(G_2,G^+_2) \in \bold K^{\exlf}_{\lambda,\mu}$;
\sn
\item[(b)]  $(G_1,G^+_1) \subseteq (G_2,G^+_2)$, i.e. $G_1 \subseteq
  G_2$ and $G^+_1 \subseteq G^+_2$;
\sn
\item[(c)]  $G_1$ is $\Theta$-indecomposable inside $G_2$;
\end{enumerate}
\end{enumerate}
\end{claim}

\begin{remark}
The exact construction inside the proof 
probably will be used in \S4, but the present version of the claim 
is enough for \S(3B).
\end{remark}

\begin{PROOF}{\ref{b18}}
\smallskip

\noindent
\underline{Stage A}:  If (A)(c)$(\alpha)$ holds then easily $\partial =
\theta^+$ is as required in (A)(c)$(\beta)$ so \wilog \,
(A)(c)$(\beta)$ holds.
\mn
\begin{enumerate}
\item[$(*)_1$]  It suffices to prove the following apparently weaker
  version of (C):
\sn
\begin{enumerate}
\item[(C)$'$]  assume in addition that
$(G_1,G^+_1) \in \bold K^{\exlf}_{\lambda,\mu}$;
  for each $\theta \in \Theta$ there are $G_2,G^+_2$ such that
  (a),(b) there hold and
\mn
\begin{enumerate}
\item[(c)$'$]  $G_1$ is $\theta$-indecomposable inside $G_2$.
\end{enumerate} 
\end{enumerate} 
\end{enumerate} 
\mn
[Why?  Let $\Theta = \{\theta_i:i < i_*\}$ and \wilog \, $i^*$ is a
limit ordinal $\le \mu$ (so possibly with repetitions and 
$2i_* = i_*$) and we choose
$(G_{2,j},G^+_{2,j})$ by induction on $j \le \theta_*(2 i_*)$ such that:
\mn
\begin{itemize}
\item  $(G_{2,j},G^+_{2,j}) \in \bold K^{\lf}_{\lambda,\mu}$ is increasing
  continuous;
\sn
\item  $(G_{2,0},G^+_{2,0}) = (G_1,G^+_1)$;
\sn
\item  if $j = i_* \iota + 2i+2$ and $i < i_*,\iota < \theta_*$ then $G_{2,j}$
is $\theta_*$-indecomposable inside $G^+$;
\sn
\item  if $j = i_* \iota + 2 i+1$ and $i < i_*,\iota < \theta_*$, 
then $(G_{2,j+1},G^+_{2,j+1}) \in \bold K^{\exlf}_{\lambda,\mu}$.
\end{itemize}
\mn
Clearly possible and by \ref{b14}(1) the pair 
$(G_{2,i_*,\theta},G^+_{2,i_*,\theta})$ is as
required.]

So we can fix $\theta \in \Theta,(G_1,G^+_1) \in \bold
K^{\exlf}_{\lambda,\mu}$ and let $\partial \le \lambda$ be such that
$\lambda^{(\partial;\theta)} = \lambda$ and $\bold c:\lambda
\rightarrow \theta$ is $\theta$-indecomposable and $\partial =
\theta$.

Now,
\mn
\begin{enumerate}
\item[$(*)_2$]   there is $(G_2,G^+_2) \in \bold K^{\lf}_{\lambda,\mu}$ above
  $(G_1,G^+_1)$ and $a_\varp \in G_2$ for $\varp < \partial$
  satisfies $\langle \square_{a_\varp}(G_1):\varp < \partial\rangle$
  are pairwise commuting and $\seb(\cup\{\square_{a_\zeta}(G_1):\zeta
  \in \partial \backslash \{\varp\}\}) \cap \square_{a_\varp}(G_1) = \{e\}$.
\end{enumerate}
\mn
[Why? Using $\gs = \gs_{\cg}$, see
\cite[2.17=Lc50(1),2.18=Lc52(2)]{Sh:312}.  That is, we can find an
increasing continuous sequence $\langle G^+_{2,\varp}:\varp
\le \partial\rangle$ of members of $\bold K_{\exlf}$ (hence with
trivial center) and $a_s \in G^+_{2,\varp +1}$ realizing in
$G^+_{2,\varp +1}$ the type $q_{\cg}(\langle \rangle,G^+_{2,\varp})$,
hence $\square_{a_\varp}(G_1) = a^{-1}_\varp G^+_{2,\varp} a_\varp$
commute with $G^+_{2,\varp}$.  So the ``pairwise commuting" follows
and the demand on the intersection holds by $G_1$ being from 
$\bold K_{\exlf}$.]
\mn
\begin{enumerate}
\item[$(*)_3$]   there is $(G_3,G^+_3) \in \bold
  K^{\lf}_{\lambda,\mu}$ above $(G_2,G^+_2)$ such that:
\sn
\begin{itemize}
\item  if $\bar b = \langle b_i:i < \partial\rangle$ is a sequence of pairwise
  commuting elements of order 2, independent, \then \, for some 
$\bar c = \langle c_\gamma:\gamma < \lambda\rangle$ and
increasing sequence $\langle \varp(i):i < \theta \rangle$ of
ordinals $< \partial$ we have $\gamma_1 < \gamma_2 < \lambda \wedge
\bold c_\theta\{\gamma_1,\gamma_2\} = i \Rightarrow 
[c_{\gamma_1},c_{\gamma_2}] = b_{\varp(i)}$.
\end{itemize}
\end{enumerate}
\mn
[Why?  As $\lambda^{\langle \partial;\theta\rangle} = \lambda$, 
there is a sequence $\langle \bar d_\alpha:
\alpha < \lambda\rangle$ such that:
\mn
\begin{itemize}
\item  for each $\alpha,\bar d_\alpha = \langle d_{\alpha,i}:i 
< \theta\rangle$, where $d_{\alpha,i} \in G_2$ are pairwise commuting,
each of order 2 for $i < \theta$ and are independent;
\sn
\item  if $\langle b_\varp:\varp < \partial\rangle$ is as in $\bullet$
  of $(*)_3$ above so $b_\varp \in G_2$, \then \, 
for some $\alpha < \lambda$ and increasing
  sequence $\langle \varp(i):\varp < \theta\rangle$ of ordinals
  $< \partial$ we have $d_{\alpha,i} = b_{\varp(i)}$ for $i < \theta$.
\end{itemize}
\mn
Now we shall take care of each by \ref{b15}, i.e. we choose
$(G_{2,\alpha},G^+_{2,\alpha}) \in \bold K^{\lf}_{\lambda,\mu}$ by
induction on $\alpha < \lambda$,
increasing continuous with $\alpha$, such that
$\seb(\{d_{\alpha,i}:i < \theta\},G_2)$ is $\theta$-indecomposable
inside $G^+_{2,\alpha +1}$.

Lastly, let $(G_3,G^+_3) = (G_{2,\lambda},G^+_{2,\lambda})$; clearly
we are done proving $(*)_3$.]

It suffices to show that $(G_3,G^+_3)$ is as required, so toward
contradiction assume $\bar H = \langle H_i:i < \theta\rangle$ is 
increasing continuous with union $G^+_3$ and 
$\langle H_i \cap G_1:i < \theta\rangle$ is not eventually constant.

Now we shall use $(*)_2$.   So for each $\varp < \partial$, for some
$j_\varp < \theta,a_\varp \in H_{j_\varp}$ hence $\langle
\square_{a_\varp}(G_1) \cap H_i:i < \theta \rangle$ cannot be
eventually constant, hence there is $b_i \in \square_{a_i}(G_1)
\backslash H_i$ and we apply $(*)_3$ above to find $\langle
c_\alpha:\alpha < \lambda\rangle$ as there.  As $\bold c:[\lambda]^2
\rightarrow \theta$ is $\theta$-indecomposable we are done.
\end{PROOF}

\begin{discussion}
\label{b18d}
We discuss some variants of \S4 when we use \ref{b18} (instead of the
version with $\partial = \theta$).

\noindent
1) For the case $|G| = \mu,\lambda = \mu^+$ we can waive $S_2$, but we
still use $\lambda = \lambda^{(\theta;\sigma)},\sigma = \aleph_0$.

\noindent
2) For the case $\lambda = \mu$ having $\theta$ copies can be done in $i
\in S_1$, this is the first step in the proof of \ref{b18}, the
$\langle a_\varp:\varp < \partial = \theta\rangle$.  For the second
step, for $i \in S_1$, we have ``tasks": $\langle b_{i,\alpha +
  \theta}:\varp < \theta\rangle$ which is the basis of a vector space
over $\bbZ/2 \bbZ$, and we add a coy using the $\theta$-indecomposable
$\bold c:[\lambda]^2 \rightarrow \theta$.  So we use $\lambda =
\lambda^{\langle \partial;\theta\rangle}$; can we weaken this?

\noindent
3) Comparing (1) and (2) above, certainly $\lambda = \mu^+$ of (1) is
improved by $\lambda = \mu$ of (2), but there is a price: instead
$\lambda^{(\partial,\theta)}$ we need $\lambda^{\langle
  \theta;\theta\rangle} = \lambda$ which almost say $\lambda \ge
2^\theta$.  The former is essentially almost always true but not the
second.  Still, if $\lambda \ge \beth_\omega$, for every regular large
enough $\theta < \beth_\omega$, we have $\lambda = \lambda^{\langle
  \theta;\theta\rangle}$.  

\noindent
4) But for \S(3B) the present version of \ref{b18} is more than
sufficient.

\noindent
5) What occurs if we know only that $\lambda^{[\theta;\theta]}  =
\lambda$?  As we naturally assume 
$\lambda^\theta > \lambda$.

We have to deal with the following:
\mn
\begin{enumerate}
\item[$(*)$]
\begin{enumerate}
\item[(a)]  $G_1 \in \bold K^{\exlf}_\lambda$;
\sn
\item[(b)]  $a_i \in G$ is of order 2 for $i < \theta$;
\sn
\item[(c)]  we try to find $G_2 \in K^{\lf}_\lambda$ extending $G_1$
  such that: if $u \in [\theta]^\theta,\langle a_i:i \in u\rangle$ are
  pairwise commuting and independent, \then \, $\seb(\{a_i:i \in
  u\},G_i)$ is $\theta$-indecomposable inside $G_2$.
\end{enumerate}
\end{enumerate}
\end{discussion}

\begin{discussion}
\label{b18e}
So we wonder:
\mn
\begin{enumerate}
\item[$(*)$]  for which $\bold c:[\lambda]^k \rightarrow \lambda$ do
  we have:
\sn
\begin{itemize}
\item  for every $G \in \bold K^{\exlf}_\lambda,\bar a \in
  {}^\lambda(G_1)$ is there $G_2 \in \bold K^{\exlf}_\lambda$ such that:
\begin{enumerate}
\item[(a)]  $G_2$ generated by $G_1 \cup \{b_\alpha:\alpha < \lambda\}$;
\sn
\item[(b)]  $G_1 \subseteq G_2$
\sn
\item[(c)]  $b_\alpha$ of order 2;
\sn
\item[(d)]  $a_{\bold c(u)} \in \seb(\{b_\alpha:\alpha \in u\},G_1)$
  for $u \in [\lambda]^k$ such that $\bold c(u) \ne 0$.
\end{enumerate}
\end{itemize}
\end{enumerate}
\end{discussion}

\begin{conclusion}
\label{b19}
1) Assume $\lambda$ uncountable, $\theta =
\cf(\theta) \le \lambda$ and let $\Theta_\lambda =
\Theta_{\lambda,\theta} =
\{\cf(\lambda),\theta\}$ except when $\lambda = \mu^+,\mu >
\cf(\mu)$ in which case $\Theta_\lambda =
\{\cf(\lambda),\cf(\mu),\theta\}$ and $\Theta'_\lambda = \{\sigma <
\lambda:\sigma = \cf(\lambda)$ and $\lambda^{\langle
  \lambda;\sigma\rangle} = \lambda$.
If $G \in \bold K^{\lf}_\lambda$, 
\then \, there is a $(\Theta'_\lambda 
\backslash \Theta_\lambda)$-indecomposable 
$H \in \bold K^{\lf}_\lambda$ extending $G$.

\noindent
2) If $G \in \bold K^{\lf}_\lambda,\lambda > \aleph_0$, \then \, we can find
$\bar H = \langle H_i:i \le \cf(\lambda)\rangle$ increasing continuous
such that $G \subseteq H_{\cf(\lambda)},i < \cf(\lambda) \Rightarrow
|H_i| < \lambda$ and $H_i$ is $(\Theta'_{\|H_2\|},
\Theta_{\|H_i\|})$-indecomposable.
\end{conclusion}

\begin{PROOF}{\ref{b19}}
We prove this by induction on $\lambda$ using \ref{b15}, \ref{b18} 
and \ref{w14},
\ref{b14}(1); well for successor of singular we use also \ref{b14}(5) 
arriving to $\lambda$.

\noindent
1) For part (1) it suffices to find $\bold c:[\lambda]^2 \rightarrow
\lambda$ which is $\theta$-indecomposable for every $\theta = \cf(\theta)
\notin \Theta_\lambda$ (as in the proof of \ref{b18}).
\medskip

\noindent
\underline{Case 1}:  $\lambda = \mu^+,\mu$ regular

Use \ref{b18} with $(\lambda,\lambda,\lambda,\Theta)$ standing for 
$(\lambda,\mu,\theta_*,\Theta)$ there and \ref{w14}(1).
\medskip

\noindent
\underline{Case 2}:   $\lambda$ a limit cardinal

As in Case 1, we can find 
$\bold c:[\lambda]^2 \rightarrow \lambda$ such that if
$\mu < \lambda$ then $\bold c \rest [\mu^{++} \backslash \mu^+]^2
\rightarrow \mu^{++}$ witness $\mu^{++} \nrightarrow
[\mu^{++}]^2_{\mu^+}$.  Now check.
\medskip

\noindent
\underline{Case 3}:  $\lambda$ a successor of a singular cardinal

Let $\lambda = \mu^+,\mu >\cf(\mu)$

Let $\langle \lambda_i:i <
\cf(\mu)\rangle$ be increasing with limit $\mu$, for $\alpha \in
[\mu,\lambda)$ let $\bar u_\alpha  = \langle u_{\alpha,i}:i <
\cf(\mu)\rangle$ be increasing continuous with union $\alpha$ such
that $u_{\alpha,0} = \emptyset,|u_{\alpha,i+1}| = \lambda^{++}_i$.

We define $\bold c:[\lambda]^2 \rightarrow \lambda$ such that:
\mn
\begin{enumerate}
\item[$(*)$]  if $\alpha \in [\mu,\lambda)$ \then \,:
\sn
\begin{enumerate}
\item[(a)]  $\bold c$ maps $[\{\mu \alpha + \varp:\varp \in
u_{\alpha,i+1} \backslash u_{\alpha,i}\}]^2$ into $u_{\alpha,i}$;
\sn
\item[(b)]  if $u \subseteq \{\mu \alpha + \varp:\varp \in
u_{\alpha,i+1} \backslash u_{\alpha,i}\}$  has cardinality
$\lambda^{++}_i$ then $\Rang(\bold c \rest [u]^2) = u_{\alpha,i}$.
\end{enumerate}
\end{enumerate}
\mn
It is easy to find such a $\bold c$.  Now assume $\langle v_\varp:\varp
< \delta\rangle$ is $\subseteq$-increasing with union $\lambda$ and
$\cf(\delta) \notin \Theta_\lambda$.  For each $\alpha \in
[\mu,\lambda)$, for every large enough $i < \cf(\mu),\lambda^{++}_i
  \ne \cf(\delta)$, so as $|u_{\alpha,i+1} \backslash u_{\alpha,i}| =
  \lambda^{++}_i$, for some $\varp = \varp(\alpha,i) < \delta$ the set
  $(u_{\alpha,i+1} \backslash u_{\alpha,i}) \cap v_\varp$ has
cardinality $\lambda^{++}_i$, hence $\Rang(\bold c \rest 
[v_\zeta]^\varp) \supseteq \Rang(\bold c \rest [(u_{\alpha,i+1}
\backslash u_{\alpha,i}) \cap v_\varp]^2) \supseteq u_{\alpha,i}$.  

Lastly, for each $\alpha$ for some $\varp(\alpha) <
  \delta,\cf(\mu) = \sup\{i:\varp(\alpha,i) < \varp(\alpha)\}$ and
  for some $\zeta < \delta$ we have $\lambda = \mu^+ = \sup\{\alpha:\alpha \in
  [\mu,\lambda)$ and $\varp(\alpha) < \zeta\}$.

\noindent
2) Should be clear.
\end{PROOF}

\begin{remark}  
\label{b22}
Claim \ref{b19} implies $\bold K^{\exlf}_\lambda$ densely has 
$(\Theta_\lambda)$-indecomposable models and may help in:
\mn
\begin{enumerate}
\item[$(a)$]  strengthening \cite[5.1]{Sh:312}, see \S4,\S5
\sn
\item[$(b)$]  proving $\bold K^{\lf}_{\beth_\omega}$ has no universal
member and see \S(3B).
\end{enumerate}
\end{remark}
\bigskip

\subsection {Universality} \label{3b}\
\bigskip

For quite many classes, there are universal members in any (large enough)
$\mu$ which is strong limit of cofinality $\aleph_0$, see
\cite{Sh:820} which include history.   
Below we investigate ``is there a universal member of
$\bold K^{\lf}_\mu$ for such $\mu$".  We prove that if
there is a universal member, e.g. in $\bold K^{\lf}_\mu$, then there is a
canonical one.

\begin{theorem}
\label{b25}
1) Let $\mu$ be strong limit of cofinality $\aleph_0$.

The following conditions are equivalent:
\mn
\begin{enumerate}
\item[$(A)$]   there is a universal $G \in \bold K^{\lf}_\mu$
\sn
\item[$(B)$]   if $H \in \bold K^{\lf}_\lambda$ is 
$\aleph_0$-indecomposable (see existence in \ref{b18}(B) + \ref{b19}) 
for some $\lambda < \mu$, \then \, there is a sequence
$\bar G = \langle G_\alpha:\alpha < \alpha_* \le \mu\rangle$ such that:
\sn
\begin{enumerate}
\item[$(a)$]   $H \subseteq G_\alpha \in \bold K^{\lf}_\mu$
\sn
\item[$(b)$]   if $G \in \bold K^{\lf}_\mu$ extend $H$, \then \, for
  some $\alpha,G$ is embeddable into $G_\alpha$ over $H$.
\end{enumerate}
\sn
\item[$(B)^+$]  We can add in (B)
\sn
\begin{enumerate}
\item[$(c)$]   if $\alpha_1 < \alpha_2 < \alpha_*$, \then \, there are
  no $G,f_1,f_2$ such that $H \subseteq G \in \bold K_{\lf}$ and
  $f_\ell$ embeds $G_{\alpha_\ell}$ into $G$ over $H$ for $\ell=1,2$.
\sn
\item[$(d)$]   $(H,G_\alpha)$ is an amalgamation pair (see Definition
  \ref{a31}(4), i.e. with $\chi_1 = \chi_2 = \mu = \|G_\alpha\|$); 
moreover if $G_\alpha \subseteq G' \in \bold K^{\lf}_\mu$, \then \, $G'$
  is embeddable into $G_\alpha$ over $H$.
\end{enumerate}
\end{enumerate}
\mn
2)  We can add in part (1):
\mn
\begin{enumerate}
\item[$(C)$]   there is $G_*$ such that:
\sn
\begin{enumerate}
\item[$(a)$]  $G_* \in \bold K^{\lf}_\mu$ is universal in 
$\bold K^{\lf}_{<\mu}$;
\sn
\item[$(b)$]  $\cE^{\aleph_0}_{G_*,<\mu}$, see below, 
is an equivalence relation with $\le \mu$ equivalence classes;
\sn
\item[$(c)$]  $G_*$ is
  $\aleph_0-\cE^{\aleph_0}_{G_*,<\mu}$-indecomposably homogeneous, see
  below.
\end{enumerate}
\end{enumerate}
\end{theorem}

\noindent
Before we prove \ref{b25},
\begin{definition}
\label{b41}
For $\theta = \cf(\theta) < \mu$ and
$\gk = \bold K_{\lf}$ or any a.e.c. with $\LST_{\gk} < \mu$, we
define for $M_* \in K^{\gk}_\mu$:
\mn
\begin{enumerate}
\item[$(A)$]   $\IND^\theta_{M_*,<\mu} = \{N:N \le_{\gk} M_*$ has
  cardinality $< \mu$ and is $\theta$-indecomposable$\}$.
\sn
\item[$(B)$]   $\cF^\theta_{M_*,<\mu} = \{f$: for some
  $\theta$-indecomposable $N = N_f \in K^{\gk}_{< \mu}$ with universe
  an ordinal, $f$ is a $\le_{\gk}$-embedding of $N$ into $M_*\}$.
\sn
\item[$(C)$]   $\cE^\theta_{M_*,<\mu} = \{(f_1,f_2):f_1,f_2 \in
  \cF^\theta_{M_*,<\mu},N_{f_1} = N_{f_2}$ and there are
  $\le_{\gk}$-embeddings $g_1,g_2$ of $M_*$ into some
  $\le_{\gk}$-extension $M$ of $M_*$ such that $g_1 \circ f_1 = g_2
  \circ f_2\}$.
\sn
\item[$(D)$]   $M_*$ is $\theta-\cE^\theta_{M_*,<\mu}$-indecomposably
homogeneous \when \,: if $f_1,f_2 \in \cF^\theta_{M_*,<\mu}$  and
$(f_1,f_2) \in \cE^\theta_{M_*,<\mu}$ and $A \subseteq M_*$ has
  cardinality $< \mu$ \then \, there is $(g_1,g_2) \in
\cE^\theta_{M_*,<\mu}$ such that $f_1 \subseteq g_1 \wedge f_2
  \subseteq g_2$ and $A \subseteq \Rang(g_1) \cap \Rang(g_2)$; it
  follows that if $\cf(\mu) = \aleph_0$ then for some $g \in
  \aut(M_*)$ we have $f_2 = g \circ f_1$.
\end{enumerate}
\end{definition}

\begin{remark}
\label{b44}
We may consider in \ref{b25} also $(A)_0 \Rightarrow (A)$ where
\mn
\begin{enumerate}
\item[$(A)_0$]   if $\lambda < \mu,H \subseteq G_1 \in \bold
  K^{\lf}_{< \mu}$ and $|H| \le \lambda$, \then \, for some $G_2$ we
  have $G_1 \subseteq G_2 \in \bold K^{\lf}_{< \mu}$ and $(H,G_2)$ is a
  $(\mu,\mu,\lambda)$-amalgamation base.
\end{enumerate}
\end{remark}

\begin{PROOF}{\ref{b23}}
It suffices to prove the following implications:
\medskip

\noindent
\underline{$(A) \Rightarrow (B)$}:

Let $G_* \in \bold K^{\lf}_\mu$ be universal and choose a sequence
$\langle G^*_n:n < \omega\rangle$ such that $G_* =
\bigcup\limits_{n} G^*_n,G^*_n \subseteq G^*_{n+1},|G^*_n| < \mu$.

Let $H$ be as in \ref{b25}(B) and let
$\cG = \{g:g$ embed $H$ into $G^*_n$ for some $n\}$.  So clearly
$|\cG| \le \sum\limits_{n} |G^*_n|^{|H|} \le \sum\limits_{\lambda <
  \mu} 2^\lambda = \mu$ (an over-kill).

Let $\langle g^*_\alpha:\alpha < \alpha_* \le \mu\rangle$ list $\cG$
and let $(G_\alpha,g_\alpha)$ be such that:
\mn
\begin{enumerate}
\item[$(*)_1$]  $(a) \quad H \subseteq G_\alpha \in \bold K^{\lf}_\mu$;
\sn
\item[${{}}$]  $(b) \quad g_\alpha$ is an isomorphism from $G_\alpha$
  onto $G_*$ extending $g^*_\alpha$.
\end{enumerate}
\mn
It suffices to prove that $\bar G = \langle G_\alpha:\alpha <
\alpha_*\rangle$ is as required in (B).  Now clause (B)(a) holds by
$(*)_1(a)$ above.  As for clause (B)(b), let $G$ satisfy 
$H \subseteq G \in \bold
K^{\lf}_{\le \mu}$, hence there is an embedding $g$ of $G$ into
$G_*$.  We know $g(H) \subseteq G = \bigcup\limits_{n} G_n$ hence $\langle
g(H) \cap G_n:n < \omega\rangle$ is $\subseteq$-increasing with union
$g(H)$; but $g(H)$ by the assumption on $H$ is
$\aleph_0$-indecomposable, hence $g(H) = g(H) \cap G^*_n 
\subseteq G^*_n$ for some $n$, so $g \rest H \in \cG$ 
and the rest should be clear.
\medskip

\noindent
\underline{$(B) \Rightarrow (B)^+$}:

What about $(B)^+(c)$? while $\bar G$ does not necessarily satisfy it,
we can ``correct it", e.g. we choose $u_\alpha,v_\alpha$ and if
$\alpha \notin \cup\{v_\beta:\beta < \alpha\}$ also $G'_\alpha$ by
induction on $\alpha < \alpha_*$ such that (the idea is that if $\beta
\in v_\alpha,G_\beta$ is discarded being embeddable into some $G'$ and
$G'_\alpha$ is the ``corrected" member):
\mn
\begin{enumerate}
\item[$(*)^2_\alpha$]
\begin{enumerate}
\item[(a)]  $G_\alpha \subseteq G'_\alpha \in \bold K^{\lf}_\mu$ if
  $\alpha \notin \cup\{v_\beta:\beta < \alpha\}$;
\sn
\item[(b)]  $u_\alpha \subseteq \alpha$ and $v_\alpha
  \subseteq \alpha_* \backslash \alpha$;
\sn
\item[(c)]  if $\beta < \alpha$ then $u_\beta =
  u_\alpha \cap \beta$ and $u_\alpha \cap v_\beta =\emptyset$;
\sn
\item[(d)]   if $\alpha = \beta +1$ then $\beta \in
  u_\alpha$ \Iff \, $\beta \notin \cup\{v_\gamma:\gamma < \beta\}$;
\sn
\item[(e)]   if $\alpha \notin \cup\{v_\gamma:\gamma < \alpha\}$, \then \,:
\begin{enumerate}
\item[$\bullet_1$]  $\gamma \in v_\alpha$ \Iff \, 
$G_\gamma$ is embeddable into $G'_\alpha$ over $H$;
\sn
\item[$\bullet_2$]   if $\gamma \in \alpha_* \backslash
 (\alpha +1) \backslash \cup\{v_\beta:\beta \le \alpha\}$ 
then $G_\gamma$ is not embeddable over
  $H$ into any $G'$ satisfying $G'_\alpha \subseteq G' \in \bold K_{\lf}$;
\end{enumerate}
\sn
\item[(f)]   if $\alpha = \beta +1$ and $\beta \notin
  u_\alpha$ then $v_\beta = \emptyset$.
\end{enumerate}
\end{enumerate}
\mn
This suffices because if we let $u_{\alpha_*} = \alpha_* \backslash
\cup\{v_\gamma:\gamma < \alpha_*\}$, 
then $\langle G'_\alpha:\alpha \in u_{\alpha_*}\rangle$ is
as required.  Why can we carry the induction?

For $\alpha=0,\alpha$ limit we have nothing to do because $u_\alpha$
is determined by $(*)^2_\alpha(c)$ and $(*)^2_\alpha(d)$.   For $\alpha =
\beta +1$, if $\beta \in \bigcup\limits_{\gamma < \beta} v_\gamma$ we
have nothing to do, in the remaining case we choose $G'_{\beta,i}
\in \bold K^{\lf}_\mu$ by induction on 
$i \in [\alpha,\alpha_*]$, increasing continuous
with $i,G'_{\beta,\alpha} = G_\alpha,G'_{\beta,i+1}$ make clause (e)
true, i.e. if $G'_{\beta,i}$ has an
extension into which $G_i$ is embeddable over $H$, then there is such
an extension of cardinality $\mu$ and choose $G'_{\beta,i+1}$ as such
an extension.

Lastly, let $G'_\alpha = G'_{\alpha,\alpha_*}$ and $u_\alpha = u_\beta
\cup \{\alpha\}$ and $v_\alpha = \{i:i \in \alpha_*,i \notin \cup
\{v_\gamma:\gamma < \beta\}$ and $G_i$ is embeddable into $G_\beta$
over $H\}$.

So clause $(B)^+(c)$ holds; and clause $(B)^+(d)$ follows from
$(B)(b)+(B)^+(c)$, so we are done.
\medskip

\noindent
\underline{$(B)^+ \Rightarrow (A)$}:

We prove below more: there is something like ``special model", i.e. part
(2), now $(C) \Rightarrow (A)$ is trivial so we are left with the following.
\medskip

\noindent
\underline{$(B)^+ \Rightarrow (C)$}:

Choose
$\bar\lambda = \langle \lambda_n:n < \omega\rangle,\lambda_n \in
\{\aleph_{\alpha +2}:\alpha \in \Ord\}$ such that $2^{\lambda_n} <
\lambda_{n+1}$ and $\mu = \sum\limits_{n} \lambda_n$.

Let $\bold K^{\slf}_{\bar\lambda}$ be the class of $\bar G$ such that:
\mn
\begin{enumerate}
\item[$(*)^3_{\bar G}$]  
\begin{enumerate}
\item[(a)]  $\bar G = \langle G_n:n < \omega\rangle$ is 
increasing (hence if $H \subseteq
  \bigcup\limits_{n} G_n$ is $\aleph_0$-indecomposable \then \, $(\exists n)(H
\subseteq G_n))$
\sn
\item[(b)]  $G_n \in \bold K_{\lf}$ has cardinality $2^{\lambda_n}$
\sn
\item[(c)]  if $H \subseteq G_n,|H| = \lambda_n$, \then \,
there is $\aleph_0$-indecomposable $H' \in \bold K^{\lf}_{\lambda_n}$, 
 e.g. as in \ref{b18} such that $H \subseteq H' \subseteq G_n$
\sn
\item[(d)] if $H \subseteq G_n$ is
$\aleph_0$-indecomposable and $H \in \bold K^{\lf}_{\le \lambda_n}$
\then \, the pair $(H,\bigcup\limits_{m} G_m)$ is an amalgamation base
(see Definition \ref{a31}(2)); moreover, it is as in $(B)^+(d)$;
\sn
\item[(e)]   if $H \subseteq G_n$ is
  $\aleph_0$-indecomposable of cardinality $\le \lambda_n,H \subseteq
  H' \subseteq G' \in \bold K^{\lf}_\mu,H'$ is
  $\aleph_0$-indecomposable\footnote{The $\aleph_0$-indecomposability
    is not always necessary, but we need it sometimes.}, 
$|H'| \le \lambda_n$ and
  $\bigcup\limits_{n} G_n,G'$ are compatible over $H$ (in $\bold
  K^{\lf}_\mu$), \then \, $G'$
  is embeddable into $\bigcup\limits_{m} G_m$ over $H$, noting the
  necessarily $H'$ mapped into some $G_m$.
\end{enumerate}
\end{enumerate}
\mn
Now we can finish because clearly:
\mn
\begin{enumerate}
\item[$(*)_4$]  if $G \in \bold K^{\lf}_{\le \mu}$ \then \, for some
 $\bar G \in \bold K^{\slf}_{\bar\lambda},G$ is embeddable into
  $\bigcup\limits_{n} G_n$;
\sn
\item[$(*)_5$]
\begin{enumerate}
\item[(a)]   if $\bar G^1,\bar G^2 \in \bold
  K^{\slf}_{\bar\lambda}$ then
$\bigcup\limits_{n} G^1_n,\bigcup\limits_{n} G^2_n$ are isomorphic;
\sn
\item[(b)]   if $\bar G^1,\bar G^2 \in \bold K^{\slf}_{\bar\lambda},
H \in \bold K^{\lf}$ is $\aleph_0$-indecomposable and $f_\ell$ embeds
$H$ into $G^\ell_n$, for $\ell=1,2$, and this diagram 
can be completed, (i.e. there are $G \in K^{\lf}_\mu$ and embedding
$g_\ell:\bigcup\limits_{k} G^\ell_k \rightarrow G_*$ such
that $g_1 \circ f_1 = g_2 \circ f_2$) \then \, there is $h$ such that:
\sn
\begin{enumerate}
\item[$(\alpha)$]  $h$ is an isomorphism from $\bigcup\limits_{k}
  G^1_k$ onto $\bigcup\limits_{k} G^2_k$;
\sn
\item[$(\beta)$]  $h \circ f_1 = f_2$;
\sn
\item[$(\gamma)$]  $h$ maps $G^1_{n+k}$ into $G^2_{n+\ell}$ for some
  $\ell$;
\sn
\item[$(\delta)$]  $h^{-1}$ maps $G^2_{n+k}$ into $G^1_{n+\ell}$ for
  some $\ell$.
\end{enumerate}
\end{enumerate}
\end{enumerate}
\end{PROOF}

\begin{discussion}
\label{b28}
We can sort out which claims of this section holds for any suitable
a.e.c.  E.g. $\bold K$ is a universal class with $|\tau_{\bold K}| <
\mu$ and the parallel of \ref{b18} holds.
\end{discussion}

\begin{claim}
\label{b47}
1) We can generalize \ref{b25}(1)(2) to any $\gk$ such that:
\mn
\begin{enumerate}
\item[(a)]  $\gk$ is an a.e.c. with $\LST_{\gk} < \mu$
\sn
\item[(b)]   every $M \in K^{\gk}_{<\mu}$ has an
$\aleph_0-\gk$-indecomposable $\le_{\gk}$-extension from $K^{\gk}_{<
  \mu}$
\sn
\item[(c)]   $\gk_{\le \mu}$ has the $\JEP$.
\end{enumerate}
\mn
2) We can in \ref{b25}(1) $(A) \Rightarrow (B) \Rightarrow (B)^+$
replace $\aleph_0$ by $\partial = \cf(\partial) = \cf(\mu)$ so $\mu$
strong limit.
\end{claim}

\begin{PROOF}{\ref{b47}}
As above.
\end{PROOF}

\begin{claim}
\label{b29}
1) For every $n \ge 2$ ($n=2$ is eough).  There is $\gs \in \gS$ such
that: if $G \subseteq H \in \bold K_{\lf}$ and $a \in G$ has order $n$
and $q_{\gs}(a,G) = \tp(\gamma \langle b_0,b_1\rangle,G,H)$, \then \,
$b_\ell$ commutes with $G$ and $[b_0,b_1] = a$.

\noindent
2) Assume $G \in \bold K_{\lf},a_\alpha \in G$ for $\alpha < \lambda$
has order 2, $a_\gamma = e_G$, \then \, there are $H,\bar b_0,\bar b_1$
such that:
\mn
\begin{enumerate}
\item[(a)]  $G \subseteq H \in \bold K_{\lf}$
\sn
\item[(b)]  $\bar b_0,\bar b_2 \in {}^\lambda H$
\sn
\item[(c)]  $H$ is generated by $\bar b_0 \cup \bar b_1 \cup G$
\sn
\item[(d)]  $[b_{0,s},b_{1,s}) = a_s$
\sn
\item[(e)]  $b_{\ell,\alpha}$ commute with
  $\cup\{b_{\iota,\alpha,\beta}:
\iota < 2,\beta \in \lambda \backslash \{\alpha\}\} \cup G$
\sn
\item[(f)]  $\langle b_{0,\alpha},b_{1,\alpha}\rangle \in {}^2 H$
  realizes $q_{\gs}(a_\alpha,G)$ and even $g_{\gs}(a_\alpha,G_\alpha)$
  where $G_\alpha = \seb_H(\{a_{\iota,\beta}:\iota < 2,\beta <
  \lambda,\beta \ne \alpha\} \cup G)$.
\end{enumerate}
\end{claim}
\bigskip

\subsection {Universal in $\beth_\omega$} \label{3C} \
\bigskip

In \S(3B) we have special models in $\bold K$ of cardinality,
e.g. $\beth_\omega$.  Our intention is to first investigate
$\gk_{\fnq}$ (the class structures consisting of a set and 
a directed family of equivalence
relations on it, each with a finite bound on the size of equivalence
classes).  So $\gk_{\fnq}$ is similar to $\bold K_{\lf}$ but seems
easier to analyze.  We consider some partial orders on $\gk = \gk_{\fnq}$.

First, under the substructure order, $\le_1 = \subseteq$, this class
fails amalgamation.  Second, another order, $\le_2$ demanding $\TV$
for countably many points, finitely many equivalence relations, we
have amalgamation.  
Third, we add: if $M \le_3 N$ then $M \le_1 N$ and the union of
$(P^n,E_d)_{d \in Q(M)}$ is the disjoint union of models isomorphic to
$(P^M,E_d)_{d \in Q(M)}$, the equivalence relation is $E_{M,N}$.  This
is intended to connect to locally finite groups.  So we may instead
look at $\{f \in \Sym(N)$: if $a \in N \backslash M$ and $a/E_{M,N} \ncong
M$ then $f \rest (a/E_{M,N}) = \id(a/E_{M,N})$; no need of representations.

The model in $\beth_\omega$ will be $\bigcup\limits_{n} M_n,\|M_n\| =
\beth_{n+1}$, gotten by smooth directed unions of members of
cardinality $\beth_\omega$ by $\bold I_n \subseteq P^{M_{n+1}}$ is a
set of representatives for $E_{M_n,M_{n+1}}$.
\bigskip

\begin{definition}
\label{b55}
Let $\bold K = \bold K_{\fnq}$ be the class of structures $M$ such
that (the vocabulary is defined implicitly and is $\tau_{\bold K}$,
i.e. depends just on $\bold K$):
\mn
\begin{enumerate}
\item[(a)]  $P^M,Q^M$ is a partition of $M,P^M$ non-empty;
\sn
\item[(b)]  $E^M \subseteq P^M \times P^M \times Q^M$ (is a
  three-place relation) and we write $a E^M_c b$ for $(a,b,c) \in E^M$;
\sn
\item[(c)]  for $c \in Q^M,E^M_c$ is an equivalence relation on $P^M$
  with $\sup\{|a/E^M_c|:a \in P^M\}$ finite (see more later);
\sn
\item[(d)]  $Q^M_{n,k} \subseteq (Q^M)^n$ for $n,k \ge 1$
\sn
\item[(e)]  if $\bar c = \langle c_\ell:\ell < n\rangle \in {}^n(Q^M)$ we
  let $E^M_{\bar c}$ be the closure of $\bigcup\limits_{\ell} E_\ell$
  to an equivalence relation;
\sn
\item[(f)]  ${}^n(Q^M) = \bigcup\limits_{k \ge 1} Q^M_{n,k}$;
\sn
\item[(g)]  if $\bar c \in Q^M_{n,k}$ then $k \ge |a/E^M_{\bar c}|$
  for every $a \in P^M$.
\end{enumerate}
\end{definition}

\begin{definition}
\label{b58}
We define some partial order on $\bold K$.

\noindent
1) $\le_1 = \le^1_{\bold K} = \le^1_{\fnq}$ is being a sub-model.

\noindent
2) $\le_3 = \le^3_{\bold K} = \le^3_{\fnq}$ is the following: $M \le_3
N$ iff:
\mn
\begin{enumerate}
\item[(a)]  $M,N \in \bold K$
\sn
\item[(b)]  $M \subseteq N$
\sn
\item[(c)]  if $A \subseteq N$ is countable and $A \cap Q^N$ is
  finite, \then \, there is an embedding of $N \rest A$ into $M$ over
  $A \cap M$ or just a one-to-one homomorphism.
\end{enumerate}
\mn
3) $\le_2 = \le^2_{\bold K} = \le^2_{\fnq}$ is defined like $\le_3$
but in clause (c), $A$ is finite.
\end{definition}

\begin{claim}
\label{b61}
1) $\bold K$ is a universal class, so $(\bold K,\subseteq)$ is an
a.e.c.

\noindent
2) $\le^3_{\bold K},\le^2_{\bold K},\le^1_{\bold K}$ are partial oders
  on $\bold K$.

\noindent
3) $(\bold K,\le^2_{\bold K})$ is an a.e.c.

\noindent
4) $(\bold K,\le^2_{\bold K})$ has disjoint amalgamation.
\end{claim}

\begin{PROOF}{\ref{b61}}
1),2),3) Easy.

\noindent
4) By \ref{b64} below.
\end{PROOF}

\begin{claim}
\label{b64}
If $M_0 \le^1_{\bold K} M_1,M_0 \le^3_{\bold K} M_2$ and $M_1 \cap M_2
= M_0$, \then \, $M = M_1 + M_2$, the disjoint sum of $M_1,M_2$
belongs to $\bold K$ and extends $M_\ell$ for $\ell=0,1,2$ and even
$M_1 \le^3_{\fnq} M$ and $M_0 \le^2_{\bold K} M_1 \Rightarrow M_2
\le^2_{\bold K} M$ \when \,:
\mn
\begin{enumerate}
\item[$(*)$]  $M = M_1 + M_2$ which means $M$ is defined by:
\sn
\begin{enumerate}
\item[(a)]  $|M| = |M_1| \cup |M_2|$;
\sn
\item[(b)]  $P^M = P^{M_1} \cup P^{M_2}$;
\sn
\item[(c)]  $Q = Q^{M_1} \cup Q^{M_2}$;
\sn
\item[(d)]  we define $E^M$ by defining $E^M_c$ for $c \in Q^M$ by
  cases:
\begin{enumerate}
\item[$(\alpha)$]  if $c \in Q^{M_0}$ then $E^M_c$ is the closure of
  $E^{M_1}_\ell \cup E^{M_2}_\ell$ to an equivalence relation;
\sn
\item[$(\beta)$]  if $c \in Q^{M_\ell} \backslash Q^{M_0}$ and $\ell
  \in \{1,2\}$ then $E^M_c$ is defined by
\begin{itemize}
\item  $a E^M_c b$ \Iff \, $a = b \in P^{M_{3-\ell}} \backslash M_0$ or
$a E^{M_\ell}_c b$ so $a,b \in P^{M_\ell}$;
\end{itemize}
\end{enumerate}
\sn
\item[(e)]  $Q^M_{n,k} = Q^{M_1}_{n,k} \cup Q^{M_2}_{n,k} \cup \{\bar
  c:\bar c \in {}^n(Q^M) \backslash ({}^n(Q^{M_1})) \cup
  {}^n(Q^{M_2})\}$.
\end{enumerate}
\end{enumerate}
\end{claim}

\begin{PROOF}{\ref{b64}}
Clearly $M$ is a well defined structure, extends $M_0,M_1,M_2$ and
satisfies clauses (a),(b),(c) of Definition \ref{b55}.  There are two
points to be checked: $a \in P^M,\bar c \in Q^M_{n,k} \Rightarrow
|a/E^M_{\bar c}| \le k$ and ${}^n(Q^M) = \bigcup\limits_{k \ge 1}
Q^M_{n,k}$
\mn
\begin{enumerate}
\item[$(*)_1$]  if $a \in P^M$ and $\bar c \in Q^M_{n,k}$ \then \,
  $|a/E^M_{\bar c}| \le k$.
\end{enumerate}
\mn
Why?  If $\bar c \in Q^M_{n,k} \backslash (Q^{M_1}_{n,k} \cup
Q^{M_2}_{n,k})$ this holds by the definition, so assume $\bar c \in
Q^{M_\iota}_{n,k}$.  If this fails, then there is a finite set $A
\subseteq M$ such that $\bar c \subseteq A,a \in A$ and letting $N = M
\rest A$ we have $|a/E^N_{\bar c}| > k$.  By $M_0 \le^1_{\bold K}
M_1,M_0 \le^3_{\bold K} M_2$ (really $M_0 \le^2_{\bold K} M_2$
suffice) there is a one-to-one homomorphism $f$ from $A \cap M_2$ into
$M_0$.  Let $B' = (A \cup M_1) \cup f(A \cap M_2)$ and $N' = M \rest B$
and let $g = f \cup \id_{A \cap M_1}$.  So $g$ is a homomorphism from
$N$ onto $N'$ and $g(a)/E^{N'}_{g(\bar c)}$ has $>k$ members, which
implies $g'(a)/E^{M_1}_{g'(\bar c)}$ has $>k$ members.  Also $g(\bar
c) \in Q^{M_1}_{n,k}$.  (Why?  If $\iota=1$ trivially, if $\iota=2$ by
the choice of $f$, contradiction to $M \in \bold K$.)]
\mn
\begin{enumerate}
\item[$(*)_2$]  if $\bar c \in {}^n(Q^M)$ then $\bar c \in
  \bigcup\limits_{k} Q^M_{n,k}$.
\end{enumerate}
\mn
Why?  If $\bar c \in M_1$ or $\bar c \subseteq M_2$, this is obvious
by the definition of $M$, so assume that they fail.  By the definition
of the $Q^M_{n,k}$'s we have to prove that $\sup\{a/E^M_{\bar c}:a \in
P^M\}$ is infinite.  Toward contradiction assume this fails for each
$k \ge 1$ there is $a_k \in P^M$ such that $a_k/E^M_{\bar c}$ has $\ge
k$ elements hence there is a finite $A_k \subseteq M$ such that 
$a_k/E^{M \rest A_k}_{\bar c}$ has $\ge k$ elements.  Let $A =
\bigcup\limits_{k \ge 1} A_k$, so $A_k$ is a countable subset of $M$
and we continue as in the proof of $(*)_1$.

Additional points (not really used) are proved like $(*)_2$:
\mn
\begin{enumerate}
\item[$(*)_3$]  $M_1 \le^3_{\bold K} M$;
\sn
\item[$(*)_4$]  $M_0 \le^2_{\bold K} M_1 \Rightarrow M_2 \le^2_{\bold K} M$;
\sn
\item[$(*)_5$]  $M_1 +_{M_0} M_2$ is equal to $M_2 +_{M_0} M_1$.
\end{enumerate}
\end{PROOF}

\begin{claim}
\label{b67}
If $\lambda = \lambda^{< \mu}$ and $M \in \bold K$ has cardinality
$\le \lambda$ \then \, there is $N$ such that:
\mn
\begin{enumerate}
\item[(a)]  $N \in \bold K_\lambda$ extend $M$;
\sn
\item[(b)]  if $N_0 \le^3_{\bold K} N_1$ and $N_0$ has cardinality $<
  \mu$ and $f_0$ embeds $N_0$ into $N$, \then \, there is an embedding
  $f_1$ of $N_1$ into $N$ extending $f_0$.
\end{enumerate}
\end{claim}
\newpage

\section {Density of being Complete in $\bold K^{\lf}_\lambda$} \label{4}

We prove here that for quite many cardinals $\lambda$, the complete $G
\in \bold K^{\exlf}_\lambda$ are dense in $(\bold
K^{\exlf}_\lambda,\bold c)$; e.g. every $\lambda$ successor of
regular satisying $\lambda = \lambda^{\aleph_0} \vee \lambda \ge
\beth_\omega$. 

\begin{discussion}
\label{c20}
We would like to prove for as many cardinals $\mu = \lambda$ or at
least pairs $\mu \le \lambda$ of cardinals we can prove 
that $(\forall G \in \bold K^{\lf}_\mu)(\exists H \in 
\bold K^{\exlf}_\lambda)(G \subseteq H \wedge H$ complete).  We
necessarily have to assume $\lambda \ge \mu + \aleph_1$.  So far we
have known it only for $\lambda = \mu^+,\mu = \mu^{\aleph_0}$, (and $\lambda =
\aleph_1,\mu = \aleph_0$, see the introduction of \cite{Sh:312}).  
We would like to prove it also for as many pairs of 
cardinals as we can and even for $\lambda = \mu$.
In this section, we shall consider using 
$\Pr_0(\lambda,\lambda,\lambda,\aleph_0)$.
\end{discussion}

\noindent
In this section, in particular in \ref{c23}(2) we rely on \cite{Sh:312}.
\begin{hypothesis}
\label{c23}
1) $\lambda > \theta = \cf(\theta) > \aleph_0$.

\noindent
2) $\bold K = \bold K_{\lf}$.

\noindent
3) $\gS$ a set of schemes consisting of all of them or is just of
cardinality $\le \lambda$, is dense and containing enough of
those mentioned in \cite[\S2]{Sh:312}.  Also $c \ell(\gS) = \gS$, i.e.
$\gS$ is closed, see \cite[1.6=La21,1.8=La22]{Sh:312} hence by
\cite{Sh:312} there is such countable $\gS$.  
Recall $G \le_{\gS} H$ means that $G \subseteq H$ and for 
every $\bar b \in {}^{\omega >}H$ for some
$\bar a \in {}^{\omega >} G$ and $\gs \in c \ell(\gS)$ we have
$\tp_{\bs}(\bar b,G,H) = q_{\gs}(\bar a,G)$.

\noindent
4) $S = \langle S_1,S_2,S_3\rangle$ is a partition of $\theta
\backslash \{0\}$ to stationary subsets, 
$S_3 \subseteq S^\theta_{\aleph_0} := \{\delta <
\theta:\cf(\delta) = \aleph_0\}$ and $\{\varp +1:\varp \in (0,\theta)\}
\subseteq S_1$ and $\zeta \in S_2 \Rightarrow \omega^2|\zeta$; 
we may let $S_0 = \{0\}$.
\end{hypothesis}

\begin{definition}
\label{c25}
Let $\bold M_1 = \bold M^1_{\lambda,\theta,\bar S}$ be the class of objects
$\bold m$ which consists of:
\mn
\begin{enumerate}
\item[$(a)$]  $G_i = G_{\bold m,i}$ for $i \le \theta$ is increasing
  continuous, $G_0$ is the trivial group with universe $\{0\},G_1 \in \bold
  K_{\le \lambda}$ and for $i \in (\theta +1) \backslash \{0,1\}$ the group
$G_i \in \bold K_\lambda$ has universe $\{\theta \alpha
  +j: \alpha < \lambda$ and $j<1+i\}$ and so $e_{G_i} =0$;
\sn
\item[$(b)$]  if $i < \theta$, \then \, we have:
\sn
\begin{enumerate}
\item[$(\alpha)$]
\begin{itemize}
\item  sequences $\bold b_i = \langle \bar b_{i,s}:s \in
  I_i\rangle,\bold a_i = \langle \bar a_{i,s}:s \in J_i\rangle$;
\sn
\item  each $\bar a_{i,s}$ is a finite sequence from $G_i$;
\sn
\item  each $\bar b_{i,s}$ is a finite sequence from $G_{i+1}$;
\sn
\item $I_i$ is a linear order of cardinality
  $\lambda$ with a first element;
\sn
\item  $J_i$ is a set or linear order of cardinality $\le \lambda$;
\sn
\item  if $i=0$ then $J_i = \lambda = I_i \subseteq \lambda$ 
and $\langle \bar b_{i,s}
  = \langle b_{i,s}\rangle:s \in I_0\rangle$ lists the members of $G_1$
  and $\bar a_{i,s} = \langle \rangle$;
\end{itemize}
\sn
\item[$(\beta)$]  $G_{i+1}$ is generated by $\cup\{\bar b_{i,s}:s \in
  I_i\} \cup G_i$;
\sn
\item[$(\gamma)$]  $\bar a_{i,s} \in {}^{\omega >}(G_i)$ and $\bar
  a_{i,\min(I_i)} = e_{G_i}$; 
\sn
\item[$(\delta)$]  $\bold c_i:[I_i]^2 \rightarrow \lambda$;
\end{enumerate}
\sn
\item[$(c)$]  [toward being in $\bold K_{\exlf}$] 
if $i \in S_1$, \then \, $J_i = I_i$ and we also have
  $\langle \gs_{i,s}:s \in I_i\rangle$ such that:
\sn
\begin{enumerate}
\item[$(\alpha)$]  $\gs_{i,s} \in \gS$;
\sn
\item[$(\beta)$]  $\tp_{\bs}(\bar b_{i,s},G_i,G_{i+1}) = 
q_{\gs_{i,s}}(\bar a_{i,s},G_i)$ so $\ell g(\bar b_{i,s}) =
n(\gs_{i,s})$ and $\ell g(\bar a_{i,s}) = k(\gs_{i,s})$;
\sn
\item[$(\gamma)$]  if $s_0 <_{I_i} \ldots <_{I_i} s_{n-1}$ then
$\tp_{\bs}(\bar b_{i,s_0} \char 94 \ldots \char 94 
\bar b_{i,s_{n-1}},G_i,G_{i+1})$ is gotten by $(\gs_{i,s_0},\bar
a_{i,s_0}),\ldots,(\gs_{i,s_{n-1}},\bar a_{i,s_{n-1}})$ by one of
the following two ways: 
\newline
\underline{Option 1}: we use the linear order 
$I_i$ on $\lambda$ so $\tp_{\qf}(\bar
  b_{i,s},G_{i,s},G_{i,t})$ is equal to $q_{\gs_{i,s}}(\bar
  a_{i,s},G_{i,s})$ where $G_{i,s}$ is the subgroup of $G_{i+1}$
  generated by $G_i \cup \{\bar b_{i,t}:t <_{I_i} s\}$, see
  \cite[\S(1C),1.28=La58]{Sh:312}; 
\newline
\underline{but}\footnote{Option 1
is useful in some generalizations to $K_{\gk}$ not closed under
products.} we choose:
\newline
\underline{Option 2}: intersect the atomic types over all
  orders on $\{\alpha_0,\dotsc,\alpha_{n-1}\}$ each gotten as in
  Option 1, so $I_i$ can be a set of cardinality $\lambda$, see
  \cite[\S3]{Sh:312}; 
\sn
\item[$(\delta)$]  $\bold c_i$ is constantly zero;
\end{enumerate}
\sn
\item[$(d)$]  [toward indecomposability] if $i \in S_2$ \then \,: 
\sn
\begin{enumerate}
\item[$(\alpha)$]  $J_i \subseteq \lambda$ and $J_i = \bigcup\limits_{\alpha <
    \lambda} J_{i,\alpha}$, disjoint union
\sn
\item[$(\beta)$]  $\langle I_{i,\alpha}:\alpha < \lambda\rangle$ is a
  partition of $I_i \subseteq \lambda$;
\sn
\item[$(\gamma)$]  $\bar a_{i,\alpha} = 
\langle a_{i,\alpha}\rangle,\bar b_{i,s} = \langle b_{i,s}\rangle$ and
$a_{i,0} = e_{G_i}$;
\sn
\item[$(\delta)$]  $G_i$ is generated by $\{a_{i,\alpha}:\alpha < \lambda\}$;
\sn
\item[$(\varp)$]  if $s < t$ and $\bold c_i\{s,t\} = \alpha <
  \lambda$ then in $G_{i+1}$ we have $[b_{i,s},b_{i,t}] =
a_{i,\alpha}$;
\sn
\item[$(\zeta)$]
\begin{enumerate}
\item[$\bullet_1$]  if $\beta_1,\beta_2 \in J_{i,\alpha}$ then
  $b_{i,\beta_1},b_{i,\beta_2}$ are commuting each of order 2;
\sn
\item[$\bullet_2$]   if $s \in I_i$ then $b_{i,s}$ commutes with $G_i$
  (in $G_{i+1}$);
\end{enumerate}
\sn
\item[$(\eta)$] 
\begin{enumerate}
\item[$\bullet_1$]   if $s \ne t \in I_{i,\alpha}$ then $\bold
  c_i\{s,t\} \in J_{i,\alpha}$ and
\sn
\item[$\bullet_2$]   if $s \in I_{i,\alpha},t \in
  I_{i,\beta}$ and $\alpha \ne \beta$ then $\bold c_i\{s,t\} = 0$;
\end{enumerate}
\sn
\item[$(\theta)$]  $\bold c_{i,\alpha} = \bold c_i \rest [I_{i,\alpha}]^2$;
\end{enumerate}
\sn
\item[$(e)$]  [against external automorphism] 
if $i \in S_3$ then $\cf(i) = \aleph_0$ and 
$\bar j_i,\langle I_{i,\alpha}:\alpha <\lambda\rangle,\langle
  a_{i,\alpha}:\alpha < \lambda\rangle$ satisfies:
\sn
\begin{enumerate}
\item[$(\alpha)$]  $\bar j_i = \langle j_{i,n}:n < \omega\rangle$ is
  increasing with limit $i$;
\sn
\item[$(\beta)$]  $a_{i,\omega \alpha +\ell} \in G_{j_{i,\ell +1}}$
  commutes with $G_{j_{i,\ell}}$ and if $\ell \ne 0$ then it has order 2, and
  $\notin G_{j_{i,\ell}}$ and $a_{i,\omega \alpha} \equiv e_{G_i}$;
  moreover:
\begin{itemize}
\item  for some infinite $v \subseteq \omega \backslash\{0\}$
  we\footnote{An alternative is $v = \omega \backslash
    \{0\},a_{i,\omega \alpha + \ell} \in \bold
    C(G_{j_{i,\ell}},G_{j_{i,\ell +1}})$.  In this case in
    \ref{c26}$(e)(\varp)$  we naturally have $c_\varp \in \bold
    C(G_{i_\varp},G_{i_{\varp +1}})$ and $\ell_0 = 1,\ell_1
    =2,\ldots$.  But then we have to be more careful in \ref{c29},
    e.g. in \ref{c29}(1) if we assume, e.g. $\lambda =
    \lambda^{\langle \theta;\theta\rangle}$ and $\theta > \aleph_1$
    all is O.K. (recalling we have guessing clubs on
    $S^\theta_{\aleph_0}$).  However, using $\gs_{\cg}$, see
    (\cite[2.17=Lc50]{Sh:312}), the present is enough here.} have $\ell
  \in \omega \backslash v \Rightarrow a_{i,\omega \alpha + \ell} =
  e_{G_i},\ell \in v \Rightarrow a_{i,\omega \alpha + \ell} \in 
\bold C(G_{j[i,\omega \alpha + \ell]},G_{j[i,\omega \alpha + \ell]+1})$,
  where:
\sn
\item $j[i,\omega \alpha + \ell) \in [j_{i,\ell},j_{i,\ell +1})$;
\end{itemize}
\sn
\item[$(\gamma)$]  $\bar I_i = \langle I_{i,\alpha}:\alpha <
  \lambda\rangle$ is a partition of $I_i$; for $s \in I_i$ let
$\alpha_i(s)$ be the $\alpha < \lambda$ such that $s \in I_{i,\alpha}$
and let $\bold c_{i,\alpha} = \bold c_i \rest [I_{i,\alpha}]^2$;
\sn
\item[$(\delta)$]  if $s,t \in I_{i,\alpha}$ then $[b_{i,s},b_{i,t}]=
  a_{i,\bold c_i\{s,t\}}$ and $\bold c_i\{s,t\} \in 
\{\omega \alpha + \ell:\ell < \omega\}$;
\sn
\item[$(\varp)$]  if $s,t \in I_i$ and $\alpha_i(s) \ne
  \alpha_i(t)$ then $[b_{i,s},b_{i,t}] = e_{G_i}$.
\end{enumerate}
\end{enumerate}
\end{definition}

\begin{convention}
\label{c25a}
If the identity of $\bold m$ is not clear, we may write $G_{\bold
  m,i}$, etc., but if clear from the context we may not add it.
\end{convention}

\begin{definition}
\label{c26}
1) Let $\bold M_2 = \bold M^2_{\lambda,\theta,\bar S}$ be the set of $\bold m
\in \bold M_1$ when we add in Definition \ref{c25}:
\mn
\begin{enumerate}
\item[(c)]
\begin{enumerate}
\item[$(\varp)$]   if $\gs \in \gS,i \in S_1,\bar a \in 
{}^{n(\gs)}(G_i)$ and $k=k(\gs)$, \then \, for $\lambda$
 elements $s \in I_i$ we have $(\gs_{i,s},\bar a_{i,s}) = 
(\gs,\bar a)$; 
\end{enumerate}
\sn
\item[(d)]
\begin{enumerate}
\item[$(\iota)$]  if $i \in S_2$ and $\langle I'_j:j < \theta\rangle$ is
  increasing with union $I_i$ \then \, for some 
$j < \theta,G_i \subseteq \seb(\{b_{i,s}:s \in I'_j\},G_{i+1})$; this
follows from $(\kappa) + (\mu)$ below; 
\sn
\item[$(\kappa)$]  if $i \in S_2$ and $\alpha < \lambda$ \then \, $\bold
  c_{i,\alpha}$, a function from $[I_{i,\alpha}]^2$ onto $J_{i,\alpha}$,
 is $\theta$-indecomposable; 
\sn
\item[$(\lambda)$]  if $i \in S_2$ and $\alpha \in \lambda$ \then \,
  $|J_{i,\alpha}| = \theta$ and $\langle b_{i,s}:s \in
  J_{i,\alpha}\rangle$ pairwise commute, each has order 2 and are
  independent;
\sn
\item[$(\mu)$]   if $i \in S_2$ and $\langle a_\varp:\varp <
  \theta\rangle \in {}^\theta(G_i)$ is a sequence of pairwise
  commuting independent elements of order 2 (as in \ref{b18}) \then \,
  for some $\alpha < \lambda,\{b_{i,\gamma}:\gamma \in J_{i,\alpha}\}
  \cap \{a_\varp:\varp < \theta\}$ has cardinality $\theta$;
\end{enumerate}
\sn
\item[$(e)$]
\begin{enumerate}
\item[$(\zeta)$]   if $\langle i_\varp:\varp <
  \theta\rangle$ is increasing continuous and $i_\varp < \theta$ and
  $c_\varp \in \bold C'(G_{i_\varp},G_{i_\varp +1})$
 has order $2$ and for transparency $c_\varp \notin
  G_{i_\varp}$ then for some 
\newline
$(i,\alpha,v,\ell_0,\ell_1,
\ldots,\varp_0,\varp_1,\ldots)$ we have:
\begin{enumerate}
\item[$\bullet_1$]  $i < \theta,\alpha < \lambda$ and $v \subseteq w
  \backslash \{0\}$ is infinite;
\sn
\item[$\bullet_2$]  $\varp_0 < \varp_1 < \ldots < \theta$
  and $1 \le \ell_0 < \ell_1 < \ldots$;
\sn
\item[$\bullet_3$]  $i = \cup\{\varp_n:n < \omega\}$;
\sn
\item[$\bullet_4$]  $j_{i,\omega \alpha + \ell_n} \le i_{\varp_n}
  < j_{i,\theta,\alpha + \ell_{n+1}}$ and $a_{\bold m,i,\omega \alpha
    + \ell_n} = c_{\varp_n}$;
\end{enumerate}
\end{enumerate}
\sn
\item[(f)]  $\langle I_i:i < \theta\rangle$ are pairwise disjoint
with union $I$.
\end{enumerate}
\mn
1A) Let $M_{1.5} = \bold M^{1.5}_{\lambda,\theta,\bar S}$ be the set of
$\bold m \in \bold M_1$ such that in part (1) weakening some clauses:
\mn
\begin{enumerate}
\item[(M)$'$]  $(\alpha) \Rightarrow (\beta)$ where
\sn
\begin{enumerate}
\item[$(\alpha)$]  if $i \in S_2,j < i,j \in S_1,I^{\cg}_j = \{s
  \in I_i: \gs_{j,s} = \gs_{\cg}\},H_{j,s} = b^{-1}_{j,s} G_j b_{j,s}$
  hence (see below) $\langle H_{j,s}:s \in I^{\cg}_j\rangle$ are
  pairwise commuting subgroups of $G_i$ and $a_s \in \CH_{j,s}
  \backslash \{e\}$ for $s \in I^{\cg}_j$
\sn
\item[$(\beta)$]  for some $\alpha < \lambda$ and sequence $\langle
  s(\varp):\varp < \theta\rangle$ of pairwise distinct members of
  $I^{\cg}_j$ we have $\{b_{i,t}:t \in J_{i,\alpha}\} \cap
  \{a_{s(\varp)}:\varp < \theta\}$ has cardinality $\theta$.
\end{enumerate}
\end{enumerate}
\mn
[Pre 17.05.30 version: Retain? (2017.10.30)

\noindent
2) Let $\bold M_3 = \bold M^3_{\lambda,\theta,\bar S}$ be the class of
$\bold m \in \bold M_2$ such that in addition:
\mn
\begin{enumerate}
\item[(g)]  
\begin{enumerate}
\item[$(\alpha)$]  the linear order $I$ is the disjoint sum of
  $\langle I_i:i < \theta,i \ne 0\rangle$, each $I_i$ of cardinality
  $\lambda$;
\sn
\item[$(\beta)$]  $\bold c$ is a function from $[I]^2$ to $\lambda$ such
  that $\bold c_{i,\alpha} = \bold c \rest [I_{i,\alpha}]^2$ for every
  $i \in S_2 \cup S_3$ and $\alpha < \lambda$;
\end{enumerate}
\sn
\item[$(h)$]  $\bold c$ witnesses
  $\Pr_0(\lambda,\lambda,\lambda,\aleph_0)$; see Definition \ref{w11}(1).
\end{enumerate}
\mn
3) Let $\bold M_{2.5} = \bold M^{2.5}_{\lambda,\theta,\bar S}$ be the class
of $\bold m \in \bold M_{1.5}$ such that in addition:
\mn
\begin{enumerate}
\item[(g),(h)]  as in part (2).
\end{enumerate}
\end{definition}

\noindent
The following definition is just to hint at what we need to get more
cardinals, not to be used in \S4.
\begin{definition}
\label{c28}
Let $\Pr_{2.5}(\lambda,\mu,\sigma,\partial,\theta)$ mean that $\theta
= \cf(\theta),\lambda \ge \mu,\sigma,\partial,\theta$ 
and some pair $(\bold c,\bar W)$ witness it, 
which means (if $\lambda = \mu$ we
may omit $\lambda$, if $\sigma = \partial \wedge \lambda = \mu$ then
we can omit $\sigma,\lambda$):
\mn
\begin{enumerate}
\item[$(a)$]  $\bar W = \langle W_i:i < \mu\rangle$ is a sequence of
  pairwise disjoint subsets of $\lambda$;
\sn
\item[$(b)$]  $\bold c:[\lambda]^2 \rightarrow \sigma$;
\sn
\item[$(c)$]  if $i < \mu$ is even and $\varp \in u_\varp \in
  [\lambda]^{< \partial}$ for $\varp \in W_i$ and $\gamma < \sigma$
  \then \, for some $\varp < \zeta < \lambda$ we have:
\sn
\begin{enumerate}
\item[$(\alpha)$]  $\varp \notin u_\zeta,\zeta \notin u_\varp$;
\sn
\item[$(\beta)$]  $\bold c\{\varp,\zeta\} = \gamma$;
\sn
\item[$(\gamma)$]  if $\xi_1 \in u_\zeta \backslash u_\varp$ and
  $\xi_2 \in u_\varp \backslash u_\zeta$ and $\{\xi_1,\xi_2\}
  \ne \{\varp,\zeta\}$ \then \, $\bold c\{\xi_1,\xi_2\} = 0$;
\sn
\item[$(\delta)$]  optional $(u_\varp,u_\zeta)$ is a $\Delta$-system
  pair (see proof);
\end{enumerate}
\sn
\item[$(d)$]  if $\langle \cU_\zeta:\zeta < \theta\rangle$ is
  $\subseteq$-increasing with union $W_i$ where $<\mu,i$ is 
odd then for some $\zeta < \theta$ we have
$\Rang(\bold c \rest [\cU_\zeta]^2) = \sigma$.
\end{enumerate}
\end{definition}

\begin{claim}
\label{c29}
1) Assume $\lambda = \lambda^{\langle \theta;\aleph_0\rangle}$, and
moreover $\lambda = \lambda^{\langle \theta;\theta\rangle}$, 
see Definition \ref{w22} recalling (see \ref{c23}) that 
$\theta = \cf(\theta) \in (\aleph_0,\lambda)$.
If $G \in \bold K_{\le\lambda}$, \then
\, there is $\bold m \in \bold M^2_{\lambda,\theta,\bar S}$ such that
$G_{\bold m,1} \cong G$.

\noindent
1A) If in addition $\Pr_0(\lambda,\lambda,\lambda,\aleph_0)$ or just
$\Pr_0(\lambda,\lambda,\aleph_0,\aleph_0)$ \then \, we
can add $\bold m \in \bold M^3_{\lambda,\theta,\bar S}$.

\noindent
1B) If $\lambda \ge 2^{\aleph_0}$ \then \, in part (1) we can
strengthen Definition \ref{c26} adding in clause 
$(e)(\varp)\bullet_1,\bullet_2$ that $v = \omega \backslash \{0\}$ hence
$\ell_0=1,\ell_1=2,\ldots$.

\noindent
2) Assume $\lambda = \lambda^{\langle \theta;\aleph_0\rangle}$ and
$\lambda = \lambda^{(\lambda;\theta)}$ \then \, there is $\bold m \in
\bold M^{1.5}_{\lambda,\theta,S}$ such that $G_{\bold m,1} \cong G$.

\noindent
3) In part (2), if in addition
$\Pr_0(\lambda,\lambda,\aleph_0,\aleph_0)$ then we can add $\bold m
\in \bold M^{2.5}_{\lambda,\theta,5}$.

\noindent
4) If $\lambda \ge \mu := \beth_\omega$ (or just $\mu$ strong limit)
\then \, for every large enough regular $\theta < \mu$, the
assumption of part (1) holds.

\noindent
5) If above $\theta = \aleph_1 < \lambda = \lambda^\theta$, \then
\, the assumption of part (1) holds.
\end{claim}

\begin{PROOF}{\ref{c29}}
See more details in the proof of \ref{i14}.

\noindent
1) For clause (d) of Definition \ref{c26} recall the proof of
\ref{b15} and the assumption $\lambda = \lambda^{\langle
  \lambda;\theta\rangle}$ (considering the phrasing of
\ref{c26}(d)$(\eta)$, why is not $\lambda =
\lambda^{\langle\lambda;\theta\rangle}$ enough?  Because in repeating
\S(3A) we have to get a set of pairwise commuting elements).  For
clause (e) of Definition \ref{c26} recall the assumption $\lambda =
\lambda^{\langle \theta;\aleph_0\rangle}$ and the proof of
\ref{b19}(1).

\noindent
1A) For clause (c) of Definition \ref{c25}, by
\cite[\S(1C),\S3]{Sh:312}  this should be clear.  

\noindent
1B) Should be similar.

\noindent
2),3)  Straightforward.

\noindent
4) By \cite{Sh:460} or see \cite[\S1]{Sh:829}.

\noindent
5) Check the definitions.
\end{PROOF}

\begin{remark}
\label{c33}
1) Concerning the use in \ref{c29}(2) of 
Definition \ref{w22} note that conceivably $(\forall \lambda >
2^{\aleph_0})(\exists \theta < \aleph_\omega)
(\aleph_0 <\theta = \cf(\theta) < \lambda
\wedge \lambda^{\langle \theta;\aleph_0\rangle} = \lambda \wedge
\lambda^{(\theta;\theta)} = \lambda)$, i.e. 
conceivably this is provable in $\ZFC$.

\noindent
2) Concerning \ref{c29} recall \ref{w14}.
\end{remark}

\begin{claim}
\label{c34}
Let $\bold m \in \bold M_1$.

\noindent
1) If $i<j \le \theta$ \then \, $G_{\bold m,i} \le_{\gS} G_{\bold
  m,j}$, see \ref{c23}(3).

\noindent
2) For every finite $A \subseteq G_{\bold m,\theta}$ there is a
sequence $\bar u = \langle u_i:i \in v\rangle$ such that:
\mn
\begin{enumerate}
\item[$(*)^1_{\bar u}$]
\begin{enumerate}
\item[(a)]   $v \subseteq \theta$ is finite and $0 \in v$ for
  notational simplicity;
\sn
\item[(b)]  $u_i \subseteq I_i$ is finite\footnote{Note that in
    \ref{c34}(2) we allow ``$u_i$ is empty".} for $i \in v$;
\sn
\item[(c)]  if $i \in v$, \then \, $\tp_{\qf}(\langle \bar b_{i,s}:
s \in u_i\rangle,G_i,G_\theta)$ does not split over
$\cup\{\bar b_{j,s}:j \in v \cap i$ and $s \in u_j\}$;
\sn
\item[(d)]   if $i \in S_1$ and $s \in u_i$ \then \,
$\bar a_{i,s} \subseteq \seb(\{\bar b_{j,s}:j \in v \cap i,s \in u_j\},G_i)$;
\sn
\item[(e)]   if $i \in S_2 \cup S_3$ and $s,t \in u_i$ \then \,
$\bar a_{i,\bold c\{s,t\}} \subseteq 
\seb(\{\bar b_{j,s}:j \in v \cap i,s \in u_j\},G_i)$;
\sn
\item[(f)]   if $A \subseteq G_{\bold m,i}$ and $i \in (0,\theta)$
 \then \, $v \subseteq i$;
\end{enumerate}
\sn
\item[$(*)_2$]  $A$ is included in $\seb(\{\bar b_{i,s}:i \in v,s \in
  u_i\},G_\theta)$.
\end{enumerate}
\mn
3) We have $\bar u = \langle u^1_i \cup u^2_i:i \in v\rangle$
satisfies $(*)_1$, i.e. $(*)^1_{\bar u}$ from part (2) holds \when \,: 
\mn
\begin{enumerate}
\item[$\oplus$]  
\begin{enumerate}
\item[(a)]  $\bar u_\ell = \langle u^\ell_i:i \in v\rangle$ for $\ell=1,2$;
\sn
\item[(b)]   we have $(*)^1_{\bar u_\ell}$ for $\ell=1,2$;
\sn
\item[(c)]  if $i \in v,s_1 \in u^1_i \backslash u^2_i$
  and $s_2 \in u^2_i \backslash u^1_i$ then $\bold c_i\{s_1,s_2\}=0$.
\end{enumerate}
\end{enumerate}
\mn
3A) If $v_1 \subseteq v_2,\bar u^2 = \langle u_i:i \in v_2\rangle,\bar
u^1 = \bar u^2 \rest v_1$ and $i \in v_2 \backslash v_1 \Rightarrow u_i
= \emptyset$ \then \, $(*)^1_{\bar u^1}
\Leftrightarrow (*)^1_{\bar u^2}$.

\noindent
4) The type $\tp_{\qf}(\langle \bar b^\ell_{i,s}:s \in u^\ell_i,\ell \in
\{1,2\}\rangle,G_i,G_{i+1})$ does not split over $\{\bar b^\ell_{j,s}:
j \in v \cap i,s \in u^\ell_j,\ell \in \{1,2\}\} \cup
\{a_{i,\alpha}\}$ \when \,:
\mn
\begin{enumerate}
\item[$(a)$]  $\bar u_\ell = \langle u^\ell_j:j \in v\rangle$;
\sn
\item[$(b)$]  $(*)^1_{\bar u_\ell}$ holds for $\ell =1,2$;
\sn
\item[$(c)$]  $i \in S_3 \cap v$;
\sn
\item[$(d)$]  $s_* \in u^1_i \backslash u^2_i,t_* \in u^2_i \backslash
  u^1_i$;
\sn
\item[$(e)$]  $\alpha = \bold c_i\{s_*,t_*\}$;
\sn
\item[$(f)$]  clause (c) from part (3) holds when $\{s_1,s_2\} \ne
  \{s_*,t_*\}$.
\end{enumerate}
\end{claim}

\begin{PROOF}{\ref{c34}}
1) By part (2). 

\noindent
2) By induction on $\min\{j < \theta:A \subseteq G_{\bold m,j}\}$.
Note that for $A \subseteq G_1$ clause $(*)^1_{\bar u}(c)$ is trivial.

\noindent
3),4) Easy, too.
\end{PROOF}

\begin{claim}
\label{c35}
If $\bold m \in \bold M_2$ or just $\bold m \in \bold M_{1.5}$, 
\then \, $G_{\bold m,\theta} \in \bold K^{\exlf}_\lambda$ 
is complete and is $(\lambda,\theta,\gS)$-full and
extend $G_{\bold m,1}$.
\end{claim}

\begin{PROOF}{\ref{c35}}
Being in $\bold K^{\lf}_\lambda$ is obvious as well as extending
$G_{\bold m,1}$; being
$(\lambda,\theta,\gS)$-full is witnessed by $\langle G_{\bold m,i}:i
<\theta\rangle,S_1$ being unbounded in $\theta$ and clauses
\ref{c25}$(c)$, \ref{c26}$(c)(\varp)$.

The main point is proving $G_{\bold m,\theta}$ is complete, so assume
$\pi$ is an automorphism of $G_{\bold m,\theta}$.

Now
\mn
\begin{enumerate}
\item[$(*)_1$]  if $\partial = \theta^+,i < \theta$ is a limit ordinal
  and $j \in S_2 \backslash (i+2)$, then $G_i$ is
  $\theta$-indecomposable inside $G_{j+1}$.
\end{enumerate}
\mn
[Why?  Toward contradiction assume $\langle H_\varp:\varp <
\theta\rangle$ is increasing with union $G_{j+1}$ but $\varp
< \partial \Rightarrow G_i \nsubseteq H_\varp$.  Recall that by
\ref{c23}(4), $i = \sup(S_1 \cap i)$ hence by \ref{c25}(c),
\ref{c26}(c) we have $G_i \in \bold K_{\exlf}$.  Also there is
$I^\bullet \subseteq I_{j+1}$ of cardinality $\lambda$ such that $s
\in I^\bullet \Rightarrow \gs_{j+1,s} = \gs_{\cg}$.  Hence $\langle
\square_{b_{j,s}}(G_i):s \in I^\bullet\rangle$ is a sequence of
pairwise commuting subgroups of $G_{i+1}$.  
For each $s \in I^\bullet,b_{j,s}$ belongs to
$H_{\varp(s)}$ for some $\varp(s) < \theta$, hence
$\langle \square_{b_{j,s}}(G_i) \cap H_\varp:\varp < \theta\rangle$
is not eventually constant and choose $c_{s,\varp} \in
\square_{b_{j,s}}(G_i)) \backslash H_\varp$ of order 2 for $\varp
< \partial$.  

As $|I^\bullet| \ge \theta$ we can find pairwise distinct $s(\varp) \in
I^\bullet$ for $\varp < \theta$.  Hence $\langle
c_{s(\varp),\varp}:\varp < \theta\rangle$ is a sequence of members of
$G_{j+1}$, pairwise commuting (recall the choice of $I^\bullet$!) each
of order two and independent.  By \ref{c26}(d)$(\eta)$ there is
$\alpha < \lambda$ such that $A = \{c_{s(\varp),\varp}:\varp <
\theta\} \cap \{a_{j,\gamma}:\gamma \in J_{j,\alpha}\}$ has
cardinality $\theta$.  By \ref{c26}(d)$(\varp)$ applied on the pair
$(j,\alpha)$, the function $\bold c_{j,\alpha}$ from
$[I_{j,\alpha}]^2$ into $J_{j,\alpha}$ is $\theta$-indecomposable.

For $\varp < \theta$ let $I_{j,\alpha,\varp} = \{s \in
I_{j,\alpha}:b_{j,s} \in H_\varp\}$, so $\langle
I_{j,\alpha,\varp}:\varp < \theta\rangle$ is $\subseteq$-increasing
with union $I_{j,\alpha}$ hence for some $\varp(*) < \theta$, the set
$\{\bold c_j\{s,t\}:s \ne t \in I_{j,\alpha,\varp(*)}\}$, in fact, is
equal to $J_{j,\alpha}$.  Now the set $\cX = \{[b_{j,s},b_{j,t}]:s \ne t
\in I_{j,\alpha,\varp(*)}\}$ is included in $H_{\varp(*)}$ by the
choice of $I_{j,\alpha,\varp(*)}$.  Moreover, recalling $[b_{j,s},b_{s,t}] =
a_{j,\bold c_j\{s,t\}}$ and the choice of $\varp(*)$, the set $\cX$ 
includes a subset
of $\{c_{s(\varp),\varp}:\varp < \theta\}$ is of cardinality
$\theta$.  But this contradicts the choice of the
$c_{s(\varp),\varp}$'s.  So $G_i$ is indeed $\theta$-indecomposable
inside $G_{j+1}$.]

So $\langle \pi(G_{\bold m,i}):i < \theta\rangle$ is $\le_{\bold
  K_{\lf}}$-increasing with union $G_{\bold m,\theta}$ hence by
$(*)_1 $ above, if $i \in S_2$ then $(\forall^\infty j<\theta)(G_{\bold
  m,i} \subseteq \pi(G_{\bold m,j}))$.  The parallel statement holds
for $\pi^{-1}$ hence $E$ is a club of $\theta$ where $E := \{i <\theta:i$ 
is a limit ordinal, hence $i=\sup(S_1 \cap i)$ and $\pi$ maps 
$G_{\bold m,i}$ onto $G_{\bold m,i}\}$; note that by the middle demand,
$i \in E \Rightarrow G_i \in \bold K_{\exlf}$.

Next we define:
\mn
\begin{enumerate}
\item[$(*)_2$]  $S^\bullet$ is the set of $i \in E \cap S_1$ such that $\pi$ 
is not the identity on $\bold C'(G_{\bold m,i},G_{\bold m,i + \omega})$.
\end{enumerate}
\medskip

\noindent
\underline{Case 1}:  $S^\bullet$ is unbounded in $\theta$

So for $i \in S^\bullet$ choose $c_i \in \bold C'(G_{\bold m,i},
G_{\bold m,i+1})$ such that $\pi(c_i) \ne c_i$.  \Wilog \, $c_i$ has
roder 2, because the set of elements of order 2 
from $\bold C(G_{\bold m,i},G_{\bold m,i+ \omega})$ 
generates it, see \cite[4.1=Ld36,4.10=Ld93]{Sh:312}.  Choose $\langle
\bold i_\varp = \bold i(\varp):\varp < \theta\rangle$ increasing, 
$\bold i_\varp \in S^\bullet$
and so as $\bold i_\varp + \omega \le \bold i_{\varp +1} \in E$ clearly
 $\pi(c_\varp) \in G_{\bold m,\bold i(\varp+1)}$.  Now we apply
\ref{c26}(e), \ref{w11}(1) and get contradiction by \ref{c34}(4) recalling
\ref{c26}(2)(h) and \ref{c25}(e); but we elaborate.

\noindent
Now we apply \ref{c26}(1)(e) (indirectly \ref{c29}(1), \ref{w22}).  So
there are $(i,\alpha,v,\ell_0,\ell_1,\ldots,\varp_0,\varp_1,\ldots)$ as
there, in particular $i \in S_3$ and here $v = \omega \backslash
\{0\}$.  Now for every $s \in I_{i,\alpha}$
we apply \ref{c34}(2), getting $\bar u_s = \langle u_{s,\iota}:\iota
\in v_s\rangle$ and let $\ell_s$ be such that $v_s \subseteq
j_{i,\omega \alpha + \ell_s}$, \wilog \, $i \in v_s,s \in u_{s,i}$.

Now consider the statement:
\mn
\begin{enumerate}
\item[$(*)_3$]  there are $s_1 \ne s_2 \in I_{i,\alpha}$ and $k$ such
  that:
\sn
\begin{enumerate}
\item[(a)]  $\bold c\{s_1,s_2\} = \ell_k$;
\sn
\item[(b)]  $\ell_k > \ell_{s_1},\ell_{s_2}$;
\sn
\item[(c)]  if $t_1 \in \cup\{u_{s_1,\iota}:\iota \in v_{t_1}
  \backslash i\},t_2 \in \cup\{u_{s_2,\iota}:\iota \in v_{t_2}
  \backslash i\}$ and $\{t_1,t_2\} \ne \{s_1,s_2\}$ 
\then \, $\bold c\{t_1,t_2\}=0$;
\newline
or for later proofs:
\sn
\item[(c)$'$]
\begin{enumerate}
\item[$(\alpha)$]  if $t_1 \in u_{s_1,i} \backslash u_{s_2,i}$ and
  $t_2 \in u_{s_2,i} \backslash u_{s_1,i}$ and
\begin{itemize}
\item  $\{t_1,t_2\} \ne \{s_1,s_2\}$ then $\bold c\{t_1,t_2\} = 0$,
  \oor \, just
\sn
\item  $t_1,t_2 \in I_{i,\alpha} \Rightarrow \bold c\{t_1,t_2\} <
  \ell_k$;
\sn
\item  $t_1,t_2 \in I_{i,\beta},\beta < \lambda;\beta \ne \alpha$ then
  $j_{i,\omega \beta + \bold c\{t_1,t_2\}} < j_{i,\omega \alpha +
    \ell(k)}$ (we use $j_{i,\omega \alpha + \ell} \in
  (j^*_{i,\ell},j^*_{i,\ell +1})$ - check);
\end{itemize}
\sn
\item[$(\beta)$]  if $\iota \in v_1 \cap v_2$ and $\iota > i, (\iota
  \in S_3),\beta < \lambda$ and $t_1 \in v_{s_1,\iota},t_2 \in
  v_{s_2,\iota}$ then $\bold c\{t_1,t_2\} = 0$.
\end{enumerate}
\end{enumerate}
\end{enumerate}
\mn
Now why is $(*)_3$ true?  This is by the choice of $\bold c$, that is,
as $\bold c$ exmplifies $\Pr_0(\lambda,\lambda,\lambda,\aleph_0)$ (in
later proofs\footnote{Notice we have used $\lambda^{(\theta;\aleph_0)}
  = \lambda$ then $\Pr_*$ from Definition \ref{i5} suffice, whereas if
  we have used $\lambda = \lambda^{\langle \theta;\aleph_0\rangle}$
  then $\Pr_{**}$ is needed because this also has to fit the version
  of $\bold M$ we use.} we use less).

Now to get a contradiction we would like to prove:
\mn
\begin{enumerate}
\item[$(*)_4$]  the type $\tp((\pi(b_{s_1}),\pi(b_{s_2})),G_{\bold m,i},
G_{\bold m,\theta})$ does not split over 
$G_{\bold m,j_{i,\omega \alpha + \ell(k)}} \cup \{c_{\bold i(\varp_k)}\}$
hence over $G_{\bold m,\bold i(\varp_k)} \cup \{c_{\bold i(\varp(k))}\}$.
\end{enumerate}
\mn
It follows from $(*)_4$ that $\tp((b_{s_1},b_{s_2}),\pi^{-1}(G_{\bold m,i}),
\pi^{-1}(G_{\bold m,\theta}))$ does not split over
$\pi^{-1}(G_{m,\bold i(\varp_k)}) \cup \{\pi^{-1}(c_{\bold i(\varp)})\})$.  But
$i(\varp_k),i \in E$ have it follows that $\pi(G_{m,i}) = G_{\bold
  m,i}$ and $\pi^{-1}(G_{\bold i(\varp_k)} = G_{\bold i(\varp_k)})$ has
$\tp((b_{s_1},b_{s_2}),G_{\bold m,i},G_{\bold m,\theta})$ does not
split over $G_{\bold i(\varp_k)} \cup\{\pi^{-1}(c_{\bold i(\varp_k)})\}$.

Now $[b_{s_1},b_{s_2}] = \pi^1([b_{s_1},b_{s_2}]) = 
\pi^{-1}(c_{\bold i(\varp_k)})$
which is $\ne c_{i(\varp_k)}$. But as $c_{\bold i(\varp_k)} \in \bold
C(G_{\bold m,\bold i(\varp_k)},G_{\bold m,\theta})$ clearly also
$\pi^{-1}(c_{\bold i(\varp_k)})$ belongs to it, hence it follows that
$\pi^{-1}(c_{\bold i(\varp_k)}) \in \seb(\{c_{\bold i(\varp_k)}\};G_\theta)$,
but as $c_{\bold i(\varp_k)}$ has order two, the latter belongs to
$\{c_{\bold i(\varp_k)},e_{G_\sigma}\}$.

However $\pi^{-1}(c_{\bold i(\varp_k)})$ too has order 2 hence is equal to
$c_{\bold i(\varp_k)}$; applying $\pi$ we get $c_{\bold i(\varp_k)} =
\pi(c_{\bold i(\varp_k)})$ a contradiction to the choice of the $c_i$'s.

[Pre 16.11.11 proof: as 
$c_{\bold i(\varp_k)} \ne G_{\bold m,j_i,\omega \alpha + \ell(k)}$ (see
above), contradiction.

Why does it hold?  By \ref{c34}(3) and $(*)_3(c)$, in later proofs by
finer versions (for more cardinals) we have to do more.]
\medskip

\noindent
\underline{Case 2}:  $i_* = \sup(S^\bullet)+1$ is $< \theta$.

Now for any $i \in S' := E \cap S_1 \backslash i_*$ 
by \cite[2.18=Lc62]{Sh:312} there is $g_i \in G_{\bold m,i+1}$ 
such that $\square^{g_i}(G_{\bold m,i}) \subseteq 
\bold C(G_{\bold m_i},G_{\bold m,i+1})$.  So if $a \in 
G_{\bold m,i}$ then $g^{-1}_i a g_i \in \bold C'(G_{m,i},G_{\bold
  m,i+1})$ and $a = g_i(g^{-1}_i a g_i)g^{-1}_i$ 
hence $\pi(a) = \pi(g_i) \pi(g^{-1}_i a
g_i) \pi(g^{-1}_i) = \pi(g_i)(g^{-1}_i a g_i) 
\pi(g_i)^{-1}$ recalling $i \notin
S^\bullet$ being $ \ge i_*$ hence $\pi(a) = (g_i \pi(g_i)^{-1})^{-1} a
g_i \pi(g^{-1}_i)$.  If for some $g$ the
set $\{i \in S':g_i=g\}$ is unbounded in $\theta$ we are easily
done, so toward contradiction assume this fails. 

But for every $\delta \in \acc(E) \cap S_1 \backslash i_*$, we can by
\ref{c34}(1)  choose a finite $\bar a_\delta \subseteq G_\delta$ and
$\gs_\delta \in \gS$ such that $\tp_{\bs}(\pi(g_\delta)g^{-1}_\delta,
G_\delta,G_\theta) 
= q_{\gs_\delta}(\bar a_\delta,G_\delta)$ and let
$i(\delta) \in E \cap \delta$ be such that $\bar a_\delta \subseteq
G_{i(\delta)}$.  

Clearly:
\mn
\begin{enumerate}
\item[$\circledast$]  if $d_1,d_2 \in G_\delta,d_2 \ne \pi(d_1)$ then
  $\tp_{\bs}(\langle d_1,d_2\rangle,\bar a_\delta,G_\delta) \ne
\tp_{\bs}(\langle d_1,\pi(d_1)\rangle,\bar a_\delta,G_\delta)$.
\end{enumerate}
\mn
[Why?  Because $\pi(d_1) = \pi(g_\delta) g^{-1}_i d_1 g_i \pi(g_\delta)^{-1}$
 and the choice of $\bar a_\delta$.]

Hence for some group term $\sigma_{d_1}(\bar x_{1+\ell g(\bar
  b_\delta)})$ we have $\pi(d_1) =
\sigma^{G_\delta}_{d_1}(d_1,\bar a_\delta)$ and $\sigma_{d_1}$ depends
only on $\tp_{\bs}(d_1,\bar a_\delta,G_\delta)$.  By Fodor Lemma for
some $i(*)$ the set $S = \{\delta:\delta \in \acc(E) \cap S_1 
\backslash i_*$ and
$i(\delta) = i(*)\}$ is a stationary subset of $\theta$.

Now we can finish easily, e.g. as $G_\delta$ for $\delta \in S$
belongs to $\bold K_{\exlf}$ and we know that it can be
extended to a complete $G' \in \bold K_{\exlf}$ or just see that all the
definitions in $\circledast$ agree and should be one conjugation.
\end{PROOF}

\begin{conclusion}
\label{c41}
Assume $\lambda > \beth_\omega$ is a successor of a regular and $G
\in \bold K^{\lf}_{\le \lambda}$ and $\theta = \cf(\theta) \in
(\aleph_0,\beth_\omega)$ is large enough and $\gS$ is as in \ref{c23}(3).

\Then \, there is a complete $(\lambda,\theta,\gS)$-full $H \in \bold
K^{\exlf}_\lambda$ extending $G$.
\end{conclusion}

\begin{PROOF}{\ref{c41}}
Fixing $\lambda$ and $\theta$ by \ref{c35} it suffices to find $\bold
m \in \bold M^3_{\lambda,\theta}$ such that $G_{\bold m,1} = G$.  As
$\lambda \ge \beth_\omega$, the assumption of \ref{c29}(1) holds for
every sufficiently large $\theta < \beth_\omega$; hence
there is $\bold m \in \bold M^2_{\lambda,\theta,\bar S}$ such that
$G_{\bold m,1}$ is isomorphic to $G$ and $\bar S$ as there.

As $\lambda$ is a successor of a regular, the assumption of
\ref{c29}(1A) holds (by \ref{w14}) hence $m \in \bold
M^3_{\lambda,\theta,\bar S}$.  So by \ref{c35} we indeed are done.
\end{PROOF}

\begin{remark}
\label{c44}
The assumption ``$\lambda > \beth_\omega$" comes from quoting
\ref{c29}(2) hence it is ``hard" for $\lambda < \beth_\omega$ to fail.
Similarly below.
\end{remark}
 
\noindent
Of course we have:
\begin{observation}
\label{c38}
If $\bold m \in \bold M_{1.5}$ \then \, $G_{\bold m,\theta}$ is
$(\lambda,\theta,\gS)$-full and extends $G_{\bold m,0}$.
\end{observation}
\newpage

\section {More uncountable cardinals} \label{5}
\bigskip

\subsection {Regular $\lambda$} \label{5A}
\bigskip

\begin{discussion}
\label{n01}
1) [Saharon: beginning of explanation of what we intend to do was lost
- FILL.]

\noindent
2) Clearly this helps because if $\bold c:"[\lambda]^2 \rightarrow
\lambda$ is $\theta$-indecomposable by \S(3A) we can find $G_3 \in
\bold K_{\lf}$ extending $G_2$ such that for each $\ell < 2$, the
subgroup $\seb(\{b_{\alpha,\ell}:\alpha < \lambda\},G_2)$ is
$\theta$-indecomposable in $G_3$ and this implies $G \in
\seb(\{a_\alpha:\alpha < \lambda\},G_1)$ is $\theta$-indecomposable in
$G_3$.

\noindent
3) How shall we prove this?  We use \cite[4.1=Ld36]{Sh:312} and more
in \cite[\S(4A)]{Sh:312}.  For the proof in \S4 to work, we need also
that $G_1 \le_{\gS} G_3$, but the proof in \S4 seems to use more.
That is, assume that in Definition \ref{c23} we have $\bar S = \langle
S_\ell:\ell=1,\dotsc,4\rangle$ for $\delta \in S_4$ we apply the
above.  We need that the case $i \in S_4$ will be similar to $i \in
S_2$, in the sense of putting together the $\tp(\bar b_{i,s},G_i,G_{i+1})$'s 
for $s \in I_1$.

\noindent
4) Why have we not used the above in \S4?  Note that for the final
result in \S4, this does not make a real change and may well be for
the other universe $\bold K$, the proof in \S4 works but not the proof
here.

\noindent
5) Now in the results of \S4, if $\aleph_0 < \theta = \cf(\theta) \le
\lambda$, we really need just $\lambda =
\lambda^{\langle \theta;\aleph_0\rangle}$.  
The point is the above helps to prove
$E = \{i:\pi(G_{\bold m,i}) = G_{\bold m,i}\}$ is a club of $\theta$.
We still have to use $\lambda = \lambda^{\langle
  \theta;\aleph_0\rangle}$.  Also the use of
$\Pr_0(\lambda,\lambda,\aleph_0,\aleph_0)$ is as before replacing it
by $\Pr^*_0(\lambda,\lambda,\aleph_0,\aleph_0)$.
\mn
\begin{enumerate}
\item[(A)]  Restricting ourselves to $\lambda = \cf(\lambda) >
  \aleph_0$, we can partition $S^\lambda_{\aleph_0}$ into
  $\lambda$-stationary sets $\langle S^+_\varp:\varp < \lambda \rangle$,
  and $\bar C = \langle C_\delta:\delta \in S^\lambda
  _{\aleph_0}\rangle$ satisfies $C_\delta \subseteq \delta =
  \sup(C_\delta),\otp(C_\delta) = \omega$ and each $\bar C \rest
  S^*_\varp$ guesses clubs and let $\bold c:[\lambda]^2 \rightarrow
  \aleph_0$ be $\bold c(\alpha,\beta) = \otp(C_\beta \cap \alpha)$ for
  $\alpha < \beta < \lambda$.
\sn
\item[(B)]  if $\lambda$ is singular, $\mu = \mu^{\aleph_0} < \lambda
  \le 2^\mu$ then as in \cite{Sh:331} but simpler, we use $S_\alpha
  \subseteq S^{\mu^+}_{\aleph_0}$ stationary for $\alpha < \lambda$,
  no one inducted in the countable union of others.
\sn
\item[(C)]  $\lambda$ singular, $\lambda = \sup\{\chi:\cf(\chi) =
  \aleph_0\}$ work as in \cite{Sh:331}.
\sn
\item[(D)]  $\cf(\chi) = \aleph_0 < \chi$ ... FILL.
  \end{enumerate}
\end{discussion}

\begin{discussion}
\label{i2}
(17.05.31)

\noindent
1) Do we need in \S4 the property $\Pr_*$ or $\Pr_{**}$ defined in
\ref{i5}, \ref{i29} below?  We need it in Case 1 of the proof.  So if
we have used $\lambda^{(\theta;\aleph_0)} = \lambda$, \then \, $\Pr_*$
suffices (as we really ``know" the sequence of elements used for
$I_{i,\alpha},i \in S_3$ in but if we use $\lambda^{\langle
  \theta;\aleph_0\rangle} = \lambda$.

This reflects on the choice of the version of $\bold M$ and proof of
completeness. 

\noindent
2) It is nicer if we allow $\sigma < \theta < \lambda$.  For this we
may have $\bold c:[\lambda]\rightarrow \lambda$ and for the demands on
$\sigma$ we ``ignore" values of $\bold c$ which are $\ge \sigma$, etc.

\noindent
3) If Definition \ref{i5}, \ref{i29} it is not enough in
$(d)(\beta)\bullet_2$ to demand only
\mn
\begin{itemize}
\item  $\bold c\{\varp_1,\zeta_1\} \ne j$.
\end{itemize}
\mn
That is, is it not enough in the proof of completeness in \S4.

\noindent
4) See \S(5B) - trying to use \cite{Sh:331}.
\end{discussion}

\noindent
Central here (see a relative $\Pr_{**}$ in \ref{i29}).
\begin{definition}
\label{i5}
(was \ref{n05})

Assume $\lambda \ge \mu \ge \kappa + \theta_0 + \theta_1,\kappa$ regular,
$\bar\theta = (\theta_0,\theta_1)$, if
$\theta_0 = \theta$ we may write $\theta_0$ and then below \wilog \,
$u_{\varp,0} = u_{\varp,1}$; if we omit $\kappa$ we mean $\kappa = \sigma$.

Let $\Pr_*(\lambda,\mu,\sigma,\bar\kappa,\theta)$ mean that there is $\bold
c:[\lambda]^2 \rightarrow \sigma$ witnessing it, which means:
\mn
\begin{enumerate}
\item[$(*)^2_{\bold c}$]  for some $\bar S$ we have:
\sn
\begin{enumerate}
\item[(a)]  $\bar S$ is a sequence of the form $\langle
  (S_\alpha,\delta^*_\alpha):\alpha < \lambda\rangle$, the $S_\alpha$
  pairwise disjoint:
\sn
\item[(b)]  $\delta^*_\alpha = \sup(S_\alpha)$ and 
$\cf(\delta^*_\alpha) \in (\aleph_0,\lambda]$;
\sn
\item[(c)]  $S_\alpha$ is stationary;
\sn
\item[(d)]  if $(\alpha)$ then $(\beta)$, where:
\sn
\begin{enumerate}
\item[$(\alpha)$]  $\alpha < \lambda,S \subseteq
  S_\alpha$ is stationary in $\delta^*_\alpha = \sup(S_\alpha)$
 and $\langle u_{\varp,0},u_{\varp,1}:\varp \in 
S\rangle$ satisfies $u_{\varp,0} \cap u_{\varp,1} = 
u_{\varp,\iota} \subseteq [\lambda]^{< \kappa_\iota}$ for
$\varp \in S,\iota < 2$; if $u_{\varp,0} = u_{\varp_i}$ we may write
$u_{\varp,\iota} = u_\varp$;
\sn
\item[$(\beta)$]  for some $A \subseteq \sigma$ of cardinality $<
  \kappa_0 + \kappa_1$, arbitrarily large $j < \sigma$ there are
  $\varp < \zeta$ from $S$ such that such that:
\sn
\begin{enumerate}
\item[$\bullet_1$]    $\bold c\{\varp,\zeta\} = j$
\sn
\item[$\bullet_2$]    if $\alpha \in u_{\varp,0} \backslash (u_{\zeta,0} \cup
  u_{\zeta,1})$ and $\beta \in u_{\zeta,1} \backslash (u_{\varp,0} \cup
  u_{\varp,1})$ and $\{\alpha,\beta\} \ne \{\varp,\zeta\}$ then 
$\bold c\{\alpha,\beta\}  \in A$
\newline
(earlier version: $(\exists j_*)(\forall j > j_*)$ and in $\bullet_2$
only $\bold c\{\alpha,\beta\} < j$
\end{enumerate}
\end{enumerate}
\sn
\item[(e)]  $\bold c \rest S_\alpha$ is $\theta$-indecomposable
  regular $\kappa \notin \{\sigma,\lambda\}$.
\end{enumerate}
\end{enumerate}
\end{definition}

\begin{claim}
\label{i8}
If $\lambda = \cf(\lambda) > \sigma^+,\lambda > \theta = \cf(\theta)
\ne \sigma = \cf(\sigma)$ and $\varp < \lambda \Rightarrow |\varp|^{<
  \sigma} < \lambda$, \then \, $\Pr_*(\lambda,\lambda,\sigma,\theta)$.
\end{claim}

\begin{PROOF}{\ref{i8}}
Let $\langle S_\alpha:\alpha < \lambda\rangle$ be a partition of
$S^\lambda_\sigma := \{\delta < \lambda:\cf(\delta) = \sigma\}$ to
stationary sets.  By \cite[Ch.III]{Sh:g} we can find $\bar C$ such
that:
\mn
\begin{enumerate}
\item[$(*)_1$]  
\begin{enumerate}
\item[(a)]  $\bar C = \langle C_\delta:\delta \in
  S^\lambda_\sigma\rangle$;
\sn
\item[(b)]  $C_\delta \subseteq \delta = \sup(C_\delta)$ and
  $\otp(C_\delta) = \sigma$ if $\theta \le \sigma$, and the ordinal
  product $\theta,\sigma$ othewise (check proof of this version);
\sn
\item[(c)]  for each $\alpha < \lambda,\bar C \rest S_\alpha$ guesses clubs.
\end{enumerate}
\end{enumerate}
\mn
We define $\bold c:[\lambda]^2 \rightarrow \sigma$ by:
\mn
\begin{itemize}
\item  if $\varp < \zeta$ are from $S_\alpha$ then $\bold
  c\{\varp,\zeta\} = \otp(C_\zeta \cap \varp)$;
\sn
\item  if $\varp < \zeta < \lambda$ and $\neg(\exists
  \alpha)(\{\varp,\zeta\} \subseteq S_\alpha)$ then $\bold
  c\{\varp,\zeta\} = 0$.
\end{itemize}
\mn
Clearly in Definition \ref{i5}, Clauses (a),(b),(c) holds; to prove
clause (d) assume:
\mn
\begin{enumerate}
\item[$(*)_2$]  $\alpha < \lambda$ and $S \subseteq S_\alpha$ is
  stationary and $u_\varp \in [\lambda]^{< \theta}$ for $\varp \in S$.
\end{enumerate}
\mn
We should find $j_*$, etc. as in the definition. 
For each $\varp \in W$ let $\xi^1_\varp = \sup(u_\varp \cap \varp)$
and $\xi^2_\varp = \sup\{C_\zeta \cap \varp:\zeta \in u_\varp \cap
S^\lambda_\sigma \backslash \{\varp\} +1$.  Now $\xi^1_\varp < \varp$
because $|u_\varp| < \theta \le \sigma = \cf(\sigma)$ and $\xi^2_\varp
< \varp$ because $\zeta \in S^\lambda_\sigma \Rightarrow \otp(C_\zeta)
= \sigma \wedge \sup(C_\zeta) = \zeta$ hence $\zeta \in
S^\lambda_\sigma \backslash \{\varp\} \Rightarrow \sup(C_\zeta \cap
\varp) < \varp$.

For $\varp \in S$ let $j_\varp$ be $\sup[\{\bold c\{\zeta_1,\zeta_2\}:
\zeta_1 \ne \zeta_2,u_\varp\}$, it is $< \varp$
for some $\xi_* < \lambda$ and $j_* < \sigma$ the set $S' := \{\varp
\in W:\xi^1_\varp,\xi^2_\varp \le \xi_*\}$ is stationary.  Recalling
$\sigma \ge \theta,|v_\varp| < \sigma$ and $|\xi_*|^{< \sigma} <
\lambda$ then for some $u_*$ the set $S'' := \{\varp \in W_1:u_\varp
\cap u_*$ and $\zeta < \varp \Rightarrow u_\zeta \subseteq \varp\}$ is
stationary. 

Let $E = \{\delta < \lambda;\delta = \sup(S'' \cap \delta)\}$, clearly
a club of $\lambda$, hence we can choose $\zeta \in E \cap S''$.  We
shall show that:
\mn
\begin{enumerate}
\item[$(*)$]  for arbitrarily large $j < \sigma$ there is $\varp \in
  S'' \cap \zeta$ such that $(j_*,j,\varp,\zeta)$ is as required in
  $(d)(\beta)$ of Definition \ref{n05}, and more.  
\end{enumerate}
\mn
This clearly suffices; but why it holds?  Let $j_0 < \sigma$.  As
$\otp(C_\zeta) = \sigma$ there is $\varp_1 \in C_\zeta$ such that
$\otp(C_\zeta \cap \varp_1) > j_0$, so as $\zeta \in E$ and the choice
of $E$ there is $\varp \in S'' \cap [\varp_1,\zeta)$.

Now
\mn
\begin{enumerate}
\item[$\oplus_{2.1}$]  $\bold c\{\varp,\zeta\} = \otp(C_\zeta \cap
  \varp) > \otp(C_\zeta \cap \varp_1) > j_0$
\sn
\item[$\oplus_{2.2}$]  assume $\alpha \in u_\varp \backslash
  u_\zeta,\beta \in u_\zeta \backslash u_\varp$ and $\{\alpha,\beta\}
  \ne \{\varp,\zeta\}$, \then \, $\bold c\{\alpha,\beta\} = 0$.
\end{enumerate}
\mn
[Why?  By the choices of $S',S''$.]
\mn
\begin{enumerate}
\item[$\oplus_{2.3}$]  $\alpha \ne \beta \in u_\varp$ or $\alpha \ne
  \beta \in u_\zeta$ then $\bold c\{\alpha,\beta\} < j_*$.
\end{enumerate}
\mn
[Why?  By the choice of $S'$.]

This is more than required.

We are left with clasue
\mn
\begin{enumerate}
\item[$(*)_3$]  assume $\alpha < \lambda,\kappa = \cf(\kappa) \notin
  \{\sigma,\lambda\},\langle S_{\alpha,i}:i < \kappa\rangle$ is
  $\subseteq$-incraesing with union $S_\alpha$.
\end{enumerate}
\mn
We should prove that for some $i < \kappa,\Rang(\bold c \rest
[S_{\alpha,i}]^2) = \sigma$.  Now if $\kappa > \lambda$ necessarily
$S_{\alpha,i} = S_\alpha$ for $i < \kappa$ large enough so \wilog \,
$\kappa < \lambda_i$.

Now
\mn
\begin{enumerate}
\item[$\oplus_{3.1}$]  for some $i < \kappa,\bar C \rest S_{\alpha,i}$
  guess clubs, i.e. $S_{\alpha,i}$ is stationary and for every club
  $E$ of $\lambda$, for some $\zeta \in S_{\alpha,i}$ we have $C_\zeta
  \subseteq E$ (e.g. for stationarily many $\zeta \in S_{\alpha,i}$).
\end{enumerate}
\mn
[Why?  See \cite[Ch.III]{Sh:g} or just if $E_i$ is a counter-example
for $S_{\alpha,i}$, then $\bigcap\limits_{i < \kappa} E_i$ is a club
of $\lambda$ such that $\zeta \in S_\alpha \Rightarrow
\bigvee\limits_{i < \kappa} \zeta \in S_{\alpha,1} \Rightarrow
\bigvee\limits_{i < \kappa} (\zeta \in S_{\alpha,i} \wedge C_\zeta
\nsubseteq E_i) \Rightarrow C_\zeta \nsubseteq E_i$, contradicting the
choice of $\bar C$ in $(*)_1$.]

Having chosen $i < \kappa$ in $\oplus_{3.1}$, let $E = \{\delta <
\lambda:\delta = \sup(S_{\alpha,i} \cap \delta)\}$ hence there is
$\zeta \in S_{\alpha,i}$ such that $C_\zeta \subseteq E$.  Letting
$\langle \beta_{\zeta,j}:j < \sigma\rangle$ list $C_\zeta$ is an
increasing cardinal, clearly there are $\varp_j$ for $j < \sigma$ such
that $\varp_j \in [\beta_{\zeta,j},\beta_{\zeta,j+1}) \cap
S_{\alpha,i}$.]

So
\mn
\begin{enumerate}
\item[$\oplus_{3.2}$]  if $j < \sigma$, then $\varp_j < \zeta$ are
  from $S_{\alpha,i} \subseteq S_\alpha$ and $\bold
  c\{\varp_j,\zeta\},\otp(\varp_j \cap \zeta)-1=(j+1)-1 = j$.
\end{enumerate}
\end{PROOF}

\begin{remark}
\label{i11}
1) For $\lambda$ singular: try as in \cite[Ch.V]{Sh:e} - divide to
coset?, i.e. $\lambda = \mu^{+ \omega}$ or $\lambda = \aleph_\delta
\rest \omega^2|\delta$ or ?

\noindent
2) Maybe weaken $\Pr_*$ as in \cite[Ch.V]{Sh:e}, i.e. not ``for every
$\alpha < \lambda$ and stationary $S \subseteq S_\alpha$" but ``for
every $\alpha$ if $S = S_\alpha$ (but maybe more than the $u_\alpha$'s).
\end{remark}

\begin{claim}
\label{i14}
(was n09) In Claim \ref{c29}, we can replace
$\Pr_0(\lambda,\lambda,\lambda,\aleph_0,\aleph_0)$ by
$\Pr_*(\lambda,\lambda,\aleph_0,\theta)$.

\noindent
1) Assume $\lambda = \lambda^{\langle \theta;\aleph_0\rangle}$, and
moreover $\lambda = \lambda^{\langle \theta;\theta\rangle}$, 
see Definition \ref{w22} recalling (see \ref{c23}) that 
$\theta = \cf(\theta) \in (\aleph_0,\lambda)$.
If $G \in \bold K_{\le\lambda}$, \then
\, there is $\bold m \in \bold M^2_{\lambda,\theta,\bar S}$ such that
$G_{\bold m,1} \cong G$.

\noindent
1A) If in addition $\Pr_*(\lambda,\lambda,\aleph_0,\theta)$ or just
$\Pr_0(\lambda,\lambda,\aleph_0,\aleph_0)$ \then \, we
can add $\bold m \in \bold M^3_{\lambda,\theta,\bar S}$.

\noindent
1B) If $\lambda \ge 2^{\aleph_0}$ \then \, in part (1) we can
strengthen Definition \ref{c26} adding in clause 
$(e)(\varp)\bullet_1,\bullet_2$ that $v = \omega \backslash \{0\}$ hence
$\ell_0=1,\ell_1=2,\ldots$.

\noindent
2) Assume $\lambda = \lambda^{\langle \theta;\aleph_0\rangle}$ and
$\lambda = \lambda^{(\lambda;\theta)}$ then there is $\bold m \in
\bold M^{1.5}_{\lambda,\theta,S}$ such that $G_{\bold m,1} \cong G$.

\noindent
3) In part (2), if in addition
$\Pr_*(\lambda,\lambda,\aleph_0,\theta)$ then we can add $\bold m
\in \bold M^{2.5}_{\lambda,\theta,5}$.

\noindent
4) If $\lambda \ge \mu := \beth_\omega$ (or just $\mu$ strong limit)
\then \, for every large enough regular $\theta < \mu$, the
assumption of part (1) holds.

\noindent
5) If above $\theta = \aleph_1 < \lambda = \lambda^\theta$, \then
\, the assumption of part (1) holds.
\end{claim}

\begin{PROOF}{\ref{i14}}
Let
\mn
\begin{enumerate}
\item[$(*)_1$]
\begin{enumerate}
\item[(a)]  $\sigma = \aleph_0$ and $\bold c:[\lambda]^2 \rightarrow
  \sigma$ and $\bar S$ witness $\Pr_*(\lambda,\lambda,\sigma,\theta)$
\sn
\item[(b)]  $\cP_1 = \{v_\alpha:\alpha < \lambda\} \subseteq
  [\lambda]^{\aleph_0}$ witness $\lambda^{\langle
    \theta;\aleph_0\rangle} = \lambda$
\sn
\item[(c)]  $\cP_2 = \{u_\alpha:\alpha < \lambda\}$ witness 
$\lambda^{(\lambda;\theta)} = \lambda$ when this is assumed and
witness $\lambda^{\langle \lambda;\theta\rangle} = \lambda$ when this
is assumed so $|u_\alpha| = \theta$.
\end{enumerate}
\end{enumerate}
\mn
1) We define $(G_i,I_i,J_i,\bold a_i,\bold b_i,\bold c_i)$ by
induction on $i < \theta$ such that:
\mn
\begin{enumerate}
\item[$(*)^2_i$]
\begin{enumerate}
\item[(a)]  the relevant parts of Definition \ref{c25} holds
\sn
\item[(b)]  for $i \in S_1$, clause $(c)(\varp)$ of Definition
  \ref{c26} holds and let $I^{\cg}_i = \{\alpha \in I_i:\gs_{i,\alpha}
  = \gs_{\cg}\}$
\sn
\item[(c)]  $\{I_{i,\alpha}:i < \theta,\alpha < \lambda\rangle$ is a
  partition of $\lambda$, each $I_{i,\alpha}$ is from $\{S_\beta:\beta
  < \lambda\}$ and $I_i = \bigcup\limits_{\alpha < \lambda}
  I_{i,\alpha}$
\sn
\item[(d)]  if $i \in S_2,\langle (\gs_{i,\alpha},a_{i,\alpha}):\alpha
  < \lambda\rangle$ is as in clause $(c)(\varp)$ of Definition
  \ref{c26}
\sn
\item[(e)]  assume $i \in S_2;j \in i \cap S_2$ and $\alpha <
  \lambda$, \then \, for some $\beta < \lambda$ we have
  $u''_{j,\alpha} = \{b_{i,s}:s \in J_{i,\beta}\}$ where
\sn
\begin{itemize}
\item  $u'_{j,\alpha} = \{\gamma < \lambda$: there is $a \in u_\alpha$
  such that $a \in (b^{-1}_{j,\gamma} G_j b_{j,\gamma}) \backslash
  \{e\}$ an $a$ has order 2 in $G_i\}$
\sn
\item  $u''_{j,\alpha} = \{a$: for some $\gamma \in u'_{j,\alpha},a
  \in (b^{-1}_{j,\gamma} G_j b_{j,\gamma}) \backslash \{e\}$ has order
  2 (in $Q_i$) is minimal (in $\lambda$) under those conditions
\end{itemize}
\sn
\item[(f)]  if $i \in S_3$, a parallel condition to clause (c) above
  using $\cP_2$ (instead $\cP_2$, i.e. $v_\alpha$ instead of
  $u_\alpha$) - FILL?
\end{enumerate}
\end{enumerate}
\mn
Clearly
\mn
\begin{enumerate}
\item[$(*)_3$]
\begin{enumerate}
\item[(a)]  we can carry the induction
\sn
\item[(b)]  the $\bold m$ defined naturally by the above belongs to
  $\bold M_1$.
\end{enumerate}
\end{enumerate}
\mn
[Why?  For $i=0$ see Definition \ref{c25}.  For $i=1,G_i$ is well
defined; for $i$ limit we let $G_i = \bigcup\limits_{j<i} G_j$ and for
$i=j+1$, by our choices in $(*)^2_j$ and Definition \ref{c25}, $G_i$
is well defined.  More elaborately, if $j \in S_1$, see \cite{Sh:312};
if $j \in S_2$, see \S(3A) and if $j \in S_3$, see xxx to choose
$\langle b_{j,s}:s \in I_{j,\alpha}\rangle$ again act as in
\ref{b15}(1), \ref{b16}(2) and we can put all of them together, see
\ref{b15}(2), \ref{b16}(2).]

So arriving to $i,G_i$ is well defined and the other objects are yet
to be chosen.
\medskip

\noindent
\underline{Case 1}:  $i \in S_1$

The number of relevant $(\gs,\bar a)$ is $\le \lambda$ if $i=1$ and
$=\lambda$ if $i \in (1,\theta)$, so there is no problem to choose
$\langle \gs_{i,s}:s \in J_i = I_i\rangle,\langle a_{i,s}:s \in I_i\rangle$.
\medskip

\noindent
\underline{Case 2}:  $i \in S_2$

Use the choice of $\cP_j$.
\medskip

\noindent
\underline{Case 3}:  $i \in S_3$

Use the choice of $\cP_2$.
\end{PROOF}

\begin{conclusion}
\label{i17}
(was n011)  Assume $\theta = \cf(\theta) \ge \aleph_1,\lambda^{\langle
  \theta;\aleph_0\rangle} = \lambda$ and $\lambda^{(\lambda;\theta)} =
\lambda,\gS \subseteq \Omega[\bold K_{\lf}]$ as in xx of cardinality
$\le \lambda$.

If $G \in \bold K^{\lf}_{\le \lambda}$ \then \, there is a
$(\lambda,\theta,\gS)$-full, complete full $H \in \bold
K^{\exlf}_\lambda$ extending $G$.
\end{conclusion}

\begin{PROOF}{\ref{i17}}
If $\lambda$ is regular, use \ref{i11}, \ref{i14}.  For $\lambda$
singular not strong limit, see \S(5B).
\end{PROOF}
\bigskip

\subsection {Singular $\lambda$} \label{5B}\
\bigskip

In \ref{i29}, \ref{i32} we define and use the alternative $\Pr_{**}$
for $\lambda$ singular not strongly limit.  $\Pr_{**}$ is slightly
stronger (and enables us to use $\lambda^{\langle
  \theta;\aleph_0\rangle} = \lambda$ instead of $\lambda^{\langle
  \theta;\aleph_0\rangle} = \lambda$, i.e. to deal with cardinals $<
2^{\aleph_0}$ (e.g. $\lambda < \aleph_\omega$ or even $<$ first fix
point in which case we can prove the pcf demand.  We also prove the ...?

\noindent
A relative of $\Pr_*$ is
\begin{definition}
\label{i29}
We define $\Pr_{**}(\lambda,\mu,\sigma,\bar\kappa,\theta)$ for
cardinals $\lambda$ as in \ref{i5}:
\mn
\begin{enumerate}
\item[(a)$'$]  $\bar S = \langle S_\alpha:\alpha < \lambda\rangle$ is
a partition of $\lambda$
\sn
\item[(d)$'$]  if $(\alpha)$ then $(\beta)$ when:
\sn
\begin{enumerate}
\item[$(\alpha)$]  $\alpha < \lambda,S=S_\alpha,u_{\varp,\iota}
  \subseteq [\lambda]^{< \kappa_\iota}$ for $\varp \in S,\iota < 2$
\sn
\item[$(\beta)$]  for some $A \subseteq \sigma$ of cardinality $<
  \kappa_0 + \kappa_1$, for every $j < \sigma$ there are $\varp <
  \zeta$ from $S$ such that:
\sn
\begin{itemize}
\item  $\bold c\{\varp,\zeta\} = j$
\sn
\item  if $\varp_1 \in u_\varp \backslash u_\zeta,\zeta_1 \in u_\zeta
  \backslash u_\varp$ and $\{\varp_1,\zeta_1\} \ne \{\varp,\zeta\}$,
  then $\bold c\{\varp_1,\zeta_1\}=0$
\end{itemize}
\end{enumerate}
\sn
\item[(e)]  $\bold c \rest [S_\alpha]^2$ is indecomposable.
\end{enumerate}
\end{definition}

\begin{claim}
\label{i32}
1) If $\lambda$ is singular, not strong limit and $\sigma < \lambda$,
\then \, $\Pr_{**}(\lambda,\lambda,\sigma,\aleph_0)$.

\noindent
2) If $\alpha < \lambda \Rightarrow |\alpha|^{< \kappa} <
\lambda,\lambda$ singular, not strongly limit, $\sigma <
\lambda,\theta = \cf(\theta) \le \sigma$, \then \,
$\Pr_{**}(\lambda,\lambda,\sigma,\bar\kappa,\theta)$. 
\end{claim}

\begin{PROOF}{\ref{i32}}
We can find by Chernikov-Shelah \cite{CeSh:1035} a $\mu < \lambda$ and
sub-tree $\cT$ of ${}^{\mu >}\mu$ (or ${}^{\mu >}2$) with $\mu$ nodes  
and $\ge \lambda$ branches (note that $\Ded(\mu) \ge \lambda$ means
$\sup\{|\lim(\cT):\cT$ a tree with $\le \mu$ nodes$\}$; but $\lambda$
is singular hence if $\mu \in [\cf(\lambda),\lambda)$, the supremum is
obtained. 

We may demand $\mu$ is ergular, though this is not essential.
Let $\langle \nu_\alpha:\alpha < \mu\rangle$ list $\cT$
 each appearing $\mu$ times.  Let $\langle \eta_\alpha:\alpha <
\lambda\rangle$ list different branches: we shall use, e.g. $\partial
= \mu^{++}$, guess clubs, i.e. let $\bar C = \langle C_\delta:\delta
\in S\rangle$ be such that $S \subseteq \{\delta
< \partial:\cf(\delta) = \cf(\mu)\},C_\delta$ a club of $\delta$ of order
type $\mu$ such that $\bar C$ guesses clubs, exists
(\cite[Ch.III,\S1]{Sh:g}).  Let $h:\mu \rightarrow
\sigma$ be such that $(\forall i < \sigma)(\forall \alpha <
\mu)(\exists^\mu \xi)(h(\xi) = i \wedge \nu_\xi = \nu_\alpha)$ for $i \,
\odd \, I_{i,\alpha} \cong \mu^{+7}$.

Recall $\id_a(\bar C)$ is the guessing club ideal, the normal version
on $\Dom(\bar C)$ (check notation).  Let $S_\alpha = [\partial
\alpha,\partial \alpha + \partial)$ and we define $\bold c:[\lambda]^2
\rightarrow \sigma$ as follows:
\mn
\begin{enumerate}
\item[$(*)$]
\begin{enumerate}
\item[(a)]  $\bold c(\{\partial \alpha + \varp,\partial \alpha +
  \zeta\}) = j$ \when \,
\sn
\begin{enumerate}
\item[$(\alpha)$]  $\alpha < \lambda$
\sn
\item[$(\beta)$]  $\varp < \zeta < \partial$
\sn
\item[$(\gamma)$]  $\zeta \in S$ hence $\cf(\zeta) = cf(\delta)$
\sn
\item[$(\delta)$]  $\nu_{\otp(C_\zeta \cap \varp)} \triangleleft
  \eta_\alpha$
\sn
\item[$(\varp)$]  $j = h(\otp(C_\zeta \cap \varp))$.
\end{enumerate}
\end{enumerate}
\end{enumerate}
\mn
Let us check Definition \ref{i5}; trivially clauses (a),(b),(c) hold.
\medskip 

\noindent
\underline{Clause (d)}:  So assume $\alpha < \lambda,S = S_\alpha$ is
 stationary in $\sup(S_\alpha) = \delta^*_\alpha = \partial \alpha
  + \partial$ or $\col(S_*) = \{\varp:\partial \alpha + \varp \in S\}
  \in \id^+_a(\bar C)$ and $u_\varp \in [\lambda]^{< \sigma}$ for
$\varp \in S$ and $j < \sigma$.

We have to find $\varp < \zeta$ from $S$ such that $c\{\varp,\zeta\} =
0$ and $[\varp_1 \in u_\varp \backslash u_\zeta \wedge \zeta_1 \in
u_\zeta \backslash u_\varp \wedge \{\varp_1,\zeta_1\} \ne
\{\varp,\zeta\} \Rightarrow \bold c\{\varp,\zeta\}=0$.

For each $\varp \in S$ let $v_\varp = \{\alpha < \lambda:[\partial
\alpha,\partial \alpha + \delta) \cap u_\varp \ne 0\}$, hence thre is
$S' \subseteq S$ such that $\col(S') \in \id_a + (\bar c)$ and
$\langle v_\varp:\varp \in S'\rangle$ is a $\Delta$-system with heart
$v_*$, for some $v_* \subseteq \lambda$.  For $\varp \in S'$ let
$u'_\varp = \{\partial \alpha + \xi:\xi < \gamma$ and for some
$\beta,\partial \beta + \xi \in u_\varp\}$, so $u'_\varp$ is a finite
subset of $[\partial \alpha,\partial \alpha + \partial)$.

Similarly there is $S'' \subseteq S'$ such that $\col(S'') \in
\id^+_a(\bar C)$ and $\langle u'_\varp:\varp \in S''\rangle$ is a
$\Delta$-system.  With hard $u_*$ such that $\varp \in S'' \Rightarrow
u'_\varp \cap \varp = u_*$ and $\zeta \in S'' \wedge \varp \in S''
\cap \zeta \Rightarrow u'_\varp \subseteq \zeta$.

Next choose $\nu \triangleleft \eta_\alpha$ (hence $\nu \in \cT$) such
that $\beta \in v_* \backslash \{\alpha\} \Rightarrow \neg(\nu
\triangleleft \eta_\beta)$, possible because:
\mn
\begin{itemize}
\item  if $\beta < \alpha$ then $\eta_\beta \in {}^{\mu>}\mu$ has
  length a limit ordinal, $\eta_\beta \notin \cT$ and $(\forall i <
  \ell g(\eta_\beta))(\eta_\beta \rest i \in \cT)$.
\end{itemize}
\mn
Next, choose $i < \mu$ such that $(\nu_i,j_i) = (\nu,j)$ and let $E :=
\{\delta < \partial:\delta$ a limit ordinal such that $S'' \cap
\delta$ is an unbounded subset of $\delta\}$, clearly $E$ is a club of
$\partial$ hence by the choice of $\bar c$.  So we can choose $\zeta
\in E \cap S''$ such that $C_\zeta \subseteq E$.
\medskip 

\noindent
\underline{Clause (e)}:  1) Simpler.

\noindent
2) Similarly only choose $\mu$ such that $\mu^{< \kappa} = \mu$.
\end{PROOF}

\begin{claim}
\label{c35}
1) Assume $\lambda$ is singular not strong limit, $\sigma =
\cf(\sigma) < \lambda,\theta = \cf(\theta) < \lambda$.  \Then \,
$\Pr_*(\lambda,\lambda,\sigma,\theta)$.

\noindent
2) Similarly with $\bar\kappa$ as in \ref{i32}.
\end{claim}

\begin{PROOF}{\ref{c35}}
Combine the proof of \ref{i8} and \ref{i32}.
\end{PROOF}

\begin{claim}
\label{i38}
Like \ref{i32} for regular $\lambda$.
\end{claim}

\begin{PROOF}{\ref{i38}}
Combine the proof of \ref{i8}, \ref{i32} (differently).
\end{PROOF}

\begin{question}
1) Does \ref{e16} fit this frame?

\noindent
2) Can we, as in \cite[4.9=La2,4.1=Ld36,pg.44,4.4=Ld38,pg.46]{Sh:312}
assume the object $\bold x$ consists of $K_0 \subseteq K_1 \subseteq
K_2$ are finite group, $K_1$ with trivial center, $\bar a_\ell$
generate $K_\ell$ for $\ell=0,1,2$ and $\bar a_0 \triangleleft \bar
a_1$.

\Then \, there is $\gs = \gs_{\cm}[\bold x]$ such that:
\mn
\begin{enumerate}
\item[(a)]  $\gs \in \Omega[\bold K_{\lf}]$
\sn
\item[(b)]  $k_{\gs} = \ell g(\bar a_1)$ and $n_{\gs} = \ell g(\bar
  a_2)$
\sn
\item[(c)]  $p_{\gs}(\bar x_{\gs}) = \tp_{\bs}(\bar
  a_1,\emptyset,K_1)$
\sn
\item[(d)]  if $(\alpha)$ then $(\beta)$ where:
\sn
\begin{enumerate}
\item[$(\alpha)$]
\sn
\begin{itemize}
\item  $G_1 \subseteq G_3$
\sn
\item  $\tp_{\bs}(\bar a'_1,\emptyset,G_1) = \tp_{\bs}(\bar
  a_1,\emptyset,K_1)$
\sn
\item  $\bar c$ realizes $q_s(\bar a'_1,G_1) \in G_3$
\sn
\item  we let $G_0 = \seb(\bar a'_1,G_3)$
\sn
\item  $G_2 = \seb(\bar a'_1 \char 94 \bar c,G_3)$
\sn
\item  $H_2 = \bold C_{G_2}(c \ell(\bar a_0),G_2)$
\sn
\item  $H_1 = \bold C_{G_0}(c \ell(\bar a_0))$
\sn
\item  $L = \bold C(G_0,G_1)$ and $H_1 = L \oplus H_0$
\end{itemize}
\sn
\item[$(\beta)$]  $\NF^3(\langle G_\ell:\ell \le 3 \rangle,\bold
  C(\bar a_b,G_1),\bold C(\bar a_0)$.
\end{enumerate}
\end{enumerate}
\mn
3) Maybe we can find a $\gs$ acting on a pair of groups $(G'_1
\subseteq G_1)$?
\end{question}
\bigskip

\subsection {Strong Limit Singular and \cite{Sh:331}} \label{5C} \
\bigskip

Here we generalize \cite{Sh:331} but ignore the
``$\theta$-indecomposable" in the definition of $\Pr_*/\Pr_{**}$.  The
results are stronger than in \S(5A) and \S(5B) coverall all $\lambda >
\aleph_0$, but the main proof tells us to change the proof in
\cite{Sh:331} in some points; not such a good choice, but also to
repeat it, is not a great one either.
\bigskip

\centerline {$* \qquad * \qquad *$}
\bigskip

It seems we remain with the case $\lambda$ strong limit.  As earlier,
it seems tempting to use \cite{Sh:331}.  A problem is there
considering $\eta \in J_\alpha,M$ or $\langle M_n,N_n:n <
\kappa\rangle$, we are consider interference from
$\Sigma\{J_\beta:\beta \ne \alpha\}$, whereas now we consider
interference from $\Sigma\{J_\beta:\beta\} \backslash \{\eta\}$ a
point neglected in \cite{Sh:331}.  If for transparency the answer to
the question above to ``only $\ne j$" is yes, then we may
\mn
\begin{enumerate}
\item[$\boxplus$]  in Definition \cite[1.5=L7.3]{Sh:331}; here $\mu =
  \aleph_0 = \kappa$ suffice
\newline
\underline{Clause (D)}:  i.e. $\ell=4$: in clause (vii) replacing
``$\nu \in P^J_\omega$" by $\nu \in P^J_\omega \cup (P^I_\omega
\backslash \{\eta\})$.
\end{enumerate}
\mn
Similarly for $(D^-),(D^+),(E)$, etc.  Does \cite[1.1=L7.1]{Sh:331}
essentially suffice?  Yes, use $\;langle \ldots,\eta(n_,\eta(n)+1
\ldots:n < \omega\rangle$!  So let us look at the various claims and
lemmas in \cite{Sh:312} which deals with some case, i.e. ``if
$\lambda$ satisfies ...".

We naturaly let $\langle (\eta_i,i_\varp):\varp < \lambda\rangle$ list
$(\eta,i)$ such that $i < \lambda,\eta \in P^{I_i}_\omega,\bold
c\{\varp,\zeta)$ is 0 if $i_\varp \ne i_\varp$ and if $\Pr_1(\ell
g(\eta_\varp \cap \eta_\zeta)),\pr$-pairing.
  
\begin{notation}
\label{i41}
1) $K^\omega_{\tr}$ is up to isomorphism the class of subtrees of
$({}^{\omega \ge}\alpha,\triangleleft)$ for some $\alpha$, expanded by
$<_{\lex}$ and $\cP_n = I \cap {}^n \alpha$.

\noindent
2) Let $K_\lambda = \{I:{}^{\omega >}\lambda \subseteq I \subseteq
{}^{\omega \ge}\lambda\},I$ has cardinality $\lambda$, expanded as above.
\end{notation}

\begin{definition}
\label{i44}
We say $I \in K^\omega_{\tr}$ is $(\mu,\kappa)$-duper-unembeddable
into $J \in K^\omega_{\tr}$ \If \, : for every regular large enough  
$\chi^*$ (for which $\{I,J,\mu,\kappa\} \subseteq \cH(\chi^*)$), for
simplicity a well ordering $<^*_{\chi^*}$ of $\cH(\chi^*)$ and
$x \in {\cH}(\chi^*)$ we have (compared to \cite[1.1=L7.1]{Sh:331} we
strengthen clause (vi) and add $M_n = N_n$ retaining both to case quotation):
\mn
\begin{enumerate}
\item[$(*)$]  there are $\eta, M_{n},N_{n}$ (for $n<\omega$) such that:
\sn
\begin{enumerate}
\item[$(i)$]  $M_{n} = N_{n} \prec M_{n+1} \prec ({\cH}(\chi^*),
\in,<^*_{\chi^*})$ 
\sn
\item[$(ii)$]  $M_{n} \cap \mu  = M_{0} \cap \mu$ and $\kappa
  \subseteq M_0$
\sn
\item[$(iii)$]  $I,J,\mu,\kappa$ and $x$ belong to $M_0$
\sn
\item[$(iv)$]  $\eta \in P^I_{\omega}$
\sn
\item[$(v)$]   for each $n$  for some $k,\eta \rest k \in M_n,
\eta \rest (k+1) \in N_n \backslash M_n$ 
\sn
\item[$(vi)$]  for each $\nu \in P^J_{\omega}$ or $\nu \in P^I_\omega
  \backslash \{\eta\}$, the sequence $|\{k<n:\{\nu \rest \ell:\ell <
  \omega\} \cap N_{\omega +1} \nsubseteq N_k\}|/n$ converge to zero.
\end{enumerate}
\end{enumerate}
\end{definition}

\begin{notation}
\label{i47}
We may write $\mu$ instead $(\mu,\mu)$ and may omit it if $\mu = \aleph_0$.
\end{notation}

\begin{remark} 
\label{i50}
1) The $x$  can be omitted (and we get equivalent definition
using a bigger $\chi^*$) but in using the definition, with $x$ it is 
more natural: we construct something from a sequence of $I$'s, 
we would like to show that there are no objects such that ... and 
$x$ will be such undesirable object in a proof by contradiction.

\noindent
2) We can also omit $<^*_\kappa$ at the price of increasing $\chi^*$.

\noindent
3) Here we use $\mu = \kappa = \aleph_0$.

\noindent
4) We can deal also with the more general properties from
\cite[1.5=L7.3]{Sh:331}. 
\end{remark}

\begin{definition}
\label{i53}
1) $K^\omega_{\tr}$ has the $(\chi,\lambda,\mu,
\kappa)$-duper-bigness property \If \,: there are $I_{\alpha} \in
(K^\omega_{\tr})_{\lambda}$ for $\alpha < \chi$ such that for 
$\alpha \ne \beta,I_{\alpha}$ is $(\mu,\kappa)$-duper unembeddable
into $I_{\beta}$.

\noindent
2) $K^\omega_{\tr}$ has the full
$(\chi,\lambda,\mu,\kappa)$-duper-bigness property \If \,: \\
there are $I_{\alpha} \in (K^\omega_{\tr})_{\lambda}$ for $\alpha <
 \chi$ such that:
\mn
\begin{enumerate}
\item[(a)]  for every $\alpha < \chi,I_{\alpha}$ is 
$(\mu,\kappa)$-duper unembeddable into 
$\sum\limits_{\beta < \chi,\beta \ne \alpha} I_{\beta}$
\sn
\item[(b)]  for eery odd $\alpha < \chi$.
\end{enumerate}
\mn
3) We may omit $\kappa$ if $\kappa = \aleph_0$.
\end{definition}

\begin{fact}
\label{x56}
1) If $I \in K^\omega_{\tr}$ is $(\mu,\kappa)$-duper$^m$-unembeddable
into $J \in K^\omega_{\tr}$
\then \,  $I$ is $(\mu,\kappa)$-duper$^\ell$-unembeddable into $J$ 
when $1 \le \ell \le m \le 7,(\ell,m) \ne (5,6),\ell,m \in \{1,2,3,4,5,6,7\}$
and when $(\ell,m) \in \{(3,4^-),(4^-,4),(4,4^+),(4^+,6),(6,6^+),
(4^+,7^-)$,\\
$(7^-,7),(7,7^+),(7^-,7^\pm),(7^\pm,7^+),(6^+,7^+),(6,7),
(6^-,6^\pm),(6^\pm,6^+)\}$.

\noindent
2) If $K^\omega_{\tr}$ has the 
$(\chi,\lambda,\mu,\kappa)$-duper$^m$-bigness property \then \, 
$K^\omega_{\tr}$ has the $(\chi,\lambda,\mu,\kappa)$-duper$^\ell$-bigness
property for  $(\ell,m)$  as above.

\noindent
3) If $K^\omega_{\tr}$ has the full
$(\chi,\lambda,\mu,\kappa)$-duper$^m$-bigness property \then \,
$K^\omega_{\tr}$ has the full
$(\chi,\lambda,\mu,\kappa)$-duper$^\ell $-bigness property for $(\ell,m)$
as above.

\noindent
4) All those properties has obvious monotonicity properties: we
can decrease $\mu,\kappa$ and $\chi$ and increase $\lambda$ (if we add
to $I$ a well ordered set in level 1, nothing happens).
\end{fact}

\begin{theorem}
\label{i57}
If $\lambda > \mu$ then $K^\omega_{\tr}$ has the full
$(\lambda,\lambda,\mu)$-duper-bigness property, hence the
$(2^\lambda,\lambda,mu)$-duper-bigness property.
\end{theorem}

\begin{PROOF}{\ref{i57}}
We have to go over the proof in \cite[\S2,\S3]{Sh:331} and say how to
adapt this.
\medskip

\noindent
\underline{Case A}:  $\lambda$ regular $> \mu$,
\cite[1.11=L7.6(1)]{Sh:312}, more \cite[2.13=L7.8I]{Sh:312}.

Obviously fine for $\Pr_*$, for $\Pr_{**}$ we need $\lambda >
\kappa^+$ to guess clubs of $S^\lambda_{\aleph_0}$, i.e. the second
citation seems to work for $K_{\tr(\theta)},\lambda > \theta^+$.
\medskip

\noindent
\underline{Case B}:  $\lambda$ singular, $\chi = \chi^\kappa < \lambda
\le 2^\chi$, \cite[1.11=L76(2)]{Sh:331}.

The proof as stated fails.  However, if $\chi = \chi^\mu$ we can let
$\gB$ be a model with universe $\chi^+$, countable vocabulary such
that $c \ell_{\gB}(\{\alpha\} \cup \mu) \ge \alpha$.  Let $\langle
u_\varp:\varp < \chi\rangle$ list $[\chi]^{\le \mu}$, even length
$\chi^+$ is O.K.  Now we partition $S^\zeta$ to $\langle
S^{\zeta,\varp}:\varp < \chi\rangle$, each a stationary subset of
$\chi^+$.

For each $(\zeta,\varp)$ let $\bar C^{\zeta,\varp} = \langle
C_\delta:\delta \in S^{\zeta,\varp}\rangle$ guess clubs.  Now if
$\delta + ^{\zeta,\varp}$, let $\cM_\delta = \{M:M \subseteq \gB,M
\cap \mu = u_\varp,\alpha \in C^{\zeta,\varp} \Rightarrow M \rest
\alpha \subseteq \gB$ and $\alpha \in C^{\zeta,\varp} \wedge \beta =
\suc(\alpha,C^{\zeta,\varp}) \Rightarrow M \cap (\alpha,\beta) \ne
\emptyset\}$.  For each $M \in \cM_{\zeta,\varp}$ let $\eta_M \in
{}^\omega \delta$ be, e.g. $\eta(n) = \min\{\beta \in M:|\beta \cap M
\cap C_\delta| = n\}$.

Lastly, $I_\alpha = ({}^{\omega >}\lambda) \vee \{\eta_\mu$: there are
$M,\zeta,\varp,\delta$ such that $\zeta < \chi,\varp < \chi,\delta \in
S^{\zeta,\varp},M \in M_\delta$ and $\zeta \in A_\alpha\}$.
\medskip

\noindent
\underline{Case C}:  $\lambda$ strong limit singular, $\cf(\lambda) >
\aleph_0$, \cite[1.11=L7.6(3)]{Sh:331}.

Considering the use of \cite[1.13=L7.6A]{Sh:331},
\cite[1.13=L7.6B]{Sh:331} this should be clear, similar changes to
case 2 (more details needed)?
\medskip

\noindent
\underline{Case D}: $\lambda$ singular and equal to
$\sup\{\theta:\theta$ singular and $\pp(\theta) >\theta^+\}$,
\cite[1.16=L7.7]{Sh:331}, add there $S_i = \{\delta:\delta \in
\theta_i,\theta_{i+1})$ and $\cf(\delta) =\aleph_0\}$.

The proof is like the regular case, the treatment of $\nu \in
P^{I_i}_\omega$ rather than $\nu \in P^{\Sigma\{I_j:j \ne
  i\}}_\omega$, as in case 1.
\medskip

\noindent
\underline{Case E}: $\chi^{\aleph_0} \ge \lambda > \mu^{+2} +
\chi^{+2}$, see \cite[2.1=L7.8]{Sh:331}, proof after
\cite[2.13=L7.8I]{Sh:331} which improves the work if:
\mn
\begin{enumerate}
\item[(a)]  we double the $\eta$'s, i.e. $\eta' = \dbl(\eta)$ where
$\eta'(n^2 + \ell) = \eta(n)$ for $\ell \in [0,(n+1)^2 = n^2]$? doubtful
\sn
\item[(b)]  we use the property of \cite[1.1=L7.1]{Sh:331}.
\end{enumerate}
\medskip

\noindent
\underline{Case F}:  $\lambda$ singular, there $\ga_\varp(\varp <
\cf(\lambda),\otp(\ga_\varp) = \omega,(\Pi
\ga_\varp,<_{J^{\bd}_\omega})$ has true cofinality, $\theta =
\sup(\ga_\varp)$, pairwise almost disjoint, see \cite[2.15=L7.9]{Sh:331}.

remark: seems to cover $\lambda = \beth_\omega(\mu_1),\mu_1 > \mu$.
The non-coverd case seemed ...?
\medskip

\noindent
\underline{Case G}:  $\lambda = \beth_\omega(\mu)$, see
\cite[2.17=L7.10]{Sh:331}.
\end{PROOF}

\begin{claim}
\label{i59}
If $K^\omega_{\tr}$ has the full
$(\lambda,\lambda,\aleph_0,\aleph_0)$-duper bigness property, \then \,
$\Pr_{**}(\lambda,\aleph_0,\aleph_0)$. 
\end{claim}

\begin{PROOF}{\ref{i59}}
Let $\langle I_\alpha:\alpha < \lambda\rangle$ witness the property,
\wilog \, ${}^{\omega >}\lambda \subseteq I_i \subseteq {}^{\omega
  \ge}\lambda$.  Let $\{(\eta_\varp,\alpha_\varp):\varp < \lambda\}$
list $\{(\eta,\alpha):\alpha < \lambda$ and $\eta \in
P^{I_\alpha}_\omega\}$, i.e. $\eta \in I_Ii \cap {}^\omega \lambda$
and let:
\mn
\begin{itemize}
\item  $S_\ell = \{\varp < \lambda:\alpha_\varp = \alpha\}$
\sn
\item  $\bold c\{\varp,\zeta\}$ is: $\pr_1(\ell g(\eta_\varp \cap
  \eta_\zeta)$ if $(\exists \alpha)(\zeta,\varp \in I_\alpha)$

\hskip40pt $0$ otherwise
\end{itemize}
\end{PROOF}
\bigskip

\centerline {$* \qquad * \qquad *$}
\bigskip

However, we still like to cover the indecomposable.

\begin{definition}
\label{i61}
1) $\bold C_{\lambda,\mu}$ is the set of $\bold c:[\lambda]^2
\rightarrow \sigma$.

\noindent
2) We say $\bold c$ or $(\bold c,\bar S)$ exemplifies
$\Pr(\lambda,\mu,\theta,\kappa,\sigma)$ \when \, $(\lambda \ge \mu \ge
\theta,\sigma,\theta$ regular):
\mn
\begin{enumerate}
\item[(a)]  $\bar S = \langle S_\alpha:alpha < \lambda\rangle$ are
  pairwise disjoint subsets of $\lambda$
\sn
\item[(b)]  $\bold d:[\lambda]^2 \rightarrow \mu$ but $\varp \in
  S_\alpha \wedge \zeta \in S_\beta \wedge \varp \pm \zeta \Rightarrow
  \bold c\{\varp,\zeta\} = 0$
\sn
\item[(c)]  if $\chi$ large eough, $\lambda \in H(\chi),\alpha <
  \lambda$ and $\iota = 1 \wedge j_\iota < \sigma$ or $\iota = 2
  \wedge j_\iota < \theta$, \then \, there is a pair $(\bar
  N,\zeta_*)$ such that:
\sn
\begin{enumerate}
\item[(i)]  $\bar M = \langle N_n:n < \omega\rangle,N_n \prec \gB =
(\cH(\chi^+),\theta,<^*_*)$ is increasing with $n$
\sn
\item[(ii)]  $M_n \cap \mu = M_0 \cap \mu$ and $\kappa +1 \subseteq
M_0$
\sn
\item[(iii)]  $\bar S,\bold c,\alpha,j_\iota \in M_0$
\sn
\item[(iv)] $\zeta_* \in S_\alpha$
\sn
\item[(v)]  $(\alpha) + (\beta) \Rightarrow (\gamma)$ where
\sn
\begin{enumerate}
\item[$(\alpha)$]   $\bar zeta
 = \langle \zeta_0,\dotsc,\zeta_n\rangle,\xi_n = \zeta_*,n < \kappa$
 and $\zeta_0,\dotsc,\zeta_{n-1} \in \lambda \backslash \{\zeta_*\}$
\sn
\item[$(\beta)$]  $\varphi(\bar x,\bar y)$ a formula in
  $\bbL(\tau_{\gB}),\bar a \in {}^{\ell g(\bar y)}M$ and $\gB \models
  \varphi [\bar\zeta,\bar a]$.
\end{enumerate}
\end{enumerate}
\end{enumerate}
\mn
3) Versions:
\mn
\begin{enumerate}
\item[$(\gamma)$]   there is $\bar zeta  = \langle \zeta'_\ell:\ell
  \le n\rangle$ such that:
\sn
\begin{enumerate}
\item[$\bullet_1$]  $\models \varphi(\bar\zeta',\bar a]$ and
  $\bar\zeta' \in {}^{n+1} M$
\sn
\item[$\bullet_2$]  if $\ell < n$ and $\zeta'_\ell = \zeta_\ell$
\sn
\item[$\bullet_3$]  if $\ell < n$ and $\zeta_\ell \notin M$ then
  $\bold c\{\zeta'_\ell,\zeta_\ell\} = 0$
\sn
\item[$\bullet^+_3$]  moreover, if $\ell_1,\ell_2 \le
  n,\zeta'_{\ell_1} \ne \zeta_{\ell_1},\zeta'_{\ell_2} \ne
  \zeta_{\ell_2}$ and $(\ell_1,\ell_2) \ne
  (\zeta_{\ell_1},\zeta_{\ell_2}) \ne (\zeta_n,\zeta'_n)$ then $\bold
  c\{\zeta_{\ell_1},\zeta'_{\ell_2}\} = 0$
\sn
\item[$\bullet_4$]  if $\iota = 1$ then $\bold c\{\zeta'_n,\zeta_*\} =
  j_1$
\sn
\item[$\bullet_5$]  if $\iota = 2$ then $\bold c\{\zeta'_n,\zeta_*\}
  \in [j_2,\theta)$.
\end{enumerate}
\end{enumerate}
\end{definition}

\begin{remark}
\label{i64}
This is parallel to Definition \ref{i57}, we could have defined a
parallel of ... ?
\end{remark}

\begin{claim}
\label{i67}
If $\Pr_{\bullet \bullet}(\lambda,\lambda,\sigma,\aleph_0,\theta)$,
\then \, $\Pr_{**}(\lambda,\lambda,\sigma,\kappa,\theta)$ (check).
\end{claim}

\begin{PROOF}{\ref{i67}}
Should be clear.
\end{PROOF}

\begin{theorem}
\label{i70}
Assume $\lambda > \mu \ge \kappa + \sigma + \theta$ and $\theta$
regular.

Then $\Pr_{\bullet \bullet}(\lambda,\mu,\sigma,\kappa,\theta)$, at
least for $\kappa = \aleph_0$ which suffices.
\end{theorem}

\begin{PROOF}{\ref{i70}}
Like the proof of \ref{i56}, i.e. \cite{Sh:331} - FILL.
\end{PROOF}
\bigskip

\subsection {Continuation} \label{5D} \
\bigskip

\begin{claim}
\label{n14}
(?) Let $G_1 \in \bold K_{\lf}$ and $a_\alpha \in G_1$ of order 2 for,
$\alpha < \lambda$.  We can find $G_2,\bar b$ such that:
\mn
\begin{enumerate}
\item[(a)]  $G_1 \subseteq G_2 \in \bold K_{\lf}$ and $G_2$ has cardinality
  $|G_1| + \lambda$;
\sn
\item[(b)]  $\bar b = \langle b_{\alpha,\iota}:\alpha < \lambda,\iota
  < 2\rangle$ is a sequence of elements of $G_2$;
\sn
\item[(c)]  $b_{\alpha,\iota} \in G_2$ and $G_2 = \seb(G_1 \cup
  \{b_{\alpha,\iota}:\alpha < \lambda,\iota < 2\},G_2)$;
\sn
\item[(d)]  $\langle b_{\alpha,\iota}:\alpha < \lambda\rangle$ is a
  sequence of pairwise commuting members of $G_2$, each of order 2;
\sn
\item[(e)]  $a_\alpha \in \seb(\{b_{\alpha,0},b_{\alpha,1}\},G_2)$;
\sn
\item[(f)]  $G_1 \le_{\gS} G_2$;
\sn
\item[(g)]   $\tp_{\qf}(b_\alpha,\langle b_{\alpha,\iota}:\alpha
  \in u,\iota <2 \rangle,G_1,G_2) = q_{\gs}(\langle a_\alpha:\alpha
  \in u\rangle,G_1)$ with $\gs$ computed from $\tp_{\qf}(\langle
  a_\alpha:\alpha \in u\rangle,G_1)$.
\end{enumerate}
\end{claim}

\begin{PROOF}{\ref{n14}}
Use \cite[\S4]{Sh:312}.
\end{PROOF}

\begin{claim}
\label{n23}
(?) Assume $G \in \bold K^{\lf}_{\le \lambda},\aleph_1 \le \theta =
\cf(\theta) \le \lambda$ and $\lambda = \lambda^{\langle
  \theta;\aleph_0\rangle}$ and
$\Pr_*(\lambda,\lambda,\aleph_0,\aleph_0)$.  \Then \, there is a
complete $H \in \bold K^{\lf}_\lambda$ extending $\nu$.
\end{claim}

\begin{PROOF}{\ref{n23}}
See \ref{a15}, or see the above, or see \ref{e34}, a copy of \ref{c35}.
\end{PROOF}
\bigskip

Our aim is to omit in Conclusion \ref{c41} the assumption ``$\lambda$
is a successor of a regular", though we retain ``$\lambda \ge
\beth_\omega$" or at least some consequences of this.  So the problem
is having a weaker colouring theorem (than
$\Pr_0(\lambda,\lambda,\lambda,\aleph_0))$ in the proof of \ref{c35},
while still having $\theta$-indecomposability (recalling $\Pr^*_0$
is from Definition \ref{w11}(4A)).

\begin{claim}
\label{n72}
1) $\Pr'_0(\lambda,\aleph_0,\aleph_0)$ \when \,:  (check $\Pr$ or $\Pr'$?)
\mn
\begin{enumerate}
\item[(a)]  $\lambda = \cf(\lambda) > \aleph_1$.
\end{enumerate}
\end{claim}

\begin{PROOF}{\ref{n72}}
Let $S \subseteq S^\lambda_{\aleph_0}$ be stationary.
Let $\bar c = \langle c_\delta:\delta: \in S\rangle$ guess clubs, that is:
\mn
\begin{enumerate}
\item[$(*)$]  $C_\delta \subseteq \delta =
  \sup(C_\delta),\otp(C_\delta) = \omega$ and for every club $E$ of
  $\lambda$ for stationarily many $\delta \in S \cap E$ we have
  $C_\delta \subseteq S$ (question: adding a regressive function is
  constant in $C_\delta \cup \{\delta\}$?).
\end{enumerate}
\mn
Given $\langle u_\varp:\varp < \lambda\rangle$ as in $\Pr_*$ for some
club $E$ of $\lambda,\varp < \delta \in E \Rightarrow u_\varp
\subseteq \delta$. 
\end{PROOF}

Maybe more suitable than $\Pr_{2.5}$ is another colouring specifically
needed here.
\begin{definition}
\label{e3}
1) Let $\Pr_{2.6.1}(\chi,\lambda,\mu,\theta,\partial,\sigma)$ means 
that some $(\bold c,\bar W)$ witness it which means that (if $\chi =
\lambda$ we may omit $\chi$ and if $\chi = \lambda = \mu$ we may omit
$\chi$ and $\lambda$):
\mn
\begin{enumerate}
\item[$(*)$]
\begin{enumerate}
\item[(a)]  $\bar W = \langle W_i:i < \chi\rangle$ is a partition
  $\lambda$ to $\chi$ sets each of cardinality $\lambda$;
\sn
\item[(b)]  $\bold c:[\lambda]^2 \rightarrow \mu$;
\sn
\item[(c)]  if $i$ is even, then $\bold c \rest [W_i]^2$ is
  $\theta$-indecomposable, that is, if $\langle \cU_\varp:i <
  \theta\rangle$ is $\subseteq$-increasing with union $W_i$ \then \,
  for some $\varp < \theta$ we have $\mu = \Rang(\bold c \rest
  [\cU_\varp]^2)$;
\sn
\item[(d)]  we have $\boxplus_1 \Rightarrow \boxplus_2$ where for
  $\ell=1,2$;
\begin{enumerate}
\item[$\boxplus_\ell$]
\begin{enumerate}
\item[$(\alpha)$]   $i$ is odd;
\sn
\item[$(\beta)$]  $\cU_\ell \in [W_i]^\lambda$ stationary;
\sn
\item[$(\gamma)$]  $\varp \in u^\ell_\varp \in [\chi]^{< \sigma}$
for $\varp \in \cU_\ell$;
\sn
\item[$(\delta)$]  $\langle u^\ell_\varp:\varp \in \cU_\ell\rangle$ 
is a $\Delta$-system with heart $u_\ell$;
\sn
\item[$(\varp)$]  $\cU_\ell \cap u_\ell = \emptyset,\cU_1 \cap
\cU_2 = \emptyset$;
\sn
\item[$(\zeta)$]   $\gamma < \partial$.
\end{enumerate}
\end{enumerate}
\end{enumerate}
\end{enumerate}
\mn
\Then \, for some $\varp_1,\varp_2$ we have for $\ell=1,2$:
\mn
\begin{enumerate}
\item[$(\alpha)$]  $\varp_\ell \in \cU_\ell$;
\sn
\item[$(\beta)$]  $\bold c\{\varp_1,\varp_2\} \in [\gamma,\partial)$;
\sn
\item[$(\gamma)$]  $u^1_{\varp_1} \backslash u_1,u^2_{\varp_2}
  \backslash u_2$ are disjoint;
\sn
\item[$(\delta)$]  if $\zeta_1 \in u^1_{\varp_1} \backslash u_1$ and
 $\zeta_2 \in u^2_{\varp_2} \backslash u_{\zeta_2}$ and $\{\zeta_1,\zeta_1\}
\ne \{\varp_1,\varp_2\}$ \then \, $\bold c\{\varp_1,\zeta_1\} = 0$.
\end{enumerate}
\mn
2) We define $\Pr_{2.6.2}(\chi,\lambda,\mu,\theta,\partial,\sigma)$
similarly but replacing clause (d) by:
\mn
\begin{enumerate}
\item[$(d)'$] we have $\boxplus_1 \Rightarrow \boxplus_2$ when:
\begin{enumerate}
\item[$\boxplus_1$]
\begin{enumerate}
\item[$(\alpha)$]  $i$ is odd;
\sn
\item[$(\beta)$]  $\cU_1 \in [W_i]^\lambda$;
\sn
\item[$(\gamma)$]  $\varp \in u^1_\varp \in [\lambda]^{< \sigma}$ for
  $\varp \in \cU_1$;
\sn
\item[$(\delta)$]  $\langle u^1_\varp:\varp \in \cU_1\rangle$ is a
  $\Delta$-system with heart $u_1$;
\sn
\item[$(\varp)$]  $\varp \in u^2_\varp \in [\lambda]^{< \sigma}$ for
  $\varp \in W_i$;
\sn
\item[$(\zeta)$]  $\gamma < \partial$;
\end{enumerate}
\sn
\item[$\boxplus_2$]  for some $\varp_1,\varp_2$ we have:
\begin{enumerate}
\item[$(\alpha)$]  $\varp_1 \in \cU_1,\varp_2 \in W_i$;
\sn
\item[$(\beta)$]  $\bold c\{\varp_1,\varp_2\} \ge \gamma$;
\sn
\item[$(\gamma)$]  $u^1_{\varp_1} \backslash u_1,u^2_{\varp_2}$ are
  disjoint;
\sn
\item[$(\delta)$]  if $\zeta_1 \in u^1_{\varp_1} \backslash
u_1,\zeta_2 \in u^2_{\varp_2}$ and $\{\zeta_1,\zeta_2\} \ne
\{\varp_1,\varp_2\}$ \then \, $\bold c\{\zeta_1,\zeta_2\} = 0$.
\end{enumerate}
\end{enumerate}
\end{enumerate}
\end{definition}

\noindent
We need the following in proving $\Pr_{2.7}(\lambda,\aleph_0,\aleph_0)$
for $\lambda > \aleph_1$.
\begin{claim}
\label{e7}
If $\lambda > \aleph_1$ is regular, \then \,
$\Pr_{2.7}(\lambda,\aleph_0,\aleph_0)$. 
\end{claim}

\begin{PROOF}{\ref{e7}}
Let $\sigma = \aleph_0,\theta = \aleph_0$ and we prove more.
\medskip

\noindent
\underline{Case 1}:  $\lambda$ successor of regular.

Follow by \ref{w8}(1).

\noindent
\underline{Case 2}:  $\lambda = \cf(\lambda) > \partial^+,\partial =
\cf(\partial),\lambda > \theta,\cf(\theta)$.
\sn
\begin{enumerate}
\item[$(*)_1$]  Let $\langle W_{i,\varp}:i < \lambda,
\varp < \lambda\rangle$ be a sequence of pairwise disjoint stationary
subsets of $\lambda$ such that:
\sn
\begin{enumerate}
\item[(a)]  if $\delta \in W_{i,\varp},i$ is even then $\cf(\delta) =
  \theta$
\sn
\item[(b)]  if $\delta \in W_{i,\varp},i$ is odd then $\cf(\delta) = \partial$.
\end{enumerate}
\end{enumerate}
\mn
Let $W_i = \bigcup\limits_{\varp < \lambda} W_{i,\varp}$ and consider:
\sn
\begin{enumerate}
\item[$(*)_2$]  for each $i < \lambda,\varp < \lambda$ 
let $W'_i = \{\delta \in W_{i,\varp}:\delta = \sup(W_i
\cap \delta)\}$ and for each $\delta \in W'_{i,\varp}$ choose $C_\delta
\subseteq \delta \cap W_i$ of order type $\cf(\delta)$, unbounded in $\delta$.
\end{enumerate}
\mn
We define $\bold c$ as follows.  FILL.
\end{PROOF}

\begin{definition}
\label{e10}
1) Let $\bold M_{2.7} = \bold M^{2.7}_{\lambda,\partial,\bar S}$ be
the class of $\bold m \in \bold M^2_{\lambda,\theta,\bar S}$ such that
$\bold c$ is as in $\Pr_{2.7}(\lambda,\theta,\partial)$ where
$\partial = \theta$, see below.

\noindent
3) We define $\Pr_{2.7}(\chi,\lambda,\theta,\partial,\sigma)$ as in
\ref{e3}(1) but replace clause (d) by:
\mn
\begin{enumerate}
\item[$(d)''$]  for some $\mu,\mu = \cf(\mu) < \theta$ if $c
  \ell:[\lambda]^{< \aleph_0} \rightarrow [\lambda]^{< \mu}$ we have:
\begin{itemize}
\item  for each $i  \, \odd,\alpha < \lambda,W_{i,\alpha} \cap \{\delta <
  \lambda:\cf(\delta) = \mu\}$ is stationary;
\sn
\item  if $\cU_1 \subseteq W_{i,\alpha}$ is stationary $\langle
  u^1_\varp:\varp \in \cU_1\rangle,\varp \in u^1_\varp \in
  [\lambda]^{< \mu}$ and $\varp \in u^2_\varp \in [\lambda]^{< \mu}$
  and $\xi < \sigma$ (maybe $\sigma = \mu$) we can find $\varp <
  \zeta,u_1,u_2$ such that:
\end{itemize}
\sn
\begin{enumerate}
\item[$(\alpha)$]  $\zeta \in \cU_1,\varp < \zeta,\varp \in W_{i,\alpha}$;
\sn
\item[$(\beta)$]  $u^1_\varp \subseteq u_1$ and $u^2_\varp \subseteq u_2$;
\sn
\item[$(\gamma)$]  $u_1,u_2$ are $c \ell$-closed;
\sn
\item[$(\delta)$]   $\varp \in u_1 \backslash u_2,
\zeta \in u_2 \backslash u_1$;
\sn
\item[$(\varp)$]  if $\varp_1 \in u_1 \backslash u_2,\zeta_1 \in u_2
  \backslash u_1$ and $\{\varp_1,\zeta_1\} \ne \{\varp,\zeta\}$ then
  $\bold c\{\varp_1,\zeta_1\} = 0$;
\sn
\item[$(\zeta)$]  $\bold c\{\varp,\zeta\} = \xi$.
\end{enumerate}
\end{enumerate}
\end{definition}

\begin{claim}
\label{e13}
Assume $\lambda > \beth_\omega$ is a successor cardinal and $\theta =
\cf(\theta) < \beth_\omega$ large enough (or just $\lambda =
\lambda^{\langle \theta;\theta\rangle}$). 

\noindent
1) If $G \in \bold K^{\lf}_{\le \lambda}$ and $\gS$ is 
as in \ref{c23}(3), \then \, there is a
complete $(\lambda,\theta,\gS)$-full $H \in \bold K^{\exlf}_\lambda$
extending $G$.

\noindent
2) $\Pr_{2.7}(\lambda,\theta,\partial)$ holds for some
$\theta,\partial \in (\aleph_0,\lambda)$.
\end{claim}

\begin{PROOF}{\ref{e13}}
By \ref{w14} and \ref{c41}, \wilog \, $\lambda = \mu^+,\mu$ singular.
Let $\theta = \partial = \cf(\theta) \in (\aleph_0,\lambda)$.  We shall
choose $\bold c$ differently such that:
\mn
\begin{enumerate}
\item[$(*)_1$]  letting $I^\ell = \cup\{I_i:i \in S_\ell\}$ for
  $\ell=0,1,2,3$ we have:
\sn
\begin{enumerate}
\item[(a)]  $\bold c \rest I^3$ witness
  $\Pr_0(\lambda,\lambda,\aleph_0,\aleph_0)$ if possible or at least
$\Pr_{2.7}(\lambda,\theta,\partial,\aleph_0)$, see \ref{e10}, for some
  $\theta,\partial$;
\sn
\item[(b)]  $\bold c \rest I^2$ satisfies:
\begin{enumerate}
\item[$(\alpha)$]  for some club $E$ of $\lambda$ such that if $\varp
  < \zeta$ are from $I^2$ and $(\varp,\zeta] \cap E \ne \emptyset$
  \then \, $\bold c\{\varp,\zeta\} = 0$;
\sn
\item[$(\beta)$]  if $i \in S_2,\bold c \rest [I_i]^2$ has range
  $\lambda$ and is $\theta$-indecomposable;
\end{enumerate}
\sn
\item[(c)]  if $\varp \in I_i,\zeta \in I_j,i \ne j$ then 
$\bold c\{\varp,\zeta\} = 0$.
\end{enumerate}
\end{enumerate}
\mn
Why exist?  For $(*)_2(b)$ do as in \ref{b9}.  Now in $(*)_1(a)$ we
have a problem as in \cite{Sh:g}, concentrate on $\Pr_1$.  So we have
to repeat the proof of \ref{c35}.

How do we choose $\bold c \rest [I^3]^2$?  We choose $\theta
= \partial = \aleph_1,\delta \in I^3 \Rightarrow \cf(\delta) =
\aleph_1$ and each $I_{i,\alpha}(i \in S_3 \, \odd,\alpha < \lambda)$ is
stationary.  Also we choose $\bar C = \langle C_\delta:\delta \in
I^3\rangle$ such that $c_\delta \subseteq \delta = \sup(c_\delta)$ and
$\otp(C_\delta) = \omega$.  We define:
\mn
\begin{enumerate}
\item[$(*)_2$]  for $\varp < \zeta$ from $I^3,\bold c\{\varp,\zeta\}$ is: 
\begin{itemize}
\item  $\otp(C_\zeta \cap \varp)$ \when \, $(\exists i,\alpha)(\varp
  \in I_{i,\alpha} \wedge \zeta \in I_{i,\alpha} \wedge \varp \in C_\zeta)$;
\sn
\item  0 \when \, otherwise.
\end{itemize}
\end{enumerate}
\mn
We have to prove clause (d) of Definition \ref{e10}.  So we define:
\mn
\begin{enumerate}
\item[$(*)_3$]  $c \ell_1:[\lambda]^{< \aleph_0} \rightarrow
  [\lambda]^{< \partial}$ is such that:
\sn
\begin{itemize}
\item  if $u \in [\lambda]^{< \aleph_0},\delta_1 < \delta_2$ belongs
  to $I^3 \cap u$ then $C_{\delta_2} \cap \delta_2 \subseteq c \ell_1(u)$;
\sn
\item  if $\alpha < \delta \in I^3,\{\alpha,\delta\} \subseteq u$ then
$c_\delta \cap \alpha \subseteq c \ell_1(u)$.
\end{itemize}
\end{enumerate}
\mn
In our case $I_{i,\alpha}$ is a stationary subset of $\{\delta <
\lambda:\cf(\delta) = \partial\}$ so we can find $\cU_1,n_*$ such
that:
\mn
\begin{enumerate}
\item[$(*)_4$]
\begin{enumerate}
\item[(a)]  $\cU_1 \subseteq I_{i,\alpha}$ is stationary;
\sn
\item[(b)]  $v_\delta(\delta \in \cU_1)$ is constantly $v_*$;
\sn
\item[(c)]  $u_\delta \cap \delta  = u_*$ where $u_s =
  \cup\{u_{s,\iota}:\iota \in v_s\}$;
\sn
\item[(d)]  $u_{**} := \cup\{C_\alpha \cap \delta:\alpha \in
  u_\delta,\alpha > \delta\}$ does not depend on $\delta$;
\sn
\item[(e)]  $\cU_1 \subseteq E$ where $E$ is from $(*)_1(b)$;
\sn
\item[(f)]  $\delta \in \cU_1 \Rightarrow \gamma_* = \min\{\gamma \in 
C_\delta:\gamma > j_{i,\omega \alpha + n(*)}\}$ for $n(*)$ such that
$j_{i,\omega \alpha + n(*)} > \sup(v_*)$.
\end{enumerate}
\end{enumerate}
\mn
If $\lambda > 2^{\aleph_0}$, \wilog \, above $(\ell_0,\ell_1,\ldots) =
(1,\ldots)$ and choose $s_2 \in \cU_1,s_1 = \min(C_\delta \backslash
(u_* \cup u_{**} \cup \gamma_*))$, see the proof in \S5.

This is not full.  We have to repeat the intersecting with a club and
adding a set to all the $u^\ell_\varp$'s and after $\omega$ steps we are
done.  See details later.  So assume $\lambda \le 2^{\aleph_0}$.

This is helpful only if the conclusion of \ref{c29}(1) holds.  In this
case $\Pr_1(\lambda,\lambda,\aleph_0,\aleph_0)$ is enough because it
implies $\Pr_0(\lambda,\lambda,\aleph_0,\aleph_0)$; see \cite{Sh:g}
but we elaborate.

Let $\langle \eta_\alpha:\alpha < \lambda\rangle$ be a sequence of
pairwise distinct members of ${}^\omega 2$ and $\bold c_0:[\lambda]^2
\rightarrow \omega$.

Let $\langle h_n:n < \omega\rangle$ list $\{h$: for some $k,h$ is a
function from ${}^k 2 \times {}^k 2$ into $\omega\}$ exemplify
$\Pr_1(\lambda,\lambda,\aleph_0,\aleph_0)$.  We define $\bold c:[\lambda]^2
\rightarrow \omega$ by: if $\varp < \zeta < \lambda,\bold
c_0\{\varp,\zeta\} = n$ and $h_n:{}^{k(n)}2 \rightarrow \omega$ then
$\bold c\{\varp,\zeta\} = h_n(\eta_\varp \rest k(n),\eta_\zeta \rest
k(n))$.

Let $\langle u_\alpha:\alpha \in \cU_0,\cU_0 \in
[\lambda]^\lambda,\alpha \in u_\alpha \in [\lambda]^{\eta(*)}$ and
$n_* < \omega$ pairwise disjoint and $h:\eta(*) \times n(*)
\rightarrow \omega\rangle$.

For every $\alpha$ let $\langle \gamma_{\alpha,\ell}:\ell <
m_\alpha\rangle$ 
list $u_\alpha$ in increasing order.  Let
$k_\alpha(\alpha)$ be minimal such that
$h(\eta_{\alpha,\ell(\alpha)},\eta_{\alpha,\ell(2)}) =
h(\ell(1),\ell(2))$. For some $\cU_1 \in [\cU_0]^\lambda,\alpha \in
\cU_1 \Rightarrow k_\alpha = k_* \wedge \bigwedge\limits_{\ell < n_*}
\eta_{\gamma_{\alpha,\ell}}  \rest k_\alpha = \eta^*_\ell$.

Now apply the choice of $\bold c_*$.
\end{PROOF}

\begin{claim}
\label{e16}
If $\mu^{++} < \lambda,\lambda$ regular \then \,
$\Pr_{2.5}(\lambda,\theta,\partial)$ for some $\theta,\partial$.  Hence the
conclusion of \ref{c41} and \ref{c44} holds.
\end{claim}

\begin{remark}
\label{c50}
Check reference inside.
\end{remark}

\begin{PROOF}{\ref{c50}}
Let $S_\alpha \subseteq \{\delta < \lambda:\cf(\delta) \le \mu^+\}$ be
pairwise disjoint with union $S$ and choose $\bar C = \langle
C_\delta:\delta \in S\rangle$ such that ($S_\alpha$ contains a club of
$\delta$) $\Rightarrow \sup(C_\delta) = \delta,C_\delta$ a set of even
ordinals\footnote{Formally we need $\langle C_\delta \cup \{0\}:\delta
  \in S_\alpha,\alpha < \lambda\rangle$ be pairwise disjoint.  Can
  arrange but also get it by renaming so $\cup\{w_i:i \, \even\}$ will be
  the set of odd ordinals and can use $\{w_{i,\alpha}:i \, \odd\}$ by
  $S'_\beta = \{\gamma$: for some $\xi,\xi^\omega +1 + \beta = \gamma$
  or $\gamma \in S_\alpha$ and is divisible by $\xi^\omega$ for every
  $\xi < \gamma\}$.} and $C_\delta \subseteq \delta \cap S_\alpha$ and
$\otp(C_\delta ) \le \mu^+$, see \cite[Ch.III]{Sh:g} for existence.

Let $\langle h_\alpha:\alpha < \mu^+\rangle$ list the functions $h$
such that for some $u = u_h \in [\mu]^{< \aleph_0},h$ is a function
from ${}^u 2$ into $\mu^+$.

Define a $\mu$-closure operation $c \ell_*$ by:
\mn
\begin{enumerate}
\item[$(*)_2$]  for $u \in [\lambda]^{\le \mu}$ Grantica question: we have: such that
$S'_\alpha = \{\delta \in S_\alpha:C_\delta$ includes a club of
  $\delta\}$ is stationary and $\bar C$ if $\alpha < \delta$ are from
  $S \cap u$ then $C_\delta \cap \alpha \subseteq c \ell(u)$].
\end{enumerate}
\mn
Choose $W_i,W_{i,\alpha}$ as in xxx above such that:
\mn
\begin{enumerate}
\item[$(*)_3$]   $\alpha < \lambda \wedge i \, \odd \Rightarrow
W_{i,\alpha} \in \{S_\beta:\beta < \lambda\}.$
\end{enumerate}
\mn
We define $\bold c$ as in \ref{c44} above using $h:\mu^+ \rightarrow
\mu^+$, see below.  Let $\alpha < \lambda,\cU_1 \subseteq S_\alpha$ be
stationary $\langle u^1_\varp:\varp \in \cU_1\rangle,\langle
u^2_\varp:\varp \in S\rangle$ and $c \ell_*:[\lambda]^{< \aleph_0}
\rightarrow [\lambda]^{\le \mu}$ and $\xi < \sigma$ and we should find
$\varp_1,u_1,\varp_2,u_2$ as in the definition.

We may try by induction on $n$ to chose $(E_n,f_n,\bar S^*_n,f_n,\bar
v_n)$ but $(\bar S^*_n,f_n,v_n)$ are 
chosen after $(E_n,\bar S_n)$ were chosen:
\mn
\begin{enumerate}
\item[$(*)_4$]
\begin{enumerate}
\item[(a)]  $E_n$ is a club of $\lambda$ decreasing with $n$;
\sn
\item[(b)]  $E_0$ is $E$ from clause (b) of $(*)_3$, see \ref{c44};
\sn
\item[(c)]  $f_n$ is a regressive function on $S_\alpha \cap E_n$;
\sn
\item[(d)]
\begin{enumerate}
\item[$(\alpha)$]  $\bar S^* = \langle S^*_\eta:\eta \in {}^n
  \lambda\rangle$ is a partition of $S \cap E_n$;
\sn
\item[$(\beta)$]  $\bar v_n = \langle v_\eta:\eta \in {}^n \lambda\rangle$;
\end{enumerate}
\sn
\item[(e)]
\begin{enumerate}
\item[$(\alpha)$]  $S^*_\eta = \{\delta \in
  S^*_\alpha:\bigwedge\limits_{\ell < \ell g(\eta)} f_\ell(\delta) =
  \eta^{(\ell)}\} \cap E_{\ell g(\eta)}$;
\sn
\item[$(\beta)$]  $v_\eta \in [\lambda]^{\le \mu}$ is increasing with $\eta$;
\end{enumerate}
\sn
\item[(f)]  $E_{n+1} \subseteq E_n$ is such that $\eta \in {}^n
  \lambda \wedge (S^*_\eta$ not stationary) $\Rightarrow E_{n+1} \cap
  S_\eta = \emptyset$.
\end{enumerate}
\end{enumerate}
\mn
\underline{Way 1}:
\mn
\begin{enumerate}
\item[$(*)$]  Let $E_1$ be a club of $\lambda$ such that:
\sn
\begin{enumerate}
\item[(a)]  $E_1 \subseteq E$ from $(*)_3(b)$;
\sn
\item[(b)]  if $\delta \in E$ and $v \in [\delta_1]^{< \aleph_0}$, then
  $c \ell_*(v) \subseteq \delta$;
\sn
\item[(c)]  $\varp < \delta \in E \Rightarrow u^2_\varp,u^1_\varp
  \subseteq \delta$ when defined.
\end{enumerate}
\end{enumerate}
\mn
Let $\langle B_{\delta,i}:i < \mu^+\rangle$ list the closure of
$C_\delta$ in increasing order.  By the choice of $\langle
C_\delta:\delta \in S_\alpha\rangle$ there is $\varp_2 \in S_\alpha$
such that $Z = \{i < \mu^+:\beta_{\delta,i+1} \in E\}$ contains a club
of $\mu^+$.  Let $v_i$ be the $\{c \ell_*,c \ell\}$-closure of
$u^1_{\varp_2} \cup \{\beta_{\delta,j}j<i\},\langle v_i:i <
\mu^+\rangle$ is $\subseteq$-increasing continuous and $v_i \in
[\lambda]^{\le \mu}$.  So for some limit $i,h(i) = \xi$ and $v_i \cap
\varp_2 \subseteq \alpha_{\varp_2,j}$.

Let $u_j = v_i,\varp_2 = \alpha_{\varp_2,i}u_1 =$ the $c
\ell_*$-closure of $u^2_{\varp_1} \cup (u_1 \cap \varp_2)$.  We have
to check the conditions:
\mn
\begin{enumerate}
\item[$\bullet_1$]  $\varp_2 \in u_1 \in [\lambda]^{\le \mu},\varp_2
  \in u_2 \in [\lambda]^{\le \mu}$.
\end{enumerate}
\mn
[Why?  Check.]
\mn
\begin{enumerate}
\item[$\bullet_2$]  $u_1,u_2$ are $c \ell_*$-closed.
\end{enumerate}
\mn
[Why?  By their choice.]
\mn
\begin{enumerate}
\item[$\bullet_3$]  $\varp_2 \in u_2 \backslash u_1$.
\end{enumerate}
\mn
[Why?  $\varp_1 \in u_1$ obviously and $u_2 \subseteq \varp_2$ because
$\varp_2 \in E_1$.]
\mn
\begin{enumerate}
\item[$\bullet_4$]  $\varp_1 \in u_1 \backslash u_2$.
\end{enumerate}
\mn
[Why?  By the choice of $\varp_1$.]
\mn
\begin{enumerate}
\item[$\bullet_5$]  $\bold c\{\varp_1,\varp_2\} = \xi$.
\end{enumerate}
\mn
[Why?  $\varp_1 = \alpha_{\varp_2,i+1}$ and $h(i) = \xi$ by the choice
of $i$ we arrive to the main point.]
\mn
\begin{enumerate}
\item[$\bullet_6$]  if $\zeta_1 + u_1 \backslash u_2,\zeta_2 \in u_2
  \backslash u_1$ and $\{\zeta_1,\zeta_2\} \ne \{\varp_2,\varp_1\}$
  then $\bold c\{\varp_2,\zeta_2\} = 0$.
\end{enumerate}
\mn
[Why?  If $\zeta_2 = \varp_2$ and $\zeta_1 \notin S_\alpha$ - trivial
by the definition of $\bold c$.  If $\zeta_2 = \varp_2$ and recall
$v_i \cap C_\delta \subseteq \alpha_{\varp_2,i}$ and so $u^2_{\varp_1}
\subseteq \alpha_{\varp_2,i+1}$ hence $u_2 \subseteq
\alpha_{\varp_1,i} s,u_2 \backslash u_1$ is disjoint to $C_\delta$.
But $(\forall \varp < \varp_2)(\bold c\{\varp,\varp_2\} \ne 0
\Rightarrow \varp \in C_\delta$, so we are done.]

If $\zeta_2 \in S \backslash \{\varp_2\}$, then $C_{\zeta_2} \cap
\delta \subseteq c \ell$ (the $c \ell$-closure of $u_2$) $\cap \delta
= v_i = \delta$ hence is disjoint to $u_1 \backslash u_2$ to which
$\zeta_2$ belongs hence $\bold c\{\zeta_1,\zeta_2\} = 0$.

If $\zeta_2 \notin S$, recalling $\zeta_2 \ge \varp_2 > \varp_1 \wedge
\varp_2 \in E_1 \wedge E \subseteq t$, by the choice of $E$ we are done.
\end{PROOF}

\begin{discussion}
\label{e19}
(16.10.29)

\noindent
1) Addition to \ref{c44}:  For successor of singular but having added
$u_* \cup u_* \cup \{\gamma_*\}$ to ? we have to close again
$u^1_\varp,u^2_\varp$.  So change as follows:
\mn
\begin{enumerate}
\item[(a)]  $I^3 \subseteq \{\delta < \lambda:\cf(\delta)
  = \partial^+\},\partial^+ < \lambda$ (can use $\partial$ instead
  $\partial = \cf(\partial) \in (\aleph_0,\lambda)$);
\sn
\item[(b)]  each $I_{i,\alpha} \, (i \in S_3,\alpha < \lambda)$ is stationary;
\sn
\item[(c)]  $\bar C = \langle C_\delta:\delta \in I^3\rangle$.
\end{enumerate}
\mn
We use the $u^1_\varp,u^2_\varp,\gamma,\aleph_0$ times and they are
of cardinality $\le \partial$.  To formulate the $\Pr_{2.8}$ have to
add $c \ell:([I^3]^{\le \partial}) \rightarrow [I^3]^{\le \partial}$
a closure operation.

\noindent
2) For inaccessible, \wilog \, Mahlo (by \cite{Sh:365}), we have to
separate into strong limit or not.  For strong limit cardinals
see \ref{e16}.  Strong
limit, let $S = \{\kappa < \lambda:\kappa$ inaccessible not
Mahlo$\}$.   Build a sequence of regressive function and closure operations.
\end{discussion}

\begin{claim}
\label{e22}
If $\lambda$ is singular not strong limit then
$\Pr_{2.7}(\lambda,\theta,\partial)$ for $\theta = \partial =
\mu^{+7},\mu$ large enough.
\end{claim}

\begin{PROOF}{\ref{e22}}
We can find by Chernikov-Shelah \cite{CeSh:1035} a $\mu < \lambda$ and
sub-tree $\cT$ of ${}^{\mu >}\mu$ (or ${}^{\mu >}2$) with $\mu$ nodes  
and $\ge \lambda$ branches (note that $\Ded(\mu) \ge \lambda$ means
$\sup\{|\lim(\cT):\cT$ a tree with $\le \mu$ nodes$\}$; but $\lambda$
is singular hence if $\mu \in [\cf(\lambda),\lambda)$, the supremum is
obtained). 

We may demand $\mu$ is regular, though this is not essential.
Let $\langle \nu_\alpha:\alpha < \mu\rangle$ list $\cT$
 each appearing $\mu$ times.  Let $\langle \eta_\alpha:\alpha <
\lambda\rangle$ list different branches: we shall use, e.g. $\partial
= \mu^{++}$, guess clubs, i.e. let $\bar C = \langle C_\delta:\delta
\in S\rangle$ be such that $S \subseteq \{\delta
< \partial:\cf(\delta) = \cf(\mu)\},C_\delta$ a club of $\delta$ of order
type $\mu$ such that $\bar C$ guesses clubs, exists
(\cite[Ch.III,\S1]{Sh:g}).  Let $h:\mu \rightarrow
\sigma$ be such that $(\forall i < \sigma)(\forall \alpha <
\mu)(\exists^\mu \xi)(h(\xi) = i \wedge \nu_\xi = \nu_\alpha)$ for $i \,
\odd \, I_{i,\alpha} \cong \mu^{+7}$ but $\bold c \rest [I_{i,\alpha}]^2$
is such that:
\mn
\begin{itemize}
\item  if $\delta \in S,\nu_\xi \triangleleft \eta_\alpha$ and
  $\gamma$ is the $\xi$-th member of $C_\delta$ then $\bold
  c\{\pi_{i,\alpha}(\gamma),\pi_{i,\alpha}(\delta)\} = h_\alpha(\xi)$.
\end{itemize}
\mn
We can continue as in \ref{c53}.  [Alternative: use $\langle
S_\varp:\varp < \varp_*\rangle$, each for one height, $\langle
\cup\{C_\delta \cup \{\delta\}:\delta \in S_\varp\}:\varp <
\varp_*\rangle$ pairwise disjoint, in each we use branches of the same
height.]

But what about the indecomposability?  Is it enough to find $\cS
\subseteq [\lambda]^{\mu^{+7}},|\cS| = \lambda$ such that if $f:\cS
\rightarrow \theta$ then for some $i < \theta,\lambda = \cup\{u \in
\cS:f(u) < i\}$?  We need $\lambda = \lambda^{\langle
  \mu^{+7},\theta\rangle}$.  So no additional requirements.
\end{PROOF}

\begin{discussion}
\label{e25}
(16.10.29) For $\lambda$ strong limit singular can we use
\cite{Sh:E59}?  Can we use \cite[\S1,\S2]{Sh:331}?  Can we use
\cite{Sh:E81}?  We may use complicated linear orders; the point is then
given $G \in \bold K^{\lf}_\lambda$ there are only $\le \lambda$ of
them realized in it.

For a group $G \in \bold K_{\lf},\bar a = {}^\lambda G$ a linear order
$I$ and $\bold c:[I]^2 \rightarrow \lambda$ we can define $H$ such
that:
\mn
\begin{enumerate}
\item[$(*)$]  
\begin{enumerate}
\item[(a)]   $H$ is generated by $G \cup \{b_{s,\ell}:s \in I,\ell < 2\}
  \cup \{c_*\}$;
\sn
\item[(b)]   $H \subseteq \theta$;
\sn
\item[(c)]   $C_*$ has order 3 and commute with even members of $G
  \cup \{b_{s,\ell}:s \in I,\ell=1,2\}$;
\sn
\item[(d)]  if $s <_I t$ then $[b_{s,1},b_{t,2}] = a_{\bold c\{s,t\}}$;
\sn
\item[(e)]   if $s <_I t$ then $[b_{s,2},b_{s,1}) = c_*$;
\sn
\item[(f)]   $b_{s,1},b_{s,2}$ commute;
\sn
\item[(g)]   $H \in K_{\lf}$ has cardinality $|G| + |I| + \aleph_0$;
\sn
\item[(h)]   $G \le_{\gs} H$;
\sn
\item[(i)]   letting $\bar b_s = (b_{s,1},b_{s,2})$ the sequence
  $\langle \bar b_s:s \in I\rangle$ is skeleton-like inside for
  quantifier free formulas inside $H$ and hence inside any $H'$ when $H
  \subseteq_{\gS} H'$.
\end{enumerate}
\end{enumerate}
\end{discussion}

\begin{discussion}
\label{e28}
For $\lambda$ strong limit singular or really any $\lambda$ we can
look again at what we need:
\mn
\begin{enumerate}
\item[$(*)$]  for any $i,\alpha < \lambda$ and $g:[\lambda]^{<
    \aleph_0}$ such that $u \subseteq g(u) = g(g(u))$ we can find $s_1
  = s_2 \in I_2$, etc.
\end{enumerate}
\mn
\underline{Problem}:  Now in \cite{Sh:331}, \cite{Sh:E82} we have
similar/non-similar.  Here we need $w$ colours.  In fact, we get them
$\eta_2 \in T_{i,\alpha}$ of length $\omega$ and $\nu_1 = \eta \rest
(n+1),\nu_2 \in \suc(\eta \rest n) \backslash \{v_1\}$ such that the
image in $\Sigma\{I_{i_2,\alpha_2}:(i_1,\alpha_1) \ne (i,\alpha)\}$ of
$(\eta,\nu_1),(\eta,\nu_2)$ are similar.  But this is not what we need.
\end{discussion}

\begin{definition}
\label{e31}
1) Let $\Pr_{2.6}(\chi,\lambda,\mu,\theta,\partial,\sigma)$ mean
that some $(\bold c,\bar W)$ witnesses it, which means that (if $\chi =
\lambda$ we may omit $\chi$ and if $\chi = \lambda = \mu$ we may omit
$\chi$ and $\lambda$):
\mn
\begin{enumerate}
\item[$(*)$]
\begin{enumerate}
\item[(a)]  $\bar W = \langle W_i:i < \chi\rangle$ is a partition
  $\lambda$ to $\chi$ sets each of cardinality $\lambda$
\sn
\item[(b)]  $\bold c:[\lambda]^2 \rightarrow \mu$
\sn
\item[(c)]  if $i$ is even, then $\bold c \rest [W_i]^2$ is
  $\theta$-indecomposable, that is, if $\langle \cU_\varp:i <
  \theta\rangle$ is $\subseteq$-increasing with union $W_i$ \then \,
  for some $\varp < \theta$ we have $\mu = \Rang(\bold c \rest
  [\cU_\varp]^2)$
\sn
\item[(d)]  if $\boxplus_1$ then $\boxplus_2$ \when \,:
\begin{enumerate}
\item[$\boxplus_1$]
\begin{enumerate}
\item[$(\alpha)$]   $i < \lambda$ is odd;
\sn
\item[$(\beta)$]  $c \ell:[\lambda]^{< \sigma}  \rightarrow [\lambda]^{< \mu}$;
\sn
\item[$(\gamma)$]  $\varp \in u_\varp \in [\lambda]^{< \sigma}$;
\sn
\item[$(\delta)$]  $\xi < \partial$;
\end{enumerate}
\sn
\item[$\boxplus_2$]  for some $(\varp,u'_\varp,\zeta,u'_\zeta)$ 
we have:
\begin{enumerate}
\item[$(\alpha)$]  $\varp < \zeta$ are from $W_i$;
\sn
\item[$(\beta)$]  $u_\varp \subseteq u'_\varp$ and $u_\zeta \subseteq
  u'_\zeta$;
\sn
\item[$(\gamma)$]  $u'_\varp$ and $u'_\zeta$ are $c \ell$-closed;
\sn
\item[$(\delta)$]   $\varp \in u'_\varp \backslash u'_\zeta$;
\sn
\item[$(\varp)$]  $\bold c\{\varp,\zeta\} \in [\xi,\partial)$;
\sn
\item[$(\zeta)$]  if $\varp_1 \in u'_\varp \backslash u'_\zeta$ and
  $\zeta_1 \in u'_\zeta \backslash u'_\varp$ but $\{\varp_1,\zeta_1\}
  \ne \{\varp,\zeta\}$ then $\bold c\{\varp_1,\zeta_1\} = 0$.
\end{enumerate}
\end{enumerate}
\end{enumerate}
\end{enumerate}
\mn
2) We define $\Pr_{2.6.2}(\chi,\lambda,\mu,\theta,\partial,\sigma)$
similarly but replacing clause (d) by:
\mn
\begin{enumerate}
\item[$(d)'$] we have $\boxplus_1 \Rightarrow \boxplus_2$ \when \,:
\sn
\begin{enumerate}
\item[$\boxplus_1$]
\begin{enumerate}
\item[$(\alpha)$]  $i$ is odd;
\sn
\item[$(\beta)$]  $\cU_1 \in [W_i]^\lambda$;
\sn
\item[$(\gamma)$]  $\varp \in u^1_\varp \in [\lambda]^{< \sigma}$ for
  $\varp \in \cU_1$;
\sn
\item[$(\delta)$]  $\langle u^1_\varp:\varp \in \cU_1\rangle$ is a
  $\Delta$-system with heart $u_1$;
\sn
\item[$(\varp)$]  $\varp \in u^2_\varp \in [\lambda]^{< \sigma}$ for
  $\varp \in W_i$;
\sn
\item[$(\zeta)$]  $\gamma < \partial$;
\end{enumerate}
\sn
\item[$\boxplus_2$]  for some $\varp_1,\varp_2$ we have:
\begin{enumerate}
\item[$(\alpha)$]  $\varp_1 \in \cU_1,\varp_2 \in W_i$;
\sn
\item[$(\beta)$]  $\bold c\{\varp_1,\varp_2\} \ge \gamma$;
\sn
\item[$(\gamma)$]  $u^1_{\varp_1} \backslash u_1,u^2_{\varp_2}$ are
  disjoint;
\sn
\item[$(\delta)$]  if $\zeta_1 \in u^1_{\varp_1} \backslash
u_1,\zeta_2 \in u^2_{\varp_2}$ and $\{\zeta_1,\zeta_2\} \ne
\{\varp_1,\varp_2\}$ \then \, $\bold c\{\zeta_1,\zeta_2\} = 0$.
\end{enumerate}
\end{enumerate}
\end{enumerate}
\end{definition}

\begin{claim}
\label{e34}
If $\bold m \in \bold M^{2.6}_{\lambda,\theta,\bar S}$ and
$\Pr_{2.6}(\lambda,\lambda,\theta,\theta,\aleph_1,\aleph_0)$, \then \,
$G_{\bold m,\theta} \in K^{\exlf}_\lambda$ is complete,
$(\lambda,\theta,\gS)$-full and extends $G_{\bold m_1}$.
\end{claim}

\begin{remark}
\label{e37}
(16.11.13) 1) Check.

Note: we define $c \ell:[\lambda]^{\le \aleph_0} \rightarrow
[\lambda]^{< \aleph_0}$ as in \ref{c34}(2).

\noindent
2) The proof below is copied!

\noindent
3) The $I_{i,\alpha} \in \{w_j:j < \lambda$ and $i \in S_2 \Rightarrow
j$ even and $i \in S_3 \Rightarrow j \, \odd\,\}$.
\end{remark}

\begin{PROOF}{\ref{e34}}
We repeat the proof of \ref{c35} but replacing ``Case 1" in the proof
by:
\medskip

\noindent
\underline{Case 1}:  $S^\bullet$ is unbounded in $\theta$

So for $i \in S^\bullet$ choose $c_i \in \bold C(G_{\bold m,i},
G_{\bold m,i+ \omega})$ such that
$\pi(c_i) \ne c_i$ hence $c_i \notin \bold C(G_{m,i})$.  
\Wilog \, $c_i$ has order 2 because recalling $G_i \in \bold K_{exlf}$
and so the set of 
elements of order 2 from $\bold C(G_{\bold m,i},G_{\bold m,i+ \omega})$ 
generates it, see \cite[4.1=Ld36,4.10=Ld93]{Sh:312}.  Choose $\langle
i_\varp = i(\varp):\varp < \theta\rangle$ increasing, $i_\varp \in S^\bullet$
and so as $i_\varp + \omega \le i_{\varp +1} \in E$ clearly
 $\pi(c_\varp) \in G_{\bold m,i(\varp+1)}$.  Now we apply
\ref{c26}(e), \ref{w11}(1) and get contradiction by \ref{c34}(4) recalling
\ref{c26}(2)(h) and \ref{c25}(e); but we elaborate.

\noindent
Now we apply \ref{c26}(1)(e) (indirectly \ref{c29}(1), \ref{c32}).  So
there are $(i,\alpha,v,\ell_0,\ell_1,\ldots,\varp_0,\varp_1,\ldots)$ as
there in particular $i \in S_3$ and here $v = \omega \backslash
\{0\}$.  Now for every $s \in I_{i,\alpha}$
we apply \ref{c34}(2), getting $\bar u_s = \langle u_{s,\iota}:\iota
\in v_s\rangle$ and let $\ell_s$ be such that $v_s \subseteq
j_{i,\omega \alpha + \ell_s}$, \wilog \, $i \in v_s,s \in u_{s,i}$.

Now consider the statement:
\mn
\begin{enumerate}
\item[$(*)_3$]  there are $s_1 \ne s_2 \in I_{i,\alpha}$ and $k$ such
  that
\sn
\begin{enumerate}
\item[(a)]  $\bold c\{s_1,s_2\} = \ell_k$
\sn
\item[(b)]  $\ell_s > \ell_{s_1},\ell_{s_2}$
\sn
\item[(c)]  if $t_1 \in \cup\{u_{s_1,\iota}:\iota \in v_{t_1}
  \backslash i\},t_2 \in \cup\{u_{s_2,\iota}:\iota \in v_{t_2}
  \backslash i\}$ and $\{t_1,t_2\} \ne \{s_1,s_2\}$ 
\then \, $\bold c\{t_1,t_2\}=0$
\newline
or for later proofs
\sn
\item[(c)$'$]
\begin{enumerate}
\item[$(\alpha)$]  if $t_1 \in u_{s_1,i} \backslash u_{s_2,i}$ and
  $t_2 \in u_{s_2,i} \backslash u_{s_1,i}$ and
\begin{itemize}
\item  $\{t_1,t_2\} \ne \{s_1,s_2\}$ then $\bold c\{t_1,t_2\} = 0$
  \oor \, just
\sn
\item  $t_1,t_2 \in I_{i,\alpha} \Rightarrow \bold c\{t_1,t_2\} <
  \ell_k$
\sn
\item  $t_1,t_2 \in I_{i,\beta},\beta < \lambda;\beta \ne \alpha$ then
  $j_{i,\omega \beta + \bold c\{t_1,t_2\}} < j_{i,\omega \alpha +
    \ell(k)}$ (we use $j_{i,\omega \alpha + \ell} \in
  (j^*_{i,\ell},j^*_{i,\ell +1})$ - check
\end{itemize}
\sn
\item[$(\beta)$]  if $\iota \in v_1 \cap v_2$ and $\iota > i, (\iota
  \in S_3)?,\beta < \lambda$ and $t_1 \in v_{s_1,\iota},t_2 \in
  v_{s_2,\iota}$ then $\bold c\{t_1,t_2\} = 0$.
\end{enumerate}
\end{enumerate}
\end{enumerate}
\mn
Now why $(*)_3$ is true?  This is by the choice of $\bold c$, that is,
as $\bold c$ exmplifies $\Pr_0(\lambda,\lambda,\lambda,\aleph_0)$ (in
later proofs we use less).

Now to get a contradiction we like to prove
\mn
\begin{enumerate}
\item[$(*)_4$]  the type $\tp((\pi(b_{s_1}),\pi(b_{s_2}),G_{\bold m,i},
G_{\bold m,\theta})$ does not split over 
$G_{\bold m,j_{i,\omega \alpha + \ell(k)}} \cup \{c_{i(\varp_k)}\}$
hence over $G_{\bold m,i(\varp_k)} \cup \{c_{i(\varp(k))}$.
\end{enumerate}
\mn
It gives that $\tp((b_{s_1},(b_{s_2}),\pi^{-1}(G_{\bold m,i}),
\pi^{-1}(G_{\bold m,\theta})$ does not split over
$\pi^{-1}(G_{m,i(\varp_k)}) \cup \{\pi^{-1}(c_{i(\varp)})$.  But
$i(\varp_k),i \in E$ has it follows that $\pi(G_{m,i}) = G_{\bold
  m,i}$ and $\pi^{-1}(G_{i(\varp_k)} = G_{i(\varp_k)})$ has
$\tp(b_{s_1},b_{s_2}),G_{\bold m,i},G_{\bold m,\theta})$ does not
split over $G_{i(\varp_k)} \cup\{\pi^{-1}(c_{i(\varp_k)})\}$.

Now $[b_{s_1},b_{s_2}] = [b_{s_1},b_{s_2}] = \pi^{-1}(c_{i(\varp_k)})$
which is $\ne c_{i(\varp_k)}$. But as $c_{i(\varp_k)} \in \bold
C(G_{\bold m,i(\varp_k)},G_{\bold m,\theta})$ clearly also
$\pi^{-1}(c_{i(\varp_k)})$ belongs to it, hence it follows that
$\pi^{-1}(c_{i(\varp_k)}) \in c \ell(\{c_{i(\varp_k)}\};G_\theta)$,
but as $c_{i(\varp_k)}$ has order two, the latter is
$\{c_{i(\varp_k)},e_{G_\sigma}\}$.

However $\pi^{-1}(c_{i(\varp_k)})$ too has order 2 hence is equal to
$c_{i(\varp_k)}$; applying $\pi$ we get $c_{i(\varp_k)} =
\pi(c_{i(\varp_k)})$ a contradiction to the choice of the $c_i$'s.

[Pre 16.11.11 proof: as 
$c_{i(\varp_k)} \ne G_{\bold m,j_i,\omega \alpha + \ell(k)}$ (see
above), contradiction.

Why does it hold?  By \ref{c34}(3) and $(*)_3(c)$, in later proofs by
finer versions (for more cardinals) we have to do more.]
\end{PROOF}

\begin{claim}
\label{e40}
1) If $\lambda$ is regular $\ge \beth_\omega$ and $\theta <
\beth_\omega$ is regular large enough, then
$\Pr_{2.6}(\lambda,\lambda,\theta,\theta,\partial,\sigma)$.

\noindent
1A) In (1) if $\lambda$ is regular, $\lambda^{\langle
  \theta;\theta\rangle} = \lambda$ the conclusion holds.

\noindent
2) If $\lambda$ is singular not strong limit and $\theta <
\beth_\omega$ is regular large enough, then
$\Pr_{2.6}(\lambda,\lambda,\theta,\theta,\partial,\sigma)$. 
\end{claim}

\begin{PROOF}{\ref{e40}}
1) Choose $E$ such that:
\mn
\begin{enumerate}
\item[$(*)_1$]
\begin{enumerate}
\item[(a)]    $E$ is a club of $\lambda$;
\sn
\item[(b)]  
\begin{enumerate}
\item[$(\alpha)$]  if $\lambda$ is a successor cardinal, e.g. $\lambda
  = \mu^+$ then $\delta \in E \Rightarrow \mu^\omega|\delta$;
\sn
\item[$(\beta)$]  if $\lambda$ is a limit cardinal, then $E$ is a set
  of limit cardinals.
\end{enumerate}
\end{enumerate}
\end{enumerate}
\mn
Let ``$\theta$ large enough" mean:
\mn
\begin{enumerate}
\item[$(*)_2$]
\begin{enumerate}
\item[(a)]    $\theta < \beth_\omega$ is regular;
\sn
\item[(b)]  $\lambda^{[\theta;\theta]} = \lambda$ or just
  $\lambda^{\langle \lambda;\lambda\rangle} = \lambda$;
\newline
(even if $\lambda$ is a successor of a singular, its cofinality is
taken care of later).
\end{enumerate}
\sn
\item[$(*)_3$]  We choose $\bar W$ such that:
\sn
\begin{enumerate}
\item[(a)]   $\bar W = \langle W_i:i < \lambda\rangle$ is a partition of $\lambda$;
\sn
\item[(b)]  if $i$ is odd, then $W_i \subseteq \{\delta <
  \lambda:\cf(\delta) = \sigma\}$ is stationary;
\sn
\item[(c)]  each $W_i$ is unbounded in $\lambda$.
\end{enumerate}
\end{enumerate}
\mn
How do we choose the $c_i$?
\mn
\begin{enumerate}
\item[$(*)_4$]  $\bold c:[\lambda]^2 \rightarrow \theta$ is such that
\sn
\begin{enumerate}
\item[(a)]   if $\varp \in W_{i,\zeta},\zeta \in W_j$ and 
$i < j < \lambda$, then $\bold c\{\varp,\zeta\}=0$;
\sn
\item[(b)]  if $\varp < \zeta < \lambda$ are from $W_i$ and $i$ is
  even and $\suc_E(\varp) < \zeta$, i.e. $(\varp,\zeta] \cap E \ne
  \emptyset$ then $\bold c\{\varp,\zeta\} = 0$;
\sn
\item[(c)]  if $\varp < \zeta$ are successive members of $E$ then
  $\bold c \rest [\{\xi:\varp < \xi < \zeta\}]^2 \rightarrow \theta$
  is $\theta$-indecomposable;
\sn
\item[(d)]  for $i$ odd:
\begin{enumerate}
\item[$(\alpha)$]  choose $\bar C = \langle C_\delta:\delta \in
  W_i\rangle,C_\delta \subseteq \delta =
  \sup(C_\delta),\otp(C_\delta) = \sigma$ guesses clubs;
\sn
\item[$(\beta)$]  if $\varp < \zeta$ are from $W_i$ then $\bold
  c\{\varp,\zeta\} = \otp(C_\zeta \cap \varp)$ when $\varp \in
  C_\delta$ and $\bold c\{\varp,\zeta\}=0$ otherwise.
\end{enumerate}
\end{enumerate}
\end{enumerate}
\mn
[Why possible?  For clause (c) do as in \S(3A) and for clause
$(d)(\alpha)$ recall \cite[Ch.III]{Sh:g}.  Now in Definition \ref{e31},
clause (a),(b),(c) should be clear and we shall prove clause (d).  So
we are given $c \ell:[\lambda]^{< \aleph_0} \rightarrow [\lambda]^{<
  \sigma}$; (in the main case - as in \ref{c34}(2)).  We shall choose
$c \ell_1$ ``extending it" (Saharon: rest of the proof is copied).]
\mn
\begin{enumerate}
\item[$(*)_3$]  $c \ell_1:[\lambda]^{< \aleph_0} \rightarrow
  [\lambda]^{< \partial}$ is such that:
\sn
\begin{itemize}
\item  if $u \in [\lambda]^{< \aleph_0},\delta_1 < \delta_2$ belong
  to $I^3 \cap u$ then $C_{\delta_2} \cap \delta_2 \subseteq c
  \ell_1(u)$;
\sn
\item  if $\alpha < \delta \in I^3,\{\alpha,\delta\} \subseteq u$ then
$c_\delta \cap \alpha \subseteq c \ell_1(u)$.
\end{itemize}
\end{enumerate}
\mn
In our case $I_{i,\alpha}$ is a stationary subset of $\{\delta <
\lambda:\cf(\delta) = \partial\}$ so we can find $\cU_1,n_*$ such
that:
\mn
\begin{enumerate}
\item[$(*)_4$]
\begin{enumerate}
\item[(a)]  $\cU_1 \subseteq I_{i,\alpha}$ is stationary;
\sn
\item[(b)]  $v_\delta(\delta \in \cU_1)$ is constantly $v_*$;
\sn
\item[(c)]  $u_\delta \cap \delta  = u_*$ where $u_s =
  \cup\{u_{s,\iota}:\iota \in v_s\}$;
\sn
\item[(d)]  $u_{**} := \cup\{C_\alpha \cap \delta:\alpha \in
  u_\delta,\alpha > \delta\}$ does not depend on $\delta$;
\sn
\item[(e)]  $\cU_1 \subseteq E$ where $E$ is from $(*)_1(b)$;
\sn
\item[(f)]  $\delta \in \cU_1 \Rightarrow \gamma_* = \min\{\gamma \in 
C_\delta:\gamma > j_{i,\omega \alpha + n(*)}\}$ for $n(*)$ such that
$j_{i,\omega \alpha + n(*)} > \sup(v_*)$.
\end{enumerate}
\end{enumerate}
\mn
If $\lambda > 2^{\aleph_0}$, \wilog \, above $(\ell_0,\ell_1,\ldots) =
(1,\ldots)$ and choose $s_2 \in \cU_1,s_1 = \min(C_\delta \backslash
(u_* \cup u_{**} \cup \gamma_*))$, see the proof in \S5.

This is not full.  We have to repeat the intersecting with a club and
adding a set to all the $u^\ell_\varp$'s after $\omega$ step we are
done.  See details later.  So assume $\lambda \le 2^{\aleph_0}$

This is helpful only if the conclusion of \ref{c29}(1) holds.  In this
case $\Pr_1(\lambda,\lambda,\aleph_0,\aleph_0)$ is enough because it
implies $\Pr_0(\lambda,\lambda,\aleph_0,\aleph_0)$; see \cite{Sh:g}
but we elaborate.

Let $\langle \eta_\alpha:\alpha < \lambda\rangle$ be a sequence of
pairwise distinct members of ${}^\omega 2$ and $\bold c_0:[\lambda]^2
\rightarrow \omega$.

Let $\langle h_n:n < \omega\rangle$ list $\{h$: for some $k,h$ is a
function from ${}^k 2 \times {}^k 2$ into $\omega\}$ exemplify
$\Pr_1(\lambda,\lambda,\aleph_0,\aleph_0)$.  We define $\bold c:[\lambda]^2
\rightarrow \omega$ by: if $varp < \zeta < \lambda,\bold
c_0\{\varp,\zeta\} = n$ and $h_n:{}^{k(n)}2 \rightarrow \omega$ then
$\bold c\{\varp,\zeta\} = h_n(\eta_\varp \rest k(n),\eta_\zeta \rest
k(n))$.

Let $\langle u_\alpha:\alpha \in \cU_0,\cU_0 \in
[\lambda]^\lambda,\alpha \in u_\alpha \in [\lambda]^{\eta(*)}$ and
$n_* < \omega$ pairwise disjoint and $h:\eta(*) \times n(*)
\rightarrow \omega$.

For every $\alpha$ let $\langle \gamma_{\alpha,\ell}:\ell <
m_\alpha\}$ list $u_\alpha$ in increasing order.  Let
$k_\alpha(\alpha)$ be minimal such that
$h(\eta_{\alpha,\ell(\alpha)},\eta_{\alpha,\ell(2)}) =
h(\ell(1),\ell(2))$. For some $\cU_1 \in [\cU_0]^\lambda,\alpha \in
\cU_1 \Rightarrow k_\alpha = k_* \wedge \bigwedge\limits_{\ell < n_*}
\eta_{\gamma_{\alpha,\ell}}  \rest k_\alpha = \eta^*_\ell$.

Now apply the choice of $\bold c_*$.

\noindent
1A) Similarly.

\noindent
2) (16.11.13) - copied)

We can find by Chernikov-Shelah \cite{CeSh:1035} a regular $\mu <
\lambda$ and tree $\cT$ with $\mu$ nodes $\mu$-levels and $\ge
\lambda$ branches, (of finitely many heights).

Let $\langle \nu_\alpha:\alpha < \mu\rangle$ list $\cT$ such that
$\nu_\alpha \triangleleft \nu_\beta \Rightarrow \alpha < \beta$ each
appearing $\mu$ times.  Let $\langle \eta_\alpha:\alpha <
\lambda\rangle$ list different branches: we shall use, e.g. $\partial
= \mu^{++}$, guess clubs, i.e. let $\bar C = \langle C_\delta:\delta
\in S\rangle$ is such that $S \subseteq \{\delta
< \partial:\cf(\delta) = \mu\},C_\delta$ a club of $\delta$ of order
type $\mu$ such that $\bar C$ guess clubs.  Let $h:\mu \rightarrow
\sigma$ be such that $(\forall i < \sigma)(\forall \alpha <
\mu)(\exists^\mu \xi)(h(\xi) = i \wedge \nu_\xi = \nu_\alpha)$ for $i \,
\odd \, I_{i,\alpha} \cong \mu^{+7}$ but $\bold c \rest [I_{i,\alpha}]^2$
is such that:
\mn
\begin{itemize}
\item  if $\delta \in S,\nu_\xi \triangleleft \eta_\alpha$ and
  $\gamma$ is the $\xi$-th member of $C_\delta$ then $\bold
  c\{\pi_{i,\alpha}(\gamma),\pi_{i,\alpha}(\delta)\} = h_\alpha(\xi)$.
\end{itemize}
\mn
Can continue as in \ref{e10}.  [Alternative: use $\langle
S_\varp:\varp < \varp_*\rangle$, each for one height, $\langle
\cup\{C_\delta \cup \{\delta\}:\delta \in S_\varp\}:\varp <
\varp_*\rangle$ pairwise disjoint, in each we use branches of same
height.]

But what about the indecomposable?  It is enough to find $\cS
\subseteq [\lambda]^{\mu^{+7}},|\cS| = \lambda$ such that if $f:\cS
\rightarrow \theta$ then for some $i < \theta,\lambda = \cup\{u \in
\cS:f(u) < i\}$?  We need $\lambda = \lambda^{\langle
  \mu^{+7},\theta\rangle}$.  So no additional requirements.

\noindent
3) See \cite[\S2]{Sh:331} - FILL.
\end{PROOF}

\newpage

\bibliographystyle{alphacolon}
\bibliography{lista,listb,listx,listf,liste,listz}

\def\germ{\frak} \def\scr{\cal} \ifx\documentclass\undefinedcs
  \def\bf{\fam\bffam\tenbf}\def\rm{\fam0\tenrm}\fi 
  \def\defaultdefine#1#2{\expandafter\ifx\csname#1\endcsname\relax
  \expandafter\def\csname#1\endcsname{#2}\fi} \defaultdefine{Bbb}{\bf}
  \defaultdefine{frak}{\bf} \defaultdefine{=}{\B} 
  \defaultdefine{mathfrak}{\frak} \defaultdefine{mathbb}{\bf}
  \defaultdefine{mathcal}{\cal} \defaultdefine{implies}{\Rightarrow}
  \defaultdefine{beth}{BETH}\defaultdefine{cal}{\bf} \def\bbfI{{\Bbb I}}
  \def\mbox{\hbox} \def\text{\hbox} \def\om{\omega} \def\Cal#1{{\bf #1}}
  \def\pcf{pcf} \defaultdefine{cf}{cf} \defaultdefine{reals}{{\Bbb R}}
  \defaultdefine{real}{{\Bbb R}} \def\restriction{{|}} \def\club{CLUB}
  \def\w{\omega} \def\exist{\exists} \def\se{{\germ se}} \def\bb{{\bf b}}
  \def\equivalence{\equiv} \let\lt< \let\gt>
\providecommand{\bysame}{\leavevmode\hbox to3em{\hrulefill}\thinspace}
\providecommand{\MR}{\relax\ifhmode\unskip\space\fi MR }
\providecommand{\MRhref}[2]{%
  \href{http://www.ams.org/mathscinet-getitem?mr=#1}{#2}
}
\providecommand{\href}[2]{#2}
\begin{thebibliography}{}

\bibitem[Moo06]{Mo06}
Justin~Tatch Moore, \emph{{A solution to the $L$ space problem}}, Journal of
  the American Mathematical Society \textbf{19} (2006), 717--736.

\bibitem[Tod85]{To}
Stevo Todor\v{c}evi\'{c}, \emph{{Remarks on Chain Conditions in Products}},
  Compositio Math. \textbf{5} (1985), 295--302.

\bibitem[Sh:e]{Sh:e}
Saharon Shelah, \emph{{Non--structure theory}}, vol. accepted, {Oxford
  University Press}.

\bibitem[Sh:g]{Sh:g}
\bysame, \emph{{Cardinal Arithmetic}}, {Oxford Logic Guides}, vol.~29, {Oxford
  University Press}, 1994.

\bibitem[Sh:E53]{Sh:E53}
\bysame, \emph{{Introduction and Annotated Contents}}, arxiv:0903.3428.

\bibitem[Sh:E59]{Sh:E59}
\bysame, \emph{{General non-structure theory and constructing from linear
  orders}}, arxiv:new.

\bibitem[Sh:E81]{Sh:E81}
\bysame, \emph{{Bigness properties for $\kappa$-trees and linear order}}.

\bibitem[Sh:E82]{Sh:E82}
\bysame, \emph{{Bounding forcing with chain conditions for uncountable
  cardinals}}.

\bibitem[Sh:88a]{Sh:88a}
\bysame, \emph{{Appendix: on stationary sets (in ``Classification of
  nonelementary classes. II. Abstract elementary classes'')}}, Classification
  theory (Chicago, IL, 1985), Lecture Notes in Mathematics, vol. 1292,
  Springer, Berlin, 1987, Proceedings of the USA--Israel Conference on
  Classification Theory, Chicago, December 1985; ed. Baldwin, J.T.,
  pp.~483--495.

\bibitem[Sh:88r]{Sh:88r}
\bysame, \emph{{Abstract elementary classes near $\aleph_1$}}, Chapter I.
  0705.4137. arxiv:0705.4137.

\bibitem[Sh:312]{Sh:312}
\bysame, \emph{{Existentially closed locally finite groups}}, preprint
  \textbf{Beyond First Order Model Theory}, arxiv:math.LO/1102.5578v3.

\bibitem[Sh:331]{Sh:331}
\bysame, \emph{{A complicated family of members of trees with $ \omega +1 $
  levels}}, arxiv:new.

\bibitem[Sh:365]{Sh:365}
\bysame, \emph{{There are Jonsson algebras in many inaccessible cardinals}},
  Cardinal Arithmetic, Oxford Logic Guides, vol.~29, Oxford University Press,
  1994, General Editors: Dov M. Gabbay, Angus Macintyre, Dana Scott.

\bibitem[Sh:460]{Sh:460}
\bysame, \emph{{The Generalized Continuum Hypothesis revisited}}, Israel
  Journal of Mathematics \textbf{116} (2000), 285--321, arxiv:math.LO/9809200.

\bibitem[ShTh:524]{ShTh:524}
Saharon Shelah and Simon Thomas, \emph{{The Cofinality Spectrum of The Infinite
  Symmetric Group}}, Journal of Symbolic Logic \textbf{62} (1997), 902--916,
  arxiv:math.LO/9412230.

\bibitem[Sh:820]{Sh:820}
Saharon Shelah, \emph{{Universal Structures}}, Notre Dame Journal of Formal
  Logic \textbf{58} (2017), 159--177, arxiv:math.LO/0405159.

\bibitem[Sh:829]{Sh:829}
\bysame, \emph{{More on the Revised GCH and the Black Box}}, Annals of Pure and
  Applied Logic \textbf{140} (2006), 133--160, arxiv:math.LO/0406482.

\bibitem[Sh:877]{Sh:877}
\bysame, \emph{{Dependent $T$ and Existence of limit models}}, Tbilisi
  Mathematical Journal \textbf{7} (2014), 99--128, arxiv:math.LO/0609636.

\bibitem[Sh:900]{Sh:900}
\bysame, \emph{{Dependent theories and the generic pair conjecture}},
  Communications in Contemporary Mathematics \textbf{17} (2015), 1550004 (64
  pps.), arxiv:math.LO/0702292.

\bibitem[Sh:906]{Sh:906}
\bysame, \emph{{No limit model in inaccessibles}}, CRM Proceedings and Lecture
  Notes \textbf{53} (2011), 277--290, arxiv:0705.4131.

\bibitem[CeSh:1035]{CeSh:1035}
Artem Chernikov and Saharon Shelah, \emph{{On the number of Dedekind cuts and
  two-cardinal models of dependent theories}}, Journal of the Institute of
  Mathematics of Jussieu \textbf{15} (2016), 771--784, arxiv:math.LO/1308.3099.

\end{thebibliography}

\end{document}